\documentstyle[12pt]{article}
\newcommand\equ[1]{{\rm (\ref{#1})}}
\newcommand\beq[1]{ \begin{equation}\label{#1} }
\newcommand{\eeq}{ \end{equation} }

\newcommand\beqa[1]{ \begin{eqnarray} \label{#1}}
\newcommand{\eeqa}{ \end{eqnarray} }
\newcommand{\beqano}{ \begin{eqnarray*} }
\newcommand{\eeqano}{ \end{eqnarray*} }

\newtheorem{theorem}{Theorem}[section]
\newtheorem{definition}{Definition}[section]
\newtheorem{proposition}{Proposition}[section]
\newtheorem{lemma}{Lemma}[section]
\newtheorem{remark}{Remark}[section]

\newcommand\dfn[1]{ \begin{definition}\label{#1} }
\newcommand\thm[1]{ \begin{theorem}\label{#1}}
\newcommand\thmtwo[2]{ \begin{theorem}[#2]\label{#1}}
\newcommand\ethm{ \end{theorem} }
\newcommand\pro[1]{ \begin{proposition}\label{#1}}
\newcommand\protwo[2]{ \begin{proposition}[#2]\label{#1}}
\newcommand\epro{ \end{proposition} }
\newcommand\lem[1]{ \begin{lemma}\label{#1}}
\newcommand\lemtwo[2]{ \begin{lemma}[#2]\label{#1}}
\newcommand\elem{ \end{lemma} }

\newcommand{\giu}{{\medskip\noindent}}
\newcommand{\Giu}{{\bigskip\noindent}}
\newcommand{\nl}{{\smallskip\noindent}}
\newcommand{\noi}{{\noindent}}

\newcommand{\torus}{ {\bf T}   }
\renewcommand{\natural}{ {\bf N}   }
\newcommand{\real}{ {\bf R}   }
\newcommand{\integer}{ {\bf Z}   }
\newcommand{\complex}{ {\bf C}   }

\renewcommand{\a }{ {\alpha}   }
\renewcommand{\b}{ {\beta}   }
\newcommand{\g}{ {\gamma}   }

\renewcommand{\d}{ {\delta}   }
\newcommand{\D}{ {\Delta}   }
\newcommand{\e }{ {\varepsilon}   }
\newcommand{\z}{ {\zeta} }

\newcommand{\th }{ {\theta}   }

\newcommand{\k}{ {\kappa}   }

\newcommand{\x }{ {\xi}   }

\newcommand{\r}{ {\rho}   }
\newcommand{\s}{ {\sigma}   }

\renewcommand{\t}{ {\tau}   }
\newcommand{\f}{ {\varphi}   }
\newcommand{\ph}{ {\phi}   }

\renewcommand{\L}{ {\Lambda}   }

\newcommand{\cA}{ {\cal A} }

\newcommand{\cD}{ {\cal D} }
\newcommand{\cE}{ {\cal E} }
\newcommand{\cF }{ {\cal F} }
\newcommand{\cH}{ {\cal H} }
\newcommand{\cI}{ {\cal I} }

\newcommand{\cM}{ {\cal M} }

\newcommand{\cO}{ {\cal O} }
\newcommand{\cP}{ {\cal P} }
\newcommand{\cQ}{ {\cal Q} }
\newcommand{\cR}{ {\cal R} }

\newcommand{\cT}{ {\cal T} }
\newcommand{\cZ}{ {\cal Z} }

\newcommand{\ba}{ {\overline a} }
\newcommand{\br}{ {\overline r} }
\newcommand{\bs}{ {\overline s} }

\newcommand{\bw}{ {\overline w} }
\newcommand{\bz}{ {\overline z} }
\newcommand{\bZ}{ {\overline Z} }

\newcommand{\hm}{ {\hat m} }
\newcommand{\hn}{ {\hat n} }
\newcommand{\ho}{ {\hat \o} }
\newcommand{\hO}{ {\hat \O} }

\renewcommand{\o}{ {\omega}   }
\renewcommand{\O}{ {\Omega}   }

\newcommand{\ii}{ {\rm i} }

\begin{document}

\title{Periodic orbits close to elliptic
tori and applications to
the three-body problem}
\bigskip\medskip


\author{
Massimiliano Berti \\
{\footnotesize Settore di Analisi Funzionale e Applicazioni
} \\
{\footnotesize 
Scuola Internazionale Superiore di Studi Avanzati (SISSA)
} \\
{\footnotesize Via Beirut 2-4,
34014 Trieste
(Italy) }\\
{\footnotesize\tt (berti@sissa.it)}
\and
Luca Biasco \\
{\footnotesize Dipartimento di Matematica} \\
{\footnotesize Universit\`a  ``Roma Tre"} \\
{\footnotesize Largo S. L. Murialdo 1, 00146 Roma (Italy) }\\
{\footnotesize\tt (biasco@mat.uniroma3.it)}
\and
Enrico Valdinoci \\
{\footnotesize Dipartimento di Matematica }\\
{\footnotesize  Universit\`a di Roma ``Tor Vergata"}\\
{\footnotesize
Via della Ricerca Scientifica 1, 00133 Roma (Italy) }\\
{\footnotesize\tt (valdinoci@mat.uniroma2.it)}
}

\date{}
\maketitle

{\bf Abstract:}{\footnotesize
\noindent \ We prove, under suitable non-resonance and
non-degeneracy ``twist'' conditions, a Birkhoff-Lewis type
result showing the existence of infinitely
many periodic solutions, with larger and larger minimal period,
accumulating onto elliptic invariant tori
(of Hamiltonian systems). We prove the
applicability of this result to the spatial planetary three-body
problem in the small eccentricity-inclination regime.
Furthermore, we find other periodic orbits
under some restrictions on the
period and the masses of the ``planets''.
The proofs are based on 
averaging theory, KAM theory and variational methods\footnote{
Supported by M.U.R.S.T. {\sl Variational Methods and Nonlinear
Differential Equations.}}.}

\Giu

\tableofcontents


\section*{Introduction}

A classical problem in Hamiltonian dynamical systems 
concerns, since the researches of Poincar\'e,
the existence 
of periodic orbits in the vicinity of invariant submanifolds.

Poincar\'e 
wrote:  
\sl ``...voici un fait que je n'ai pu d\'emontrer rigoureusement,
mais qui me parait pourtant tr\`es vraisemblable.
\'Etant donn\'ees des \'equations de la forme d\'efinie dans
le n. 13\footnote{The Hamilton's equations.} 
et une solution particuli\`ere quelconque de ces \'equations,
one peut toujours trouver une solution p\'eriodique (dont la p\'eriode
peut, il est vrai, \^etre tr\`es longue), telle que la diff\'erence entre
les deux solutions soit aussi petite qu'on le veut, pendant un temps aussi
long qu'on le veut.'' 
\rm ([P], Tome 1, ch. III, a. 36).

This conjecture was often quoted by Birkhoff as a main 
motivation for his works. 
In the thirties, Birkhoff and Lewis [B]-[BL]-[L]
established the existence of
infinitely many periodic solutions in a neighborhood of
an elliptic equilibrium whose linear frequencies
are sufficiently {\it non-resonant}. 
This Theorem also requires a non-degeneracy
-{\it ``twist''}- condition, involving finitely many
Taylor coefficients of the Hamiltonian at the equilibrium,
and implying the system to be genuinely non-linear.
In addition, if the Hamiltonian is sufficiently smooth, KAM theory
ensures, in a neighborhood of the equilibrium small enough, the
existence of invariant Lagrangian tori filling up a set of
positive Lebesgue measure, see [P\"o].
Furthermore, close to 
any KAM torus, it has been proved in [CZ] the existence
of infinitely many others periodic orbits with
larger and larger minimal period accumulating to the torus itself
(as a consequence, the closure of the periodic orbits
has positive Lebesgue measure).
The result of [CZ] is proved considering the normal form Hamiltonian
which describes the dynamics near each  torus,
checking the ``twist condition'',
and, then,
applying the Birkhoff-Lewis type Theorem of [Mo].

Actually, Birkhoff and Lewis established also the existence
of infinitely many periodic solutions close to a non-constant
periodic elliptic solution.

The question of the existence of periodic orbits with larger and larger
minimal period clustering to
{\it elliptic lower dimensional}
invariant tori (of dimension greater than one)
has not been yet investigated.
An invariant torus is called elliptic, or
linearly stable, if the linearized system along the torus,
possesses purely imaginary eigenvalues.
\\[1mm]
\indent
A main motivation for studying such problem is Celestial Mechanics,
which, indeed, inspired the whole 
development of KAM theory
(Arnold [A] devoted one of the fundamental papers of this theory
to the planar $n$-body problem).

Consider, in particular, the
{\it non-planar planetary three body problem}, namely
one ``star'' and two ``planets''
interacting through a Newtonian gravitational field
in the three dimensional space.
The masses of the planets are regarded as small
parameters.
According to Poincar\'e and
Delaunay this system is described by a nearly-integrable
Hamiltonian on a eight dimensional phase space,
equipped with real-analytic action-angle variables.
Such system turns out to be {\it properly degenerate}, i.e. the integrable
limit (in which the three body-problem is described by two decoupled
and integrable two-body systems) depends only on two action variables.
Hence, in the integrable limit
all bounded motions lie on two-dimensional invariant
tori supporting periodic or
quasi-periodic motions according the ratio between the two frequencies
(related, by Kepler's law, to the major semi-axis
of the two limiting Keplerian ellipes) is a rational
or irrational number.
Furthemore, in the small eccentricity-inclination regime
(of astronomical interest) these unperturbed two-tori are elliptic:
the spatial planetary three body-problem calls
for a perturbation theory for continuing
{\it elliptic lower-dimensional} tori.
\\[1mm]
\indent
In the last years, an exhaustive perturbation theory for elliptic tori
has indeed been developed by many authors, see
[M], [E], [W], [K],  [P\"o1], [P\"o2], [Bo] and [XJ].
The persistence of elliptic tori is ensured
requiring ``Melnikov non-resonance conditions'' among the
frequencies and further non-degeneracy conditions.
\\[1mm]
\indent
In light of these results, the spatial planetary three
body problem (in the small
eccentricity-inclination regime) has been recently reexamined in [BCV],
where the persistence of two-dimensional elliptic
invariant tori -with diophantine frequencies-  has been proved, provided the
Keplerian major semi-axes belong to a two dimensional set of
positive measure.

Previous works on the spatial planetary three-body problem
are [JM], for large inclinations 
(in this case the two-tori are hyperbolic), and [LR]-[R], for
maximal dimensional tori  (the proofs are computer assisted).

\giu

In the present paper we first prove,
under suitable non-resonance and non-degeneracy ``twist'' conditions,
a general Birkhoff-Lewis type result showing the existence
of infinitely many periodic solutions, with larger and larger
minimal period, accumulating onto elliptic lower dimensional 
invariant tori,
see {\it Theorem \ref{theorem 1}} for the precise statement.

Furthermore we prove the applicability of Theorem \ref{theorem 1}
to the spatial planetary three-body problem, showing the existence of
infinitely many periodic solutions 
accumulating on the elliptic KAM tori of [BCV],
see {\it Theorem \ref{theorem 1 - 3 body}}.
Such periodic orbits revolve
around the star close to Keplerian ellipses with small
eccentricities and small non-zero  mutual inclinations.
The verification of the ``twist'' condition (Lemma \ref{nonde3})
is based on a KAM analysis.

Finally, in {\it Theorems \ref{theorem 2}, \ref{cor6}} we prove the
existence of other periodic solutions
of the spatial planetary three body-problem,
the ``periodic analogue'' of the elliptic tori of [BCV],
see {\it Theorem \ref{Theo}} for a more general  statement.

\giu

\noi
We now present a simplified version of our results.

\Giu

{\bf Periodic orbits accumulating on elliptic tori.}
As already said, the persistence of elliptic invariant tori is
ensured assuming ``Melnikov non-resonance conditions'' among the frequencies.
In particular, under the ``second order
Melnikov conditions'' (see precisely \equ{CONDIZ NON RISON}),
the surviving tori are still elliptic and
the normal form Hamiltonian describing the dynamics in a
neighborhood is
$$
\cH_* (\cI_*,\f_*,Z_*,\overline Z_*)=\o
\cdot \cI_* +\O Z_*\cdot \overline Z_*+ \sum_{2|k|+|a + \overline
a|\geq 3} R_{k, a, \overline a}^*(\f_*) \cI_*^k Z_*^a \overline
Z_*^{\overline a}\ ,
$$
where $ (\cI_*, \f_*) \in \real^n \times
\torus^n $ and $ ( Z_*, \overline Z_* ) \in \complex^{2m} $.
In these coordinates $ {\cal T}:=\{\cI_*=0 $, $\f_*\in\torus^n$,
$Z_*=\bZ_*= 0\}$ is the invariant and elliptic torus.
$ \o \in
\real^n $ are the torus frequencies and
$ \O \in \real^m $ the elliptic one's.

We refer to [JV] and [BHS] for many results concerning
the dynamics of an Hamiltonian system
close to an elliptic torus, in particular for the
existence, under appropiate non-resonance and non-degeneracy conditions,
of Cantor families of other invariant tori of
dimensions greater than $ n $.

Also for proving the existence of periodic orbits 
accumulating onto the elliptic torus ${\cal T}$ we need to
assume suitable non-resonance and non-degeneracy ``twist'' conditions.
Roughly (see 
Theorem \ref{theorem 1} for the detailed statement)
we have:
\begin{theorem}\label{theorem R}
Under suitable non-resonance and non-degeneracy conditions, the
Hamiltonian system generated by $ \cH_* $
affords infinitely many periodic solutions, with larger and larger
minimal period,
accumulating on the invariant elliptic torus $ { \cal T } $.
\end{theorem}

\noi
The precise 
Theorem \ref{theorem 1} and a description of its proof
are  given in Section \ref{sec-2}.
Theorem \ref{theorem 1} implies, in particular, the result of [CZ]
for maximal dimensional tori.

\Giu

{\bf The planetary, spatial three body-problem.}
Let $ \e > 0 $ denote the small parameter measuring the ratios between the
masses of the planets and the mass of the star (see \equ{masse}).
The existence of two-dimensional elliptic
invariant tori 
has been proved in [BCV]. 
Roughly (see Theorem \ref{quasip} for a precise statement and proof)
their result states:
\giu

\noi {\bf Theorem ([BCV])} {\sl The spatial planetary three-body
problem affords, for $ \e $ small en\-ough, a family ${\cal F}$ of
two-dimensional, elliptic invariant tori, travelled with
diophantine quasi-periodic frequencies, provided the osculating
Keplerian major semi-axes belong to a two-dimensional set of
density close to one (as $\e$ tends to $0$). } \giu

We will prove, in section \ref{physics}, that: 

\begin{theorem}\label{theorem 1 - 3 body}
The spatial planetary three-body problem affords, for $\e$ small enough,
infinitely many
periodic solutions, with larger and larger minimal period,
accumulating onto each elliptic invariant torus
of the family $ { \cal F}$.
\end{theorem}

The proof of Theorem \ref{theorem 1 - 3 body} boils down to check that
the non-resonance and non-degeneracy twist conditions
of Theorem \ref{theorem 1}
are fulfilled for the spatial planetary three-body problem.
This task is accomplished by estimates based on a
careful KAM analysis.
\\[2mm]
\indent
Finally, in section \ref{section theorem 2}, we prove the
existence of other periodic orbits 
of the spatial planetary three-body problem,
the ``periodic analogue'' of the elliptic tori of [BCV].
More precisely,
let consider the {\sl integrable} 
(i.e. $\e=0$) three-body
problem when 
the two planets revolve around the star
along circular orbits 
without interacting ({\sl decoupled} two-body systems).
Some of these motions will be periodic with,
say, minimal period $T$.
Taking into account the mutual interaction between the planets
(i.e. $ \e > 0 $)
we prove the existence of, at least, two,
slightly deformed, $T$-periodic orbits.
The parameter $\e $ must belong to a suitable interval of values
depending on the period $ T $.

\begin{theorem}\label{theorem 2}
There exist $T_0>0$ and  functions $0<\underline\e (T)<\overline\e
(T)$, defined for $T\in [T_0,\infty)$ and decreasing to zero
as $ T \to + \infty $
such that, for any circular periodic orbits of
the decoupled three-body problem with minimal period $ T \geq T_0 $,
for all $ \underline\e (T)\leq \e\leq \overline\e (T)$,
there exist at least two, $\e$-close,
geometrically distinct $T$-periodic orbits of the
spatial planetary three-body problem.
\end{theorem}

The following Theorem can be seen as an $ \e $-fixed version
of Theorem \ref{theorem 2},
requiring, for fixed small values of the masses of the planets, 
restrictions on the period $ T $.

\begin{theorem}\label{cor6}
There exist $ \e_1 > 0 $
and functions $ 0 < \underline{T} (\e) <\overline{T} (\e)$, defined for
$ \e \in (0,\e_1] $ and tending to infinity as $ \e \to 0 $, such
that, for all $ 0 < \e \leq \e_1 $ and for any circular periodic orbit
of the decoupled three-body problem with minimal period
$\underline{T} (\e)\leq T \leq\overline{T} (\e)$, there exist at
least two, $\e$-close, geometrically distinct $T$-periodic orbits of the
spatial planetary three-body problem.
\end{theorem}

A more general Theorem implying
Theorems \ref{theorem 2} and \ref{cor6} is given in
Theorem \ref{Theo}.

\giu

The paper is organized as follows. In Section \ref{sec-2}, we give
the detailed statement of Theorem \ref{theorem R}, namely
Theorem \ref{theorem 1}, and we discuss its proof. Section
\ref{Utilities} collects some number theoretical
Lemmata
used in Section \ref{section theorem 1}
for defining the ``non-resonant'' periods $T$.
Section \ref{Sec:NorForm}
is devoted to a normal form averaging result.
The proof of Theorem \ref{theorem 1} is addressed
in Section \ref{section theorem 1}.
In Section \ref{physics}, after recalling the
Poincar\'e-Delaunay formulation
of the spatial planetary three-body problem,
we prove, first, a finer version of the KAM result of [BCV] and, finally,
Theorem \ref{theorem 1 - 3 body}.
Section \ref{section theorem 2} contains the
proof of Theorems \ref{theorem 2}, \ref{cor6} and \ref{Theo}.

\Giu

\noindent
{{\bf Acknowledgements:} {\it We thank P. Bolle and
L. Chierchia for many interesting discussions.
Part of this paper was written when the 
last two authors were visiting SISSA in Trieste.
They thanks SISSA for its kind hospitality.}}

\Giu

{\bf Notations:} We denote by $O( \x )$ a real analytic
function whose norm is bounded by $ C \x $, for a suitable
constant $ C > 0 $ and $\forall \ 0 <\x\leq \x_0$. ${\rm Mat}(n
\times n, \real )$ (resp. ${\rm Mat}(n\times n, \complex)$) 
is the set of $ n \times n $ matrices with real (resp.
complex) coefficients and  ${\bf 1}_n$ the  $n\times n$
identity matrix. $B_r$ denotes the (closed) ball of radius $r$ centered
at $0$, $ B^n_r $ the closed ball of radius $r$ centered at $0$ using 
the $| \cdot |_2 $ norm,  
and  $D_{\rho}^d$ the complex $d$-ball 
of radius $ \rho $ centered at $0$. 
$\sharp A$ is the cardinality of the set $A$. ${\rm gcd}$ is 
the greatest common divisor and ${\rm lcm}$ the least common multiple. 
We set $e_i := (0, \ldots, 1, \ldots, 0) \in
\integer^n $ the $i$-th unit versor 
and $| a |_1 := \sum_{i=1}| a |_1 $ norm. Through the
paper $C_i, c_i$, const. will denote positive constants possibly depending
$n$, $ m $, $ \o $, $\O $, $\gamma$, $\tau$, $ r $, $ s $, $ \r $,
$ R^* $, etc.

\section{Periodic orbits accumulating on elliptic tori}\label{sec-2}

We now give the precise statement
concerning the existence of periodic orbits close to elliptic tori
(Theorem \ref{theorem 1}).
Let consider
Hamiltonians of the form 
\beq{ham (cal H)}
\cH_* (\cI_*,\f_*,Z_*,\overline Z_*)=\o \cdot \cI_* +\O Z_*\cdot
\overline Z_*+ \sum_{2|k|+|a + \overline a|\geq 3} R_{k, a,
\overline a}^*(\f_*) \cI_*^k Z_*^a \overline Z_*^{\overline a},
\eeq
where $ (\cI_*, \f_*) \in \real^n \times \torus^n $ are
action-angle variables and $ ( Z_*, \overline Z_* ) \in
\complex^{2m} $ are called the normal (or elliptic) coordinates.
The phase space $ \real^n \times \torus^n \times \complex^{m}
\times \complex^{m} $ is equipped with
the symplectic form\footnote{We denote the imaginary unit by
$\ii$ (not to be confused with $i$ often used as an index).}
$ d \cI_* \wedge d \f_* + \ii d Z_* \wedge d \overline Z_*$.

$ \o \in \real^n $ is the frequency vector and $ \O :=
$diag$(\O_1,\dots,$ $\O_m) $ is the $ m \times m $ diagonal matrix
of the normal (or elliptic) frequencies. Hence $ \O Z_* \cdot
\overline Z_* $ denotes $ \sum_{ 1\leq j\leq m} \O_j Z_{*j}
\overline Z_{*j} $. We will often identify the diagonal matrix
$\O$ with the vector $(\O_1, \ldots, \O_m) \in \real^m $. For the
multi-indexes $ k \in \natural^n, $ $ a, \ba \in \natural^m $ we
define $ | k | := \sum_{i=1}^n k_i,$ $ | a | := \sum_{j = 1}^m
a_j,$ $ | \ba |:=\sum_{j=1}^m \ba_j. $

The Hamiltonian $ \cH_* $ is assumed to be real analytic for
$(\cI_*,\f_*,Z_*,\overline Z_*)\in
D^n_{r_*}\times\torus^n_{s_*}\times D^{2m}_{\r_*}\subset
\complex^{2n+2m}$ for some positive constants $ r_ *, s_*, \r_*.$
\footnote{We have used the following standard notations: for a
given set $ A \subset \real^d $, $ d \in \natural^+ $ and $ \d, s
> 0 $, we  denote by $ A_\d^d := \{ z \in {\bf C}^d \ | \ {\rm
dist} ( z, A ) < \d \}$ and  $ {\bf T}^d_s := \{ z \in {\bf C}^d \
| \ | {\rm Im } \ z_j | < s, \ \forall\  1 \leq j \leq d \}$.
Moreover $D^d_\d:=\{ z\in\complex^d\ |\ |z|<\d \}.$
} Real analytic means that $\cH_* $ is real analytic in the real
and imaginary parts of $ Z_* $. The functions $ R_{k, a, \overline
a}^*(\f_*) $ can be expanded in Fourier series as
\begin{equation}\label{1'}
R_{k, a, \overline a}^*(\f_*)=\sum_{\ell\in \integer^n} R_{k, a,
\overline a,\ell}^* \, e^{\ii \ell\cdot \f_*}\,.
\end{equation}
\indent Since $ \cH_* ( \cI_*, \f_*, Z_*, \overline Z_*) \in \real
$ for all $ (\cI_*, \f_*) \in \real^n \times \torus^n $ and for
all $Z_*\in\complex^m$, the Taylor-Fourier coefficient $
R^*_{k,a,\ba, \ell} $ satisfy \beq{real}
{\overline{R_k^*}}_{a,\ba, \ell}= R^*_{k,\ba,a, -\ell}\,. \eeq

The frequency vector $ (\o, \O) $ is assumed to satisfy the
``second order Melnikov non-resonance conditions''
\begin{equation}\label{CONDIZ NON RISON}
|\o\cdot\ell+\O \cdot h|\geq \frac{\g}{1+|\ell |^\tau}\ ,\qquad
\forall\ \ell\in\integer^n\ , h\in\integer^m\ ,|h| \leq 2 \ ,
(\ell,h) \neq (0,0)\,
\end{equation}
for some positive constants $ \gamma, \tau $.
We will use such condition in order to perform
the averaging procedure of Section \ref{Sec:NorForm}.

By condition \equ{CONDIZ NON RISON}, the frequency $\o $ is
rationally independent\footnote{i.e. $\o \cdot \ell \neq 0, \
\forall  \ell \in \integer^n$.} (actually Diophantine), while
$(\o, \O) $ can possess some resonance relations which can be
described as follows. Eventually reordering the frequencies $\O_j
$, there exists $ 0 \leq \hm \leq m $ such that\footnote{if $\hm =
m $ we set $\ho  := \o $.} $\ho := (\o, \O_{\hm+1} , \ldots ,
\O_{m} )$ is rationally independent on $ \integer^{\hn} $, where $
\hn := n + m - \hm $, but $( \ho, \O_j )$ is rationally dependent
for all $1 \leq j \leq \hm $ ($\hm $ is the number of resonances).
This means that there exist $ M_j \in \natural^+ $ and $a_j \in
\integer^\hn $ such that\footnote{We denote by ${\rm gcd}$ and
${\rm lcm}$ the greatest common divisor and the least common
multiple, respectively.}
\beq{Mj}
M_j \O_j = a_j \cdot \ho \qquad
{\rm with} \quad  {\rm gcd}\,(a_{j1}, \ldots, a_{j\hn}, M_j)=1,
\qquad \forall 1 \leq j \leq \hm.
\eeq

The Hamiltonian system\footnote{Setting
$ Z_* = (P_* + \ii Q_* ) / \sqrt{2} $,
$\overline Z_* = ( P_* - \ii Q_*) / \sqrt{2} $ with
$(P_*, Q_*) \in \real^{2m} $ the last two equations are equivalent to
$ \dot{Q}_*= \partial_{P_*} \cH,
\dot{P}_*= - \partial_{Q_*} \cH$.}
\begin{equation}\label{Heq}
\dot{\cI}_*=-\partial_{\f_*}\cH_*\ ,\quad
\dot{\f}_*=\partial_{\cI_*}\cH_*\ ,\quad
\dot{Z}_*=\ii\partial_{\bZ_*}\cH_*\ ,\quad
\dot{\bZ}_*=-\ii\partial_{Z_*}\cH_*\
\end{equation}
possesses the elliptic invariant torus
$$ {\cal T} \,:=\, \Big\{
(\cI_*,\f_*, Z_*,\overline Z_*) \in \real^n\times\torus^n\times
\complex^{2m}\; \Big| \; \cI_* =0,\;Z_*=\overline Z_* =0 \Big\}
$$ supporting the quasi-periodic solutions $(0 , \f_{*0} + \o t  ,
0, 0) $.

We define the
symmetric ``twist'' matrix $ \cR \in {\rm Mat}(n\times n, \real)$
\begin{equation}\label{R2}
\cR_{i,i'} := (1 + \d_{i,i'}) R_{e_i+e_{i'},0,0,0}^* - \sum_{1
\leq j \leq m \atop{\ell \in \integer^n}}
\frac{1}{\o\cdot\ell+\O_j} \left(
R_{e_i,e_j,0,\ell}^*R_{e_{i'},0,e_j,-\ell}^* +
R_{e_i,0,e_j,-\ell}^*R_{e_{i'},e_j,0,\ell}^* \right)\ ,
\end{equation}
where $ R_{k,a,\overline a,\ell}^* $ are
the Fourier coefficients, introduced in
\equ{1'}, of $ R_{k,a,\overline a}^* ( \f_* ) $
 and
$\d_{i,i'} := 1 $ if $i = i' $ and $0 $ if $i\neq i' $. We also
define the matrix $ {\cal Q} \in {\rm Mat}(m \times n, \real) $ as
\beqa{Q} \cQ_{ji}  &:=&  R_{e_i,e_j,e_j,0}^* - \sum_{1\leq i'\leq
n\atop{\ell\in\integer^n}} \frac{\ell_{i'}}{ \o\cdot\ell+\O_j}
\left( R_{e_{i},e_{j},0,\ell}^* R_{e_{i'},0,e_{j},-\ell}^* +
R_{e_{i},0,e_{j},-\ell}^*
R_{e_{i'},e_{j},0,\ell}^*\right) \nonumber\\
& - & \sum_{1\leq j'\leq m\atop{\ell\in\integer^n}} \frac{1}{
\o\cdot\ell+\O_{j'}} \left( R_{0,e_{j},e_{j}+e_{j'},-\ell}^*
R_{e_{i},e_{j'},0,\ell}^* + R_{0,e_{j}+e_{j'},e_{j},\ell}^*
R_{e_{i},0,e_{j'},-\ell}^* \right)\ . \eeqa for  $1 \leq j \leq m,
\   1 \leq i \leq n.$ We stress that $\cR$ and $\cQ$ are {\em
real} matrices by \equ{real}.\label{real2}

\giu

\noi We now give the precise 
statement of Theorem \ref{theorem R}.

\begin{theorem}\label{theorem 1}
Let the Hamiltonian $ \cH_* $ in \equ{ham (cal H)} be real
analytic in $ D^n_{r_*} \times \torus^n_{s_*} \times D^{2m}_{\r_*}
$ and satisfy condition  \equ{CONDIZ NON RISON}. Assume the
``twist'' condition ${\rm det} \ \cR \neq 0$ 
and one of the following ``non-resonance'' conditions:
\begin{itemize}
\item $(a)$ one of the three cases below holds:
{\begin{itemize} \item (i) $m=1,2\qquad$ {\rm (low number of
elliptic directions)}, \item (ii) $\hm=0\qquad\ \,$ {\rm (no
resonances
among 
$(\o, \O)$)}, \item(iii) $ M_j \geq \hm\geq 1, \quad \forall 1
\leq j \leq \hm$\ ;
\end{itemize}}\item $(b)$ $\a:=\displaystyle{\min_{1\leq j\leq m}} |(\O -
\cQ \cR^{-1} \o)_j| > 0$\ ;
\end{itemize}
where $\hat m$ and $M_j$ are defined in \equ{Mj}.

Then, there exist $ \eta_0, C_0, C_1 > 0 $
\footnote{Such constants may depend on $ n
$, $ m $, $ \o $, $ \O $, $ \gamma$, $ \tau $, $r_*$, $ s_* $,
$\r_*$, $ R^*_{k,a,\ba, \ell} $,
$\hat m$, $\cQ$, $\a$.}
such that: $ \forall \eta \in (0, \eta_0 ] $, there exist an open
set of periods $ \cA_\eta\subseteq
[\frac{1}{\eta^2},\frac{1}{\eta^2}+C_0]$ with measure greater than
$1/C_0$ such that $\forall T\in\cA_\eta$ there exist
$k:=k(T)\in\integer^n$, $\widetilde\o:=\widetilde\o(T)\in\real^n$,
with $ \widetilde\o T = 2 \pi k $,  $ | \widetilde\o - \o | \leq C_1
\eta^2 $, and at least $ n $ geometrically distinct $T$-periodic
solutions $ \z_\eta (t) =$ $(\cI_{*\eta}(t),$
$\f_{*\eta} (t),$ $Z_{*\eta}(t),$
$\overline Z_{*\eta}(t)) $ of the Hamiltonian system
\equ{Heq} satisfying
\begin{itemize}
\item  $(i)$ $\sup_{t\in\real} \Big( |\cI_{*\eta} (t)|+
|Z_{*\eta}(t)| + | \overline Z_{*\eta}(t)|\Big)
\leq C_1 \eta^2 $ ; \item
$(ii)$ $\sup_{t\in \real} | \f_{*\eta} (t) - ( \f_{*\eta} (0) +
\widetilde \o t) | \leq C_1 \eta$.
\end{itemize}
In particular, the closure of the family of periodic orbits
$\z_\eta$ of \equ{Heq} contains the elliptic torus ${\cal T}$.
Moreover the minimal period $T_{min}$ of $\z_\eta $ satisfies
$ T_{min} \geq T^{1/(\tau+1)}/C_1 = O( \eta^{-2 / (\tau+1)}  ).$

\end{theorem}

\giu

{\bf Idea of the proof of Theorem \ref{theorem 1}.}
Since $ ( \o, \O ) $ satisfy the
second order Melnikov non-resonance conditions
\equ{CONDIZ NON RISON},
in view of an averaging procedure (Proposition \ref{averaging}),
the normal form Hamiltonian $ \cH_* $ is,
in a suitable set of coordinates
$(I, \phi, z, \bz ) \in \real^n \times \torus^n \times \complex^{2m}$,
and sufficiently close to the torus $ { \cal T }$,
a small perturbation of the integrable Hamiltonian
$$
H_{int} := \o \cdot I + \frac{\eta^2}{2} \cR I \cdot I +
\O z\bz  + \eta^2 {\cal Q} I \cdot z \bz,
$$
where 
$ \eta > 0 $ is a small rescaling parameter measuring the
distance from ${\cal T}$. 
The Hamiltonian system generated by $ H_{int} $
possesses the elliptic tori $ {\cal T}(I_0) := \{ I = I_0, \ \phi
\in \torus^n , \ z =  
0 \} $
supporting the linear flow $ t \to \{ I_0 , \ \phi_0 + (\o +
\eta^2 \cR I_0) t, 0  
\} $. On the normal space the dynamics is
described by $ \dot{z} = \ii \O_\eta (I_0) z $ 
where $\O_\eta (I_0) :=$  $\O +\eta^2
{\cal Q} I_0 $ is the vector of the {\it ``shifted elliptic
frequencies''}.
Our task is to find periodic solutions
bifurcating from the ones of $ H_{int} $.

The system is properly nonlinear 
by the ``twist condition'' det $ \cR \neq 0 $
(see remarks \ref{twist} and \ref{intera} for a comment).
Such condition
involves only finitely many Taylor coefficients $ R_{k, a,  \ba}^* $
of the normal form Hamiltonian $ {\cH}_* $
and ensures, in particular, 
that the frequencies $ \o + \eta^2 \cR I_0 $ vary with the actions
$ I_0 $.
When
$$
\widetilde \o :=
\o+\eta^2 \cR I_0 \in \frac{1}{T} 2 \pi \integer^n,
$$
$ {\cal T}(I_0) $ is a {\it completely resonant} torus,
supporting the
family of $T$-periodic motions $ {\cal P} :=$ $ \{ I(t) = I_0,$
$\phi (t) = \phi_0 + \widetilde \o t$, $z(t) = 
0 \}$
(completely resonant frequencies 
$ \widetilde \o \in (2 \pi / T) \integer^n $ 
always exists for some $ T = O ( \eta^{-2}) $ and $I_0 = O(1)$,
see \equ{I0}, \equ{kT}). 
The family $ { \cal P} $, diffeomorphic to $ \torus^n $,
will not persist in its
entirety for the complete Hamiltonian system due to resonances
among the oscillations.

The {\it key point} to continue some periodic solutions of the
family $ { \cal P} $ 
is to choose properly the ``1-dimensional parameter'' $ T $
(the period) and the actions $ I_0 $:
the period $ T $ and the
``shifted elliptic frequencies'' $ \O_\eta (I_0) $ must satisfy
a suitable {\it non-resonance property}, see \equ{Mmeno1}.

In Lemma \ref{a} (resp. \ref{b}) we actually show the existence,
assuming condition $(a)$ (resp. $(b)$) in Theorem \ref{theorem 1},
of ``non resonant'' periods $T$. 
Note that in condition $(a)$
the ``coupling matrix'' $ \cQ $
does not play any role
(in particular for $ m = 2 $ which is the case of
the spatial three body problem).
If conditions $(a)$-$(ii)$-$(iii)$ are not verified,
the system possesses ``lower order'' resonances
(see equation \equ{Mj}), which are an obstruction
for finding ``non-resonant'' periods $T$.
In this case one needs, in order to ``move away'' from the resonances,
to evaluate the matrix $ \cQ $ and check
condition ($b$).
Conditions $(a)$-$(b)$ are, indeed, sharp:
if violated,
it is not possible, in general, to find any ``non-resonant''
period $ T $,
see remark \ref{linint}.

After the previous construction,
the proof of Theorem \ref{theorem 1} is based on a
Lyapunov-Schmidt reduction, variational in nature, inspired to
[ACE]-[AB]. The non-resonance property \equ{Mmeno1}
and the ``twist condition'' det$\cR \neq 0$,
allow to build,
by means of the Contraction Mapping Theorem,
a suitable family of pseudo-$T$-periodic solutions
branching off the family ${\cal P}$,
see Lemma \ref{pseudo (1)}.
Finally, by a variational argument, we will select at least
$n$ ($=$ cat $\torus^{n-1}$)
geometrically distinct $T$-periodic solutions from them.
\giu


Finally, we recall that, for strictly convex nearly integrable Hamiltonian
systems, [BK] have proved, via a variational
argument, the existence of periodic orbits with arbitrarily large
minimal periods. We have not tried to extend their approach
in the present case, since the unperturbed Hamiltonian $ H_{int} $
is not-strictly convex (it is linear in the elliptic
action variables $z \bz $).
Moreover the previously described
phenomenon of resonances between the period $T$ and the
``shifted elliptic frequencies'' is hardly recognizable
by a purely variational approach
(which works well for finding periodic orbits
near maximal dimensional tori and problably near
hyperbolic one's). 

\giu

\section{Study of the resonances and technical Lemmata}\label{Utilities}

In this section we first collect some number theoretical
Lemmata that will be used in section \ref{section theorem 1}
for defining the ``non-resonant periods'' $ T $
(Lemma \ref{tau} will be used in Lemma \ref{a}).

\lem{congruenze}
Let $ a_1,\dots,a_n \in \integer,$ $ M \in \natural $ with
{\rm {\rm gcd}\,}$ ( a_1, \dots, a_n, M) = 1.$ Then,
$\forall\, b\in\integer$, the
congruence
$a_1 k_1+\dots +a_n k_n\equiv b$ \ ${\rm mod} M $ has $M^{n-1}$ solutions.
\elem
\noindent {\bf Proof:}
Recall a well known result from number theory:
let $a,b,M\in\integer$ and $c:={\rm {\rm gcd}\,}(a,M)$;
then the congruence $ak\equiv b$ ${\rm mod} M$, $k\in\integer$, has
solution
if and only if\footnote{
$c|b$ means that there exists $n\in\integer$ such that $b=cn.$}
$c|b$ and in this case it has  exactly
$c$ solutions. It means that there exist integers
$ 0 \leq k_1 < \dots < k_c \leq M - 1 $
such that $ a k_h - b \in M \integer $ $\forall 1\leq h\leq c,$
but  $ak-b\notin M\integer$ $\forall\, k\in\{ 0,\dots,M-1 \}$, $k\neq k_h$
 $\forall 1\leq h\leq c.$

We now prove the Lemma by induction over $n.$
The case $n=1$  is the above mentioned result.
We now suppose the statement true for $n$ and prove it for $n+1.$
Let $c:={\rm {\rm gcd}\,}(a_{n+1},M).$ Our congruence is equivalent to
$a_{n+1}k_{n+1}\equiv (b-\sum_{i=1}^n a_i k_i)$ ${\rm mod} M $
which has exactly $c$ solutions if and only if
$c|(b-\sum_{i=1}^n a_i
k_i)$.
Hence
we have only to prove that the number of integers
$0\leq k_1,\dots ,k_n\leq M-1$
for which $c|(b-\sum_{i=1}^n a_i k_i)$ is exactly $M^n/c=c^{n-1}(M/c)^n$.
This amounts to prove that
the number of integer solutions of
$ \sum_{i=1}^n a_i k_i\equiv b\, {\rm mod} c$, is exactly $c^{n-1}.$
This is true, by the inductive hypotesis,
if ${\rm {\rm gcd}\,}(a_1,\dots,a_n,c)=1.$
The last equality is actually true since, by hypotesis,
${\rm {\rm gcd}\,}(a_1,\dots,a_{n+1},M)=1;$ hence
if $d|a_1,\dots,a_n,c$, since $c|a_{n+1},M$, then $d|a_{n+1},M$
and it follows $d=1.$
\hfill$\Box$

\lem{ho}
Let $\hm,\hn\in\natural^+.$
Let $\ho\in\real^\hn$ rationally independent. Let $\hO\in\real^\hm$
with $M_j\hO_j=a_j\cdot\ho$ for $M_j\in\natural$,
$a_j\in\integer^\hn$.
Suppose that {\rm {\rm gcd}\,}$(a_{j1},\dots,a_{j\hn},M_j)=1$.
Let $M:={\rm {\rm lcm}\,}(M_1,\dots,M_\hm)$
and $K:=\{ 0,\dots,M-1 \}^\hn.$
Suppose $M_j\geq \hm,$
$\forall\, 1\leq j\leq \hm$.
Moreover, in the case $\hm=1$, suppose also $M_1\geq 2.$
Then there exists $k\in K$
such that $a_j\cdot k/M_j\notin\integer$, $\forall\, 1\leq j\leq \hm$.
\elem
\noindent {\bf Proof:}
Let $\L_j:=\{ k\in K\ {\rm s.t.} \
a_j\cdot k\in M_j\integer \}$,
$\L_j^*:=\{ k\in K_j\ {\rm s.t.} \
a_j\cdot k\in M_j\integer \}$, where
$K_j:=\{ 0,\dots,M_j-1 \}^\hn.$
Since $M_j | M$ we have\footnote{We denote by $\sharp A$
the cardinality of the set $A$.}
$\sharp \L_j= (M/M_j)^\hn\sharp\L_j^*.$
Using Lemma \ref{congruenze}
with the substitutions
$ n \to \hn, $ $b\to 0$, $a_i\to a_{ji},$ $M\to M_j$
we obtain that $\sharp \L_j^*=M_j^{\hn-1}.$
Hence we have $\sharp \L_j=M^\hn/M_j.$
Observing that $\{ 0\}\in\L_j$ $\forall\, 1\leq j\leq \hm,$
we have
\beq{numeri}
\sharp\bigcup_{j=1}^\hm (\L_j\setminus\{ 0\}) \leq
\sum_{j=1}^\hm\sharp (\L_j\setminus\{ 0\}) =
\sum_{j=1}^\hm \left( \frac{M^\hn}{M_j}-1\right)=
-\hm + M^\hn \sum_{j=1}^\hm\frac{1}{M_j}.
\eeq
If $ \hm \geq 2 $, since $ M_j \geq \hm $ by hypothesis, then
$ -\hm + M^\hn \sum_{j=1}^\hm\frac{1}{M_j}
\leq -\hm + M^\hn <  -1 + M^\hn $. If  $ \hm = 1 $, since
$ M_1 = M \geq 2$ by hypothesis,
$ -\hm + M^\hn \sum_{j=1}^\hm \frac{1}{M_j} =
- 1 + M^{\hn -1} <  -1 + M^\hn $. In both cases, from
\equ{numeri}, we get
$$
\sharp\bigcup_{j=1}^\hm (\L_j\setminus\{ 0\}) <  M^\hn - 1 =
\sharp (K\setminus\{ 0\})
$$
and the conclusion of the Lemma follows
from $\sharp\cup_{j=1}^\hm \L_j<\sharp K$.
\hfill$\Box$

\giu

In view of the next Lemma we recall the definition of the
ergodization time of a torus with linear flow.
For any vector $ \xi \in {\bf R}^\hn $ the ergodization time
$ T_{\rm {erg}}( \xi, \e) $
required to fill $ {\bf T}^\hn $ within $ \e > 0 $
is defined as
\beq{ergtime}
T_{\rm {erg}}( \xi, \e) := \inf \Big\{ t\in {\bf R}_+ \ \Big| \
\forall x \in {\bf R}^\hn,\ \forall\, 1\leq i\leq\hn,
\  {\rm dist}(x_i, A_i+[0,t] \xi_i + {\bf Z} ) \leq \e  \Big\}
\eeq
where  $A$ is some point of
${\bf R}^\hn$.
$ T_{\rm {erg}}( \xi, \e) $ is clearly independent of the choice of
$ A $.

If $ ( \o, \O ) $ are rationally independent or are
resonant only at a sufficiently high order,
namely if $ ( \o, \O)  \cdot p \neq 0 $, $ \forall p \in {\bf Z}^{n+m} $
with $ 0 < |p| \leq  a \slash \delta $ for some
constant $ a:= a_{n,m} $, then
the trajectories of the linear flow
$ \{ (\o, \Omega) t  \}_{t \in {\bf R}} $
will make an arbitrarly fine $\delta$-net of $ \torus^{n+m} $,
see Theorem 4.1 of [BBB]. This non-resonance
assumption on $(\o, \O )$ is clearly sufficient
to prove \equ{disto}-\equ{distO} below.
In the present case, however,
the weaker non-resonance assumptions \equ{2}, \equ{tulkas}
are sufficient.

\lem{tau}
Let $ M_j, a_j, \hm $ as in \equ{Mj} and suppose that
\beq{2}
M_j \geq \hm \qquad \forall\, 1 \leq j \leq \hm\ .
\eeq
In the case $\hm=1=M_1$  suppose also
\beq{tulkas}
\exists\ n+1\leq i_*\leq \hn=n+m-1\, \qquad {\rm such\ \ \ that}\qquad
a_{1 i_*}\neq 0\ .
\eeq
$\forall 0<\d\leq 1/ (2\b )$, where $\b:=2 {\rm max}_{1\leq j\leq \hm}
|a_j|_1$, and $\forall\, t_0>0,$ there exists $\tau\geq t_0$ such that
\beqa{disto}
{\rm dist} (\o_i\tau,\integer) &\leq& \d\qquad \forall\, 1\leq i\leq n\ ,\\
{\rm dist} (\O_j\tau,\integer) &\geq&
\min \left\{
\frac{1}{2\b}\ ,\ \frac{1}{4\max_{1\leq j\leq \hm} M_j}
\right\}=: d_0
\qquad \forall\, 1\leq j\leq m\ .
\label{distO}
\eeqa
Moreover  $ \tau - t_0 \leq T_{\rm {erg}} ( \ho/M, \d/M ) $
where $ M:={\rm {\rm lcm}\,}(M_1, \ldots , M_\hm ) $.

\noindent
If $ ( \o, \O) $ are rationally independent, namely $ \hm = 0 $,
{\rm \equ{disto}--\equ{distO}} are still verified
setting $\b:=2,$ $d_0:=1/4$, $M:=1.$

\elem
\noindent {\bf Proof:}
We first consider the cases $\hm\geq 2$ or $\hm=1<M_1.$
We are in the hypoteses of
Lemma \ref{ho}. Hence  there exists
$k\in K:=\{ 0,\dots,M-1 \}^\hn$
such that $a_j\cdot k/M_j\notin\integer$, $\forall\, 1\leq j\leq \hm$.
Let\footnote{We define
$\langle (y_1,\dots,y_n) \rangle:=(\langle y_1\rangle,\dots,
\langle y_n\rangle)$ where the function
$\langle \cdot \rangle : \real \to [-1/2,1/2) $
is defined as $\langle y \rangle:=y$ for $y\in[-1/2,1/2)$
and it is $1$-periodically extended for $y\in\real.$ Note that,
for $ y \neq (1 \slash 2) + \integer $,
$ \langle y \rangle = y - [y] $ where $[y] $ is the closest
integer to $ y $.}
$$x:= \Big\langle
\frac{1}{M} \Big( k + \frac{1}{\b}\sum_{i=n+1}^\hn e_i \Big) -
\frac{\ho}{M}t_0 \Big\rangle \in [-1/2,1/2)^\hn\ .
$$

Since $ \ho / M $ is rationally independent, its linear
flow ergodizes the torus $ \real^\hn/\integer^\hn $.
Let $ t_0 \leq $ $\tau \leq $
$ t_0 +  T_{\rm {erg}} ( \ho/M, \d/M ) $
the first instant for which
$\exists h\in\integer^\hn$ with $|\frac{\ho}{M}(\tau-t_0)-x-h|_{\infty}
\leq $ $\d/M.$
If $y:=\ho\tau-k-\frac{1}{\b}\sum_{i=n+1}^\hn e_i-Mh'$ then,
by a suitable choice of $h'\in \integer^\hn$, we obtain by the
construction above that
$|y|_{\infty}\leq \d.$
Hence, by definition of $ \hat \o $,
$\forall\, 1 \leq i\leq n,$
$\o_i\tau = \ho_i \tau =
y_i+k_i-Mh_i'$ and \equ{disto} holds.\giu

\noindent
Moreover,
from the resonance relation \equ{Mj}, $\forall\,
1\leq j\leq \hm,$
\beq{OJt}
\O_j\tau=\frac{a_j\cdot\ho\tau}{M_j}=
\frac{a_j\cdot y}{M_j}+\frac{a_j\cdot k}{M_j}+
\frac{1}{\b M_j}\sum_{i=n+1}^\hn a_{ji}+\frac{M}{M_j}a_j\cdot h'\ .
\eeq
We observe that  $a_j\cdot k/M_j\notin\integer$ implies
${\rm dist}(a_j\cdot k/M_j,\integer)\geq 1/M_j$.
Also,
$|a_j\cdot y|\leq
|a_j|_1\d$,
$|\sum_{i=n+1}^\hn a_{ji}|\leq |a_j|_1$.
Hence, collecting these observations and \equ{OJt}, we have
\beq{d1}
{\rm dist} (\O_j\tau,\integer)\geq
\frac{1}{M_j}-\frac{|a_j|_1\d}{M_j}-\frac{|a_j|_1}{\b M_j}
\geq \frac{1}{4 M_j}\qquad \forall\, 1\leq j\leq \hm,
\eeq
recalling also that $ 0 < \d \leq (1/2 \beta )$ and the definition
of $ \beta $.
Finally, by definition of $ \hat \o $,
for $ \hm + 1 \leq j \leq m $, we have
$\O_j \tau = $ $\ho_{n+j-\hm} \tau =$
$y_{n+j-\hm}+$ $k_{n+j-\hm}+1/\b + $ $Mh'_{n+j-\hm}$, which implies
\beq{d2}
{\rm dist} (\O_j\tau,\integer)\geq \frac{1}{\b}-\d\geq \frac{1}{2\b}\qquad
\forall\, \hm+1\leq j\leq m\ ,
\eeq
and \equ{distO} follows from \equ{d1}-\equ{d2}.

We now consider the case $\hm=1$, $M_1=1. $ We have $ M = M_1 = 1 $.
We set $a:=a_1$ to simplify the notation.
Let
$$
x:= \Big\langle \sum_{i=n+1}^\hn b_i e_i-\ho t_0 \Big\rangle\ ,
$$
where
$$
b_i := \left\{
\matrix{\displaystyle{
\frac{{\rm sign}(a_i)}{2\sum_{i=n+1}^\hn |a_i| }}
& \qquad{\rm if}\quad a_i\neq 0\ ,\cr
\displaystyle{\frac{1}{2|a|_1}} &  \qquad{\rm if}\quad a_i = 0\ .
}
\right.
$$
Note that  $| b_i | \geq 1 \slash ( 2|a|_1) $ for all $ i $.
Proceeding exactly as above
we find $ t_0 \leq \tau \leq t_0 + T_{\rm {erg}} ( \ho, \d) $
and an $ h'\in\integer^\hn$ such that, defining
$$
y:= \ho \tau - \sum_{i=n+1}^\hn  b_i e_i -h',
$$
we have
$|y|_\infty\leq \d$ and
\equ{disto} holds.
The proof of \equ{distO} is slightly different.
For $j=1$ we have
$$
\O_1\tau = a\cdot\ho \tau = a\cdot y+
a\cdot
\sum_{i=n+1}^\hn
b_i
e_i
+ a\cdot h'
=
a\cdot y + \frac{1}{2} +a\cdot h'
\ ,
$$
where, in the last equality we have used \equ{tulkas}
abd the definition of $ b_i $.
Since $ a \cdot h' \in \integer $ and
$ |a \cdot y | \leq | a |_1 \d \leq 1/4 $ we have
\beq{manwe}
{\rm dist} (\O_1\tau,\integer)\geq
\frac{1}{2}-\frac{1}{4}=\frac{1}{4}\ .
\eeq
Instead, for $2\leq j\leq m$ we have
$
\O_j\tau = \ho_{n+j-1} \tau  = y_{n+j-1}+  b_{n+j-1}
+h'_{n+j-1}
$
which implies
\beq{ulmo}
{\rm dist} (\O_j\tau,\integer)\geq
b_{n+j-1} -\d\geq
\frac{1}{2|a|_1}-\d=
\frac{1}{\b} -\d\geq \frac{1}{2\b}
\qquad
\forall\, 2\leq j\leq m\ ,
\eeq
and \equ{distO} follows from \equ{manwe}-\equ{ulmo}.

We finally consider the case
$ \hm = 0 $,  which is the simplest one
since the linear flow of $\ho=(\o,\O)$ ergodizes the whole
torus $\torus^{n+m}.$ Let define
$$
x:= \Big\langle
\frac{1}{2}\sum_{i=n+1}^{n+m} e_i-\ho t_0
\Big\rangle.
$$
There exists
$ t_0 <  \tau < t_0 + T_{\rm erg}(\ho, \d) $ and
$ h \in \integer^{n+m} $ such that
$ | \ho (\tau - t_0 ) - x - h |_\infty \leq \d $,
with $ 0 < \d \leq 1/4$. Arguing as before
we get $ | \ho \tau - (1 \slash 2) \sum_{i=n+1}^{n+m} e_i - h' |_{\infty}
\leq \d $ for a suitable $h' \in \integer^{n+m} $ and the estimates
\equ{disto}-\equ{distO} follow. The Lemma is proved.
\hfill$\Box$

\giu

The next two Lemmata, whose proof is omitted, will be used
in the proofs of Theorems \ref{theorem 1} and
\ref{theorem 1 - 3 body}, for constructing the
pseudo periodic solutions through
 the Contraction Mapping Theorem.

\lem{lemma 0}
Let $(X, | \cdot |_X) $ and $(Y, | \cdot |_Y)$
be Banach spaces, $ L: Y \longrightarrow X$
be a linear bounded operator and $P:
X \longrightarrow Y $ be a $C^1$ map. Assume that
\beq{d0}
 \d_0 \geq 2 \, |L(P(0))|_X \,
\eeq
and
\beq{DP}
 \sup_{x\in B_{\d_0}} |DP (x)| \leq \frac{1}{2 \, |L|}\,.
\eeq
Then, there exists a unique $x^\star\in B_{\d_0}$ such
that $x^\star = L(P(x^\star) )$.
\elem

\noindent The next Lemma defines a suitable
``Green operator'' $ L $
associated to the linear system \equ{eq **} below.

\lem{lemma 1}
Let $ T > 0 $, $\O\in {\rm Mat}(m\times m,\,\real)$,
$M \in {\rm Mat}(n\times n,\,\real)$ and define
$\cM:= {\bf 1}_m - e^{\ii\O T}
\in {\rm Mat}(m\times m,\,\complex)$.
Assume that $ M $ and $ \cM $ are invertible.
Let
\begin{eqnarray}
Y &:=& C \Big( [0,T], \, \real^n \times \torus^n \times
\complex^{m}\Big) \,, \label{spaceY} \\
X &:=& \Big\{ (J,\psi,z)\in Y
\,\; {\rm s.t.}\,\; \psi(0) = \psi(T) = 0 \, , z(0)= z(T) \,
\Big\}. \label{spazY}
\end{eqnarray}
$ X $ and $ Y $, endowed with the sup-norm
$ \| (\hat J,\hat \psi,\hat z)  \|:=$
$\sup_{t \in [0,T] } ( |\hat J (t)|,$ $|\hat \psi (t)|,$ $|\hat z (t)| ),$
are Banach spaces.

For any $(\hat J,\hat\psi, \hat z)\in Y $ define the constants
\begin{eqnarray*}
\a(\hat J,\hat \psi)&:=& -\frac{1}{T} \Big(
\int_0^T \int_0^s \hat J(\th)\,d\th\,ds + M^{-1}
\int_0^T \hat \psi(\th)\, d\th
\Big)\,\in \real^n\,,\\
\b(\hat z) &:=& \cM^{-1} e^{\ii\O T}
\,\int_0^T e^{-\ii \O \th} \hat z (\th)\,
d\th\,
\in\complex^m\,
\end{eqnarray*}
and the linear ``Green'' operator $ L : Y \longrightarrow X $ by
\beq{LGreen}
 L\left(
\matrix{
\hat J \cr
\hat \psi \cr
\hat z
}\right) \;:=\;
\left(\matrix{
\a(\hat J,\hat\psi) + \int_0^t \hat J(s)\,ds\cr
\cr
M\a t+M\,\int_0^t \int_0^s \hat J(\th)\, d\th\,ds+
\int_0^t \hat \psi(\th)\,d\th
\cr
\cr
e^{\ii\O t}\Big( \b (\hat z)
+\int_0^t e^{-\ii \O \th} \hat
z(\th)\,d\th\Big)
}\right).
\eeq
The Green operator $ L $ satisfies
\begin{equation}\label{eq *}
|L|\leq C \,\Big( |M^{-1}|+
|M|\,T^2+ |M|\,|M^{-1}|\,T + |\cM^{-1}| T \Big)
\, ,
\end{equation}
for a suitable constant $C > 0$ depending on
$\sup_{t\in [0,T]}|e^{\ii\O t}|$, $n$ and $m$.

Setting $(J,\psi,z) = L(\hat J,\hat \psi,\hat z) $,
then $(J,\psi,z) \in C^1 $ and
$$
\left\{ \matrix{
\dot J &=& \hat J, \nonumber \cr
\dot \psi - M\,J & = & \hat\psi, \nonumber \cr
\dot z -\ii\O z &=& \hat z\label{eq **}
\,.
} \right.
$$
\elem

\noindent The straightforward proof is omitted. We only point out that
$ L(Y) \subseteq X $ by definition of the constants
$\a(\hat J,\hat\psi)$ and $\b(\hat z)$.\giu

\section{Normal form around an elliptic torus}\label{Sec:NorForm}

In order to study the dynamics of the Hamiltonian system \equ{Heq}
in a small neighborhood of ${\cal T}$ it is a convenient device to
perform the following
rescaling
\begin{equation}\label{resca}
\cI_*:=\eta^2 \cI\,,\qquad \f_*:= \f\,,\qquad Z_*:=\eta Z\,,\qquad
\overline Z_*:=\eta \overline Z,
\end{equation}
where $ \eta > 0 $ is a positive small parameter.
Note that the linear tranformation \equ{resca} preserves the torus
$ \cT $ and that a domain of order one in the new variables
$(\cI,\f,Z,\overline Z)$ correspond to a domain in the old variables
$(\cI_*,\f_*,Z_*,\overline Z_*)$ that shrinks towards $\cT$
for $\eta $ tending to zero.

The new Hamiltonian $ \cH (\cI_,\f,Z,\overline Z) =
\eta^{-2} \cH_* ( \eta^2 \cI, \f, \eta Z, \eta \overline Z) $ writes
\begin{equation}\label{ham (2)}
\cH (\cI_,\f,Z,\overline Z)=\o \cdot \cI +\O Z\cdot \overline Z+
\sum_{2|k|+|a + \overline a|\geq 3} \eta^{2|k|+|a + \overline a|-2}
R_{k, a, \overline a}^*(\f)
\cI^k Z^a \overline Z^{\overline a}\,
\end{equation}
and it is analytic on
$$
D^n_{r_*/ \eta^2}\times\torus^n_{s_*}\times
D^{2m}_{\r_*/ \eta}\,.
$$
In order to find periodic solutions with large period
close to the elliptic torus ${\cal T}$ we will apply
in section \ref{section theorem 1} a finite dimensional reduction of
Lyapunov-Schmidt type.
For this purpose (see remark \ref{remdt}) we need first to perform
a symplectic change of variables which eliminates
from the Hamiltonian $ \cH $ defined in (\ref{ham (2)})
the following terms\footnote{For short,
in the sequel we will often omit the summation over $\ell \in \integer^n$.}
\begin{eqnarray}
&&
\eta\sum_{1\leq i\leq n\atop{|a+\ba|=1}}
R_{e_i, a, \overline a}^*(\f)
\cI_i Z^a \overline Z^{\overline a}\ ,\label{pb1}\\
&&
\eta^2
\sum_{|k|=2, \ell\neq 0}
R_{k, 0, 0,\ell}^* \cI^k e^{\ii\ell\cdot\f}\ , \label{pb22}\\
&&
\eta^2
\sum_{1\leq i\leq n\atop{{|a+\ba|=2
}\atop{a\neq \ba\,{\rm or}\, \ell\neq 0}}}
R_{e_i, a,\ba,\ell}^*
\cI_i Z^a \bZ^\ba e^{\ii\ell\cdot\f}\ ,\label{pb21}\\
&&
\eta^3\sum_{2|k|+|a+\ba|=5\atop{|a+\ba|=1}}
R_{k, a, \overline a}^*(\f)
\cI^k Z^a \overline Z^{\overline a}\ .\label{pb3}
\end{eqnarray}
This task will be
accomplished in the next Proposition.
The term \equ{pb22} will be ``averaged out'' since $ \o $ is diophantine;
the terms \equ{pb1} and \equ{pb3} using
the first order Melnikov non-resonance conditions (namely
conditions \equ{CONDIZ NON RISON} for $|h| \leq 1 $),
and, finally, the term
\equ{pb21} using the second order Melnikov non-resonance
conditions \equ{CONDIZ NON RISON}.

\begin{proposition} {\bf(Averaging)} \label{averaging}
Let the Hamiltonian $\cH $, defined in \equ{ham (2)}, be
real analytic on $D^n_{r}\times \torus^n_{s}\times D^{2m}_{\r}$ and
satisfy the second-order Melnikov
non-resonance conditions \equ{CONDIZ NON RISON}.
Then, for $ \eta $ small enough, depending on $n$, $m$, $\o$, $\O$,
$\gamma$, $\tau$, $r$, $s$, $\r$,
there exists an analytic canonical change of coordinates $ \Phi : $
$ (I, \phi, z, \overline{z}) \to ( \cI, \f ,  Z, \overline Z)$,
$\eta$-close to the identity,
$$
\Phi\,: \, D^n_{r/2}\times
\torus^n_{s/2}\times D^{2m}_{\r/2}
\,\longrightarrow
\,
D^n_{r}\times
\torus^n_{s}\times D^{2m}_{\r}\,,
$$
transforming the Hamiltonian $ \cH $ into
the Hamiltonian\footnote{
We denote
$\cQ I\cdot z\bz:=
\sum_{1\leq j\leq m\atop{1\leq i\leq n}} \cQ_{ji} I_i z_j \bz_j$.
}
\beqa{ham (****)}
H = {\cal H}\circ
\Phi &=& \o \cdot I+ \O z\bz+
\eta \sum_{|a+\ba|=3} R^*_{0,a,\ba}(\phi)z^a \bz^\ba +\nonumber\\
&+&
\eta^2\Big[ \frac{1}{2} \cR I \cdot I +
\cQ I\cdot z\bz +
\sum_{|a+\ba|=4} R_{0,a,\ba}(\phi)z^a \bz^\ba
\Big]+
\nonumber\\
&+&
\eta^3
\sum_{2|k|+|a+\ba|=5\atop{|a+\ba|=3,5}} R_{k,a,\ba}(\phi) I^k z^a \bz^\ba
+ O(\eta^4)\ ,
\eeqa
for suitable $ R_{k, a,\ba} \in \complex $,
$ R_{k,a,\ba}(\phi)$ analytic on $\torus^n_{s/2}$, and
where $ \cR $ is the symmetric twist matrix defined in \equ{R2}
and $\cQ \in {\rm Mat}(m \times n, \real) $ is defined in \equ{Q}.
\end{proposition}

\noindent {\bf Proof:}
We rewrite the Hamiltonian (\ref{ham (2)}) in the form
\begin{equation}
\cH (\cI_,\f,Z,\overline Z)=\sum_{j\geq 0}\eta^j R_*^{(j)}(\cI_,\f,Z,\overline Z)\
\end{equation}
with
\begin{equation}\label{FOURIER}
R_*^{(0)}:=\o \cdot \cI +\O Z\cdot \overline Z \qquad {\rm and} \qquad
R_*^{(j)}:=
\sum_{2|k|+|a+\ba|=j+2}
R_{k,a,\ba,\ell}^*
\cI^k Z^a \bZ^{\ba}e^{\ii\ell\cdot\f}\
\end{equation}
where we have expanded in Fourier series the functions
$R_{k,a,\ba}^*(\f)$ as in \equ{1'}.

We want to define a canonical change of variables
$\Phi$ from the new variables $\z:=(I,\phi,z,\bz)$
to the old ones $\cZ:=(\cI,\f,Z,\bZ)$, as
the flow at time 1 of a suitable Hamiltonian
$\chi$. Precisely $\Phi:=\Phi_{\chi}^1,$
where $\Phi_{\chi}^t(\z_0):=\z(t)$ is the unique solution of
the Hamiltonian system
\begin{equation}
\dot{I}(t)=-\partial_{\phi} \chi(\z(t))\ ,\quad
\dot{\phi}(t)=\partial_{I} \chi(\z(t))\ ,\quad
\dot{z}(t)=\ii\partial_{\bz} \chi(\z(t))\ ,\quad
\dot{\bz}(t)=-\ii\partial_{z} \chi(\z(t))\
\end{equation}
with initial conditions $ \z(0) = \z_0 $.

The Lie operator (Poisson brackets),
acting on a function $g:=g(\z)$, is defined as
\begin{equation}
L_{\chi}g:=\{ g,\chi\}:=\partial_{\phi}g\partial_{I}\chi-
\partial_{I}g\partial_{\phi}\chi
+\ii\partial_{z}g\partial_{\bz}\chi
-\ii\partial_{\bz}g\partial_{z}\chi.
\end{equation}
For every integer $ j \geq 0 $ we also set
\begin{equation}
L_{\chi}^0 g:=g\ , \qquad L_{\chi}^j g:=L_{\chi}L_{\chi}^{j-1} g\ .
\end{equation}
The new Hamiltonian $ H = \cH \circ \Phi $ can be developed,
for all $ j_0 \in \natural^+$,  as
\begin{equation}\label{sviluppoLie}
H:=\cH\circ\Phi=\sum_{j=0}^{j_0}\frac{1}{j!}L_{\chi}^j \cH+
\frac{1}{j_0!}\int_0^1\ (1-\xi) L_{\chi}^{j_0+1}\cH\circ\Phi_\chi^\xi\ d\xi\ .
\end{equation}
We look for a Hamiltonian $\chi $ of the form
\begin{equation}\label{chi1}
\chi =\sum_{j=1}^{j_0}\eta^j \chi^{(j)}\ ,
\quad {\rm where }\quad
\chi^{(j)}:=\chi^{(j)}(I,\phi,z,\bz):=
\sum_{2|k|+|a+\ba|=j+2}
\chi_{k,a,\ba,\ell}
I^k z^a \bz^{\ba}e^{\ii\ell\cdot\phi}\
\end{equation}
will be chosen later on.
For the sake of simplicity we will often {\em omit}
in the notation the summation over $\ell\in \integer^m$.
Inserting (\ref{chi1}) in (\ref{sviluppoLie})
we obtain
\begin{equation}\label{Hprima}
H=\sum_{d=0}^{j_0}\eta^dR^{(d)}+
\sum_{d=1}^{j_0}\eta^d\sum_{j=0}^{j_0}\frac{1}{j!}
\sum_{h=0}^{d-1}\
\sum_{{h+i_1+\dots+i_j=d}\atop{1\leq i_1,\dots i_j\leq d}}
L_{\chi^{(i_1)}}\dots
L_{\chi^{(i_j)}}R_*^{(h)}+\ O(\eta^{j_0+1})\
\end{equation}
Denoting $[ \ \cdot \ ]_d$ the $d$-th order in $\eta$,
we obtain
\begin{equation}\label{ordine1}
[H]_1=R_*^{(1)}+L_{\chi^{(1)}}R_*^{(0)}
\end{equation}
and
\begin{equation}\label{ordined}
[H]_d=R_*^{(d)}+L_{\chi^{(d)}}R_*^{(0)}+
\sum_{j=0}^{j_0}\frac{1}{j!}
\sum_{h=0}^{d-1}\
\sum_{{h+i_1+\dots+i_j=d}\atop{1\leq i_1,\dots i_j\leq d-1}}
L_{\chi^{(i_1)}}\dots
L_{\chi^{(i_j)}}R_*^{(h)}\,.
\end{equation}
Observe that
\begin{eqnarray}
&&\sum_{\scriptsize{
\begin{array}{c}
2|k|+|a+\ba|=h+2 \\
2|k'|+|a'+\ba'|=h'+2
\end{array}}
}
\{ I^k z^a \bz^{\ba}e^{\ii\ell\cdot\phi}\ ,\
I^{k'} z^{a'} \bz^{\ba'}e^{\ii\ell'\cdot\phi} \}=\\
&&\sum_{2|k''|+|a''+\ba''|=h+h'+2}
c_{k'',a'',\ba'',\ell''} I^{k''} z^{a''} \bz^{\ba''}e^{\ii\ell''\cdot\phi}\ ,
\end{eqnarray}
for suitable constants $c_{k'',a'',\ba'',\ell''}$,
which are explicitly given by the
following formula:
\begin{eqnarray}\label{broop}
\{ I^k z^a \bz^{\ba}e^{\ii\ell\cdot\phi}\ ,\
I^{k'} z^{a'} \bz^{\ba'}e^{\ii\ell'\cdot\phi} \} & = & \ii\Big[
\sum_{i=1}^n(\ell_i k_i'-\ell_i' k_i)I^{k+k'-e_i}z^{a+a'}\bz^{\ba+\ba'}+
\\
&+ &
\sum_{j=1}^m(a_j\ba_j'-\ba_j a_j')I^{k+k'}z^{a+a'-e_j}\bz^{\ba+\ba'-e_j}
\Big] e^{\ii(\ell+\ell')\cdot\phi}\nonumber\,.
\end{eqnarray}
Thus, we obtain that
\begin{equation}
L_{\chi^{(i_1)}}\dots
L_{\chi^{(i_j)}}R^{(h)}_*=
\sum_{2|k|+|a+\ba|=h+i_1+\dots+i_j+2}
c_{k,a,\ba,\ell}I^k z^a \bz^{\ba}e^{\ii\ell\cdot\phi}
\nonumber
\end{equation}
for suitable constants $c_{k,a,\ba,\ell}$.
Hence,
\begin{equation}\label{ordinedbis}
[H]_d=L_{\chi^{(d)}}R^{(0)}_*+
\sum_{2|k|+|a+\ba|=d+2} R_{k,a,\ba,\ell}I^k z^a \bz^{\ba}e^{\ii\ell\cdot\phi}
\end{equation}
for suitable $R_{k,a,\ba,\ell}$ with
\begin{equation}
 R_{k,a,\ba,\ell}:= R_{k,a,\ba,\ell}
(\chi^{(1)},\dots,\chi^{(d-1)},R^{(0)}_*,\dots,R^{(d)}_*)\ ;
\nonumber
\end{equation}
we note that, by means of \equ{ordine1}
and \equ{ordinedbis}
\begin{equation}\label{R=R*}
2|k|+|a+\ba|=3\qquad \Longrightarrow\qquad
R_{k,a,\ba,\ell}=R_{k,a,\ba,\ell}^*\ ,
\end{equation}
(recall also the setting in \equ
{FOURIER}).
We also evaluate
\begin{equation}\label{Lphid}
L_{\chi^{(d)}}R^{(0)}_*=
\{ \o\cdot I  + \O z\bz,\chi^{(d)}\} =
\sum_{2|k|+|a+\ba|=d+2}
-\ii(\o\cdot\ell+\O\cdot(a-\ba))\chi_{k,a,\ba,\ell}
I^k z^a \bz^{\ba}e^{\ii\ell\cdot\phi}\ .
\end{equation}

We define the following ``resonant'' set:
\begin{equation}\label{S}
S:=S_1\cup S_2\cup S_3\subset \Big\{
(k,a,\ba,\ell)\quad {\rm s.t.}\quad 3\leq 2|k|+|a+\ba|\leq 5\ ,\quad\ell\in
\integer^n \Big\}\ ,
\end{equation}
with $S_2=S_2^0\cup S_2^1\cup S_2^2$ and
\beqano
S_1 &:=& \{ k = 0, \ |a+\ba|=3,\ \ell\in\integer^n \} \\
S_2^0 &:=& \{ k=0,\ |a+\ba|=4,\ \ell\in\integer^n \} \\
S_2^1 &:=& \{ k=e_i,\ a=\ba=e_j,\ 1\leq i\leq n,\ 1\leq j\leq m,\ \ell=0 \} \\
S_2^2 &:=& \{ |k|=2,\ a=\ba=0, \ \ell=0 \} \\
S_3 &:=& \{ 2|k|+|a+\ba|=5,\ |a+\ba|=3,5, \ \ell\in\integer^n \} \ .
\eeqano
Let $ j_0 = 3 $.
$ \forall \ 1 \leq d \leq 3 $, $ \forall (k, a, \ba, \ell ) $ such that
$ 2|k|+|a+\ba|=d+2, \ell \in \integer^n $, one can check,
by condition \equ{CONDIZ NON RISON}, that
\beq{Spicden}
(k,a,\ba,\ell)\notin S\ \ \ \Longrightarrow\ \ \
|\o\cdot\ell+\O(a-\ba)|\geq\frac{\g}{1+|\ell|^\tau}\ ,
\eeq
and hence we can define
\begin{equation}\label{chi2}
\chi_{k,a,\ba,\ell}:=\left\{
{\begin{array}{ll}
0 & {\rm if } \ \ (k,a,\ba,\ell)\in S\\
-\ii[\o\cdot\ell+\O(a-\ba)]^{-1}R_{k,a,\ba,\ell}\qquad & {\rm otherwise}
\end{array}}\right.
\end{equation}
In light of this construction, using
(\ref{ordinedbis}) and (\ref{Lphid}),
we have\footnote{We denote by $\Pi_{{\cal S}}$,
where ${\cal S}\subset \natural^{n+2m}\times\integer^n$,
the projection on ${\cal S}$ $i.e.$
$$
\Pi_{{\cal S}}\Big( \sum_{k,a,\ba,\ell} c_{k,a,\ba,\ell}I^kz^a \bz^\ba
e^{\ii\ell\cdot\phi}\Big)
:=\sum_{(k,a,\ba,\ell)\in{\cal S}} c_{k,a,\ba,\ell}I^kz^a \bz^\ba
e^{\ii\ell\cdot\phi}\ .
$$}
\beq{HdPSd}
[H]_d=
\Pi_{S_d}\Big( \sum_{2|k|+|a+\ba|=d+2} R_{k,a,\ba,\ell}I^kz^a \bz^\ba
e^{\ii\ell\cdot\phi}\Big)
=\sum_{(k,a,\ba,\ell)\in{S_d}} R_{k,a,\ba,\ell}I^kz^a \bz^\ba
e^{\ii\ell\cdot\phi}\ .
\eeq
Define
$$
R_{k,a,\ba}:=R_{k,a,\ba,0}\ ,\qquad
R_{k,a,\ba}(\phi):=\sum_{\ell\in\integer^n} R_{k,a,\ba,\ell}
e^{\ii \ell\cdot\phi} \ .
$$
By recurrence, using \equ{ordinedbis},
it is possible to evaluate the terms $R_{k,a,\ba}$ explicitely.

From \equ{ordined} and \equ{R=R*} we have
$$
[H]_2=\Pi_{S_2}\left( \frac{1}{2}\{\{  R^{(0)}_*,\chi^{(1)}\},\chi^{(1)}\}
+ \{  R^{(1)}_*,\chi^{(1)}\}+R^{(2)}_*\right)  \ .
$$
Noting that $$ \{ R^{(0)}_*,\chi^{(1)}\}=-R^{(1)}_*+\Pi_{S_1}
R^{(1)}_*\ ,$$
we have
$$
[H]_2=\Pi_{S_2}\left( \frac{1}{2}
\{  R^{(1)}_*+\Pi_{S_1}R^{(1)}_*,\chi^{(1)}\}+R^{(2)}_*\right)  \ .
$$
In order to prove \equ{ham (****)} we only need to show that
\beqa{R00}
\frac{1}{2} \cR I \cdot I
& = &
\Pi_{S_2^2} [H_2] =\frac{1}{2}\Pi_{S_2^2}\left(
\{  R^{(1)}_*,\chi^{(1)}\}\right)+\Pi_{S_2^2} R^{(2)}_*
  \ ,
\\
\cQ I\cdot z\bz
& = &
\Pi_{S_2^1} [H_2] =
\Pi_{S_2^1}\left( \frac{1}{2}
\{  R^{(1)}_*+\Pi_{S_1}R^{(1)}_*,\chi^{(1)}\}\right)
+\Pi_{S_2^1} R^{(2)}_* \nonumber\\
& = &
\Pi_{S_2^1}
 \left(\left\{
\frac{1}{2}
  \Pi_{S_1^c} R^{(1)}_*
+
\Pi_{S_1}R^{(1)}_*,
\chi^{(1)}\right\}\right)
+\Pi_{S_2^1} R^{(2)}_*
\ ,
\label{R12}
\eeqa
where $S_1^c := \{ k=e_i,\, |a+\ba|=1,\, \ell\in\integer^n,\,
1\leq i\leq n  \}.$

We first prove \equ{R00}.
Using (\ref{chi1}), (\ref{chi2})
and (\ref{broop}),
we have
\begin{eqnarray}
\Pi_{S_2^2} [H_2]
& = &
\sum_{j=1}^m
\sum_{\scriptsize{
\begin{array}{c}
|k|=|k'|=1 \\
|a+\ba|=|a'+\ba'|=1\\
a+a'=\ba+\ba'=e_j
\end{array}}
}
\,\frac{1}{2}\,\cdot\,
\frac{R_{k,a,\ba,\ell}^*R_{k',a',\ba',-\ell}^*}{
-\o\cdot\ell+\O(\ba-a)}
(a_j\ba_j'-\ba_j a_j') I^{k+k'}\nonumber\\
& + &
\sum_{|k|=2}R_{k,0,0,0}^*I^k
\nonumber\\
& = &
\sum_{j=1}^m
\sum_{\scriptsize{
\begin{array}{c}
1\leq i,i'\leq n \\
|a+\ba|=|a'+\ba'|=1\\
a+a'=\ba+\ba'=e_j
\end{array}}
}
\,\frac{1}{2}\,\cdot\,
\frac{R_{e_i,a,\ba,\ell}^*R_{e_{i'},a',\ba',-\ell}^*}{
-\o\cdot\ell+\O(\ba-a)}
(a_j\ba_j'-\ba_j a_j') I_i I_{i'}\nonumber\\
& + &
\sum_{i,i'=1}^n \frac{1 + \d_{i,i'}}{2}R_{e_i+e_{i'},0,0,0}^*I_i I_{i'}
\nonumber\\
& = &
\frac{1}{2}\sum_{i,i'=1}^n
\sum_{j=1}^m
\left(
\frac{R_{e_i,e_j,0,\ell}^*R_{e_{i'},0,e_j,-\ell}^*}{
-\o\cdot\ell-\O_j}-
\frac{R_{e_i,0,e_j,\ell}^*R_{e_{i'},e_j,0,-\ell}^*}{
-\o\cdot\ell+\O_j}
\right)
\nonumber\\
& + &
\sum_{i,i'=1}^n \frac{1 + \d_{i,i'}}{2}R_{e_i+e_{i'},0,0,0}^*I_i I_{i'}
= \frac{1}{2} \cR I \cdot I \ ,\label{R00ok}
\end{eqnarray}
where $\d_{i,i'} = 1 $ if $i = i' $ and $0 $ if $ i \neq i' $.
Here, we observe that, since $|a+\ba|=|a'+\ba'|=1$, if
$a+a'=\ba+\ba'=e_j$, then
$(a,\ba,a',\ba')=
(e_j,0,0,e_j)$ or $(0,e_j,e_j,0)$.
So \equ{R00} directly follows from \equ{R00ok}.

We now prove \equ{R12}.
Using again (\ref{chi1}), (\ref{chi2})
and (\ref{broop}),
we have
\beqano
\Pi_{S_2^1} [H_2]
& = &
-\ii\, \Pi_{S_2^1}
 \Bigg(\Bigg\{
\frac{1}{2}
  \Pi_{S_1^c} R^{(1)}_*
+
\Pi_{S_1}R^{(1)}_*\ ,\\
&\empty&\qquad\qquad
\sum_{1\leq i'\leq n\atop{1\leq j'\leq m}}
\left(
\frac{R_{e_{i'},e_{j'},0,\ell'}^*}{
\o\cdot\ell'+\O_{j'}}z_{j'}\,
+
\,
\frac{R_{e_{i'},0,e_{j'},\ell'}^*}{
\o\cdot\ell'-\O_{j'}}\bz_{j'}
\right)
I_{i'} e^{\ii \ell'\cdot\phi}
\Bigg\}\Bigg)
+\Pi_{S_2^1} R^{(2)}_*
\eeqano
\beqa{R12ok}
& = &
\frac{1}{2}\Pi_{S_2^1}
\sum_{\scriptsize{
\begin{array}{c}
1\leq i,i'\leq n \\
1\leq j'\leq m\\
|a+\ba|=1
\end{array}}
}
(\ell_{i'}I_i+\ell_i I_{i'})
R_{e_i,a,\ba,\ell}^* z^a \bz^\ba
\left(
\frac{R_{e_{i'},e_{j'},0,-\ell}^*}{
-\o\cdot\ell+\O_{j'}}z_{j'}
+
\frac{R_{e_{i'},0,e_{j'},-\ell}^*}{
-\o\cdot\ell'-\O_{j'}}\bz_{j'}
\right)
\nonumber\\
&\empty&
+\Pi_{S_2^1}
\Bigg(
\sum_{1\leq i\leq n\atop{1\leq j,j'\leq m}}
\left\{
R_{0,e_j,e_j+e_{j'},\ell}^* \bz_{j'} z_j\bz_j e^{\ii \ell\cdot\phi}\ ,
\frac{R_{e_{i'},e_{j'},0,\ell'}^*}{
\o\cdot\ell'+\O_{j'}}z_{j'} I_{i'} e^{\ii \ell'\cdot\phi}
\right\}
\nonumber\\
&\empty& \qquad\qquad  +
\sum_{1\leq i\leq n\atop{1\leq j,j'\leq m}}
 \left\{
R_{0,e_j+e_{j'},e_j,\ell}^* z_{j'} z_j\bz_j e^{\ii \ell\cdot\phi}\ ,
\frac{R_{e_{i'},0,e_{j'},\ell'}^*}{
\o\cdot\ell'-\O_{j'}} \bz_{j'} I_{i'} e^{\ii \ell'\cdot\phi}
\right\}
\Bigg)+\Pi_{S_2^1} R^{(2)}_*
\nonumber\\
& = &
\frac{1}{2}
\sum_{1\leq i,i'\leq n\atop{1\leq j\leq m}} \ell_{i'}
\Bigg(
\frac{
R_{e_{i},0,e_{j},\ell}^*
R_{e_{i'},e_{j},0,-\ell}^*}{
-\o\cdot\ell+\O_{j}}+
\frac{
R_{e_{i},e_{j},0,\ell}^*
R_{e_{i'},0,e_{j},-\ell}^*}{
-\o\cdot\ell-\O_{j}}
\nonumber\\
& & \qquad\qquad\qquad\qquad
+\frac{
R_{e_{i},0,e_{j},\ell}^*
R_{e_{i'},e_{j},0,-\ell}^*}{
-\o\cdot\ell+\O_{j}}+
\frac{
R_{e_{i},e_{j},0,\ell}^*
R_{e_{i'},0,e_{j},-\ell}^*}{
-\o\cdot\ell-\O_{j}}
\Bigg)
I_i z_j \bz_j
\nonumber\\
&\empty &
-\sum_{1\leq i\leq n\atop{1\leq j,j'\leq m}}
\frac{1}{
\o\cdot\ell+\O_{j'}}
\left(
R_{0,e_{j},e_{j}+e_{j'},-\ell}^*
R_{e_{i},e_{j'},0,\ell}^* +
R_{0,e_{j}+e_{j'},e_{j}\ell}^*
R_{e_{i},0,e_{j'},-\ell}^*
\right)
I_i z_j \bz_j
\nonumber\\
&\empty&
+
\Pi_{S_2^1}
\sum_{2|k|+|a+\ba|=4}
R_{k,a,\ba,\ell}^* I^k z^a \bz^\ba e^{\ii\ell\cdot\phi}
\ .
\eeqa

So \equ{R12} directly follows from \equ{R12ok}.

The estimate on the new analyticity radii follows
from \equ{Spicden} and from the fact that
$\chi=O(\eta)$ (see \equ{chi1}).
\hfill$\Box$

\section{Periodic orbits winding along the torus}\label{section theorem 1}

\indent
In this section we prove the existence
of periodic solutions of longer and longer (minimal) period
shrinking closer and closer to the elliptic torus $ { \cal T} $
(Theorem \ref{theorem 1}).
\\[1mm]
\indent
By Proposition \ref{averaging} the Hamiltonian \equ{ham (2)}
can be transformed, thanks to the second order Melnikov
non-resonance conditions \equ{CONDIZ NON RISON}, into
the Hamiltonian
\beqa{ham (3)}
H (I,\phi,z,\overline z)
&=&
\o\cdot I+\O z\bz+
\eta \sum_{|a+\ba|=3} R^*_{0,a,\ba}(\phi)z^a \bz^\ba \nonumber\\
&+&
\eta^2\Big[ \frac{1}{2} \cR I \cdot I + {\cal Q}I \cdot z \overline{z} +
 \sum_{|a+\ba|=4} R_{0,a,\ba}(\phi)z^a \bz^\ba
\Big]
\nonumber\\
&+&
\eta^3
 \sum_{2|k|+|a+\ba|=5\atop{|a+\ba|=3,5}} R_{k,a,\ba}(\phi) I^k z^a \bz^\ba
+O(\eta^4)\ ,
\eeqa
analytic on
$$ D^n_{r_*/(2\eta^2)}\times\torus^n_{s_*/2}\times
D^{2m}_{\r_*/(2\eta)}\,.$$

The Hamilton's equations of motion induced by the Hamiltonian \equ{ham (3)}
\begin{equation}\label{Hequat}
\dot{I} = - \partial_{\phi } H,\quad
\dot{\phi} = \partial_{ I } H, \quad
\dot{z} = \ii \partial_{\overline z } H,\quad
\dot{\overline z} = - \ii \partial_{z} H, \
\end{equation}
can be written as
\begin{eqnarray}\label{equation of motion}
\dot I &=& -
\eta \sum_{|a+\ba|=3} \partial_\phi R^{*}_{0,a,\ba}(\phi)z^a \bz^\ba-
\eta^2 \sum_{|a+\ba|=4} \partial_\phi R_{0,a,\ba}(\phi)z^a \bz^\ba+\nonumber\\
&&-\eta^3
\sum_{2|k|+|a+\ba|=5\atop{|a+\ba|=3,5}}
\partial_\phi R_{k,a,\ba}(\phi) I^k z^a \bz^\ba
+O(\eta^4)\,,\nonumber\\
\dot \phi_i &=& \o_i +
\eta^2\Big[ (\cR I)_i + \sum_{1\leq j \leq m} {\cal Q}_{ji} z_j \bz_j
\Big]+
\nonumber\\
&&\eta^3
 \sum_{2|k|+|a+\ba|=5\atop{|a+\ba|=3,5}}
R_{k,a,\ba}(\phi) k_i I^{k-e_i} z^a \bz^\ba
+O(\eta^4)\,,
\nonumber\\
\dot z_j &=& \ii\O_j z_j +
\ii\eta \sum_{|a+\ba|=3} R^*_{0,a,\ba}(\phi)\ba_j z^a \bz^{\ba-e_j}+
\nonumber\\
&&
\ii\eta^2\Big[ ({\cal Q}I)_j z_j  +
 \sum_{|a+\ba|=4} R_{0,a,\ba}(\phi) \ba_j z^a \bz^{\ba-e_j}
\Big]+
\nonumber\\
&&
\ii\eta^3
 \sum_{2|k|+|a+\ba|=5\atop{|a+\ba|=3,5}} R_{k,a,\ba}(\phi) I^k
\ba_j z^a \bz^{\ba-e_j}+O(\eta^4)\,,
\end{eqnarray}
for $i=1,\dots, n$, $j=1,\dots ,m$.
\\[1mm]
\indent
Critical points of the Hamiltonian action functional
\beq{actionfun}
{\cal A}(I(t), \phi(t), z(t), \bz (t)  ) :=
\int_0^T I \cdot \dot \phi + \ii z \dot \bz
- H ( I(t), \phi(t), z(t), \bz (t) ) \,dt,
\eeq
in a suitable space of $T$-periodic functions, are $T$-periodic
solutions of \equ{Hequat}.

\subsection{The pseudo periodic solutions}\label{pseudo 1}

We will find periodic solutions of the Hamiltonian system
\equ{Hequat} close to periodic solutions of the integrable Hamiltonian
\beq{Hint}
H_{int} := \o \cdot I + \frac{\eta^2}{2} \cR I \cdot I +
\O z\bz  + \eta^2 {\cal Q} I \cdot z \bz.
\eeq
The manifold
$ \{ z = 0 \} $ is invariant for the Hamiltonian system
generated by $ H_{int}$\footnote{We recall
the usual notation for the vector $\cQ^T z\bz\in\real^n$
whose components are
$(\cQ^T z\bz)_i:=\sum_{j=1}^m\cQ_{ji}I_i z_j \bz_j$, $1\leq i\leq n$ and for
the $m\times m$ diagonal matrix $\O + \eta^2 {\cal Q} I:=$ $ {\rm diag}
$ $\Big(\O_1 + \eta^2 ({\cal Q} I)_1,\dots ,
\O_m + \eta^2 ({\cal Q} I)_m \Big)$. }
\beq{HSint}
\dot{I} = 0, \qquad
\dot{\phi} = \o + \eta^2 (\cR I + {\cal Q}^T z {\overline z} ), \qquad
\dot{z} = \ii (\O + \eta^2 {\cal Q} I ) z
\eeq
and it is completely filled up by the invariant tori
$ {\cal T}(I_0) := \{ I = I_0, \ \phi \in \torus^n , \ z = 0 \} $.
On $ { \cal T}(I_0) $ the flow is $t \to
\{ I_0 , \ \phi_0 + (\o + \eta^2 \cR I_0) t, 0 \} $
and in its normal space it is determined by
$ \dot{z} = \ii \O_\eta (I_0) z $
where $\O_\eta (I_0) $ is the $ m \times m $ diagonal matrix  with
{\em real}
coefficients associated to the vector
of the ``shifted elliptic frequencies''
\begin{equation}\label{omega eta}
\O_{\eta}(I_0) := \O +\eta^2 {\cal Q} I_0 \ ,
\end{equation}
For any $I_0\in \real^n$, $ \O_\eta ( I_0 )$ is a real matrix,
since $ \cQ $ is real (see
page \pageref{real2}).

For suitable $T>0$, $ I_0 \in \real^n $, $k\in\integer^n$, namely
when
\beq{otilde}
\widetilde \o := \o+\eta^2 \cR I_0 = \frac{1}{T}
2 \pi k \in \frac{1}{T} 2 \pi \integer^n,
\eeq
the torus $ {\cal
T}(I_0) $ is {\it completely resonant}, supporting the
family of $T$-periodic motions \beq{unpert} {\cal P} := \Big\{
I(t) = I_0, \quad \phi (t) = \phi_0 + \widetilde \o t, \quad z(t)
= 0 \Big\}. \eeq The family $ { \cal P} $ will not persist in its
entirety for the complete Hamiltonian system \equ{Hequat}.
However, the non-resonance property \equ{Mmeno1} below, between
the period $T$ and the ``shifted elliptic frequenies'' $ \O_\eta
(I_0) $, is sufficient to prove the persistence of at least $ n $
geometrically distinct $T$-periodic solutions of the Hamiltonian
system \equ{Hequat}, close to $ { \cal P} $. Precisely, the
required {\it non resonance property} is
\beq{Mmeno1}
\cM := \cM (I_0,
T) := {\bf 1}_m - e^{\ii \O_\eta (I_0) T} \quad {\rm is \
invertible} \quad {\rm and} \quad | \cM^{-1}(I_0, T) |\leq\,{\rm
const}.
\eeq
Our aim is then to find $ I_0 $ and $ T $ so that
\equ{otilde} and \equ{Mmeno1} hold: we will define $ I_0 := I_0
(T)$ in dependence on the ``$1$-dimensional parameter'' $T$ in
such a way that \equ{otilde} is identically satisfied and then we
will find  $ T $ so that the non resonance property \equ{Mmeno1}
holds. Moreover, for our perturbative arguments of Lemma
\ref{pseudo (1)}, we want $ I_0 = O ( 1 ) $ and $ T < 2 \slash
\eta^{2} $.

Define, for $ T \geq 1 \slash \eta^2 $,
\beqa{I0}
I_0 &:=& I_0(T)
:= -\frac{2\pi}{\eta^2 T} \cR^{-1} \Big\langle \frac{\o T}{2 \pi}
\Big\rangle \ ,\\
k &:=& k(T) = \frac{\o T}{2\pi}- \Big\langle \frac{\o T}{2\pi}
\Big\rangle\ ,\label{kT}
\eeqa
where $\langle (x_1,\dots,x_n) \rangle:=(\langle
x_1\rangle,\dots, \langle x_n\rangle)$ and the function
$\langle \cdot \rangle : \real \to [-1/2,1/2) $
is defined as $\langle x \rangle:=x$ for $x\in[-1/2,1/2)$
and it is $1$-periodically extended for $ x \in \real.$
Notice that $I_0\in \real^n$ since $\cR$ is a {\em
real} matrix (see page \pageref{real2}). With the choice
\equ{I0},\equ{kT}, $ \o T + \eta^2 \cR I_0(T) = 2 \pi k $, and
then \equ{otilde} holds. In addition, for $ T \geq 1 \slash \eta^2
$, $ I_0 (T) = O ( 1 ) $. Moreover \beq{OeT} \O_\eta T = \O_\eta
(I_0 (T)) T = \O T+\eta^2 \cQ I_0 T = 2\pi\left(
\O\frac{T}{2\pi}-\cQ \cR^{-1} \Big\langle \o \frac{T}{2\pi}
\Big\rangle \right). \eeq In order to  prove the non-resonance
property  \equ{Mmeno1}, we note that \beq{stimaM} |\cM^{-1} |\leq
\frac{1}{\displaystyle{\min_{1\leq j\leq m}} |e^{\ii \O_{\eta
j}T}-1|}\leq \frac{1}{\displaystyle{\min_{1\leq j\leq m}}
|\sin(\O_{\eta j}T)|}\leq \frac{2}{\displaystyle{\min_{1\leq j\leq
m}} {\rm dist} (\O_{\eta j}T,2\pi\integer) }\ . \eeq

\lem{a} Suppose that  condition $(a)$ of Theorem \ref{theorem 1}
hold. Define, for $ \hm \geq 1 $, \beq{delta} d_0 := \min \left\{
\frac{1}{4\displaystyle{\max_{ 1 \leq j \leq \hm} |a_j|_1}},
\frac{1}{4\displaystyle{\max_{ 1 \leq j \leq \hm} M_j}} \right\},
\ \ \d := \min \left\{ \frac{d_0}{2\displaystyle{\max_{1\leq j\leq
m} |(\cQ\cR^{-1})_j|_1}}, \frac{1}{4\displaystyle{\max_{1\leq
j\leq \hm}}|a_j|_1} \right\} \eeq where $ (\cQ\cR^{-1})_j $ is the
$j$-row of the matrix $ \cQ \cR^{-1} $, and, for $ \hm = 0 $,
\beq{deltabis} d_0 := \frac{1}{4}\ ,\qquad \qquad   \d :=  \min
\left\{ \frac{1}{8\displaystyle{\max_{1\leq j\leq m}
|(\cQ\cR^{-1})_j|_1}}\ , \frac{1}{4} \right\}. \eeq Let $ M :=\,
{\rm lcm}\, (M_1,\dots,M_\hm),$ for $\hm\geq 1$, and $M=1$ for $
\hm = 0 .$ Finally let\footnote{ The ergodization time $T_{erg}$
was defined in \equ{ergtime}.} \beq{mandos} T_e:=
T_{erg}\left(\frac{\ho}{M},\frac{\d}{M}\right), \qquad \Theta :=
\min \left\{ \frac{\d}{4\displaystyle{\max_{1\leq i\leq
n}|\o_i|}}\ , \frac{d_0}{8\displaystyle{\max_{1\leq j\leq
m}}|\O_j|} \right\}. \eeq Then $\forall\, t_0\geq 0,$ there exists
an interval ${\cal J}  \subset  [t_0 -2\pi\Theta, t_0 +  2\pi T_e
+ 2\pi\Theta ]$ of length, at least $ 4\pi \Theta$, such that
$\forall\, T\in {\cal J}$,
 $ \cM^{-1} := \cM^{-1}( I_0 (T), T ) $ satisfies
\beq{d2b}
|\cM^{-1} |\leq \frac{4}{\pi d_0}
\ .
\eeq
\elem
\noindent {\bf Proof:}
In order to define the periods $ T $ we want to use Lemma \ref{tau}.
Let us verify that its hypothesis are fulfilled.
If $(\o, \O)$ are rationally independent, i.e.
$ \hm = 0 $, we can apply it directly. If $ \hm = 1 $ we observe that
\equ{2} holds by definition, since $M_j\geq 1=\hm.$
If $ \hm = 1 = M_1 $ we have also to show that \equ{tulkas}
is satisfied.
This is true: indeed, by contradiction, if \equ{tulkas}
were false then $\O_1=\sum_{i=1}^n a_{1i}\o_i$
violating the {\sl first} order Melnikov condition
(namely \equ{CONDIZ NON RISON} with $|h|\leq 1$).
If $ m = 1 $ then $ \hm = 0, 1 $ and we are in one of
the previous cases. We finally consider the case $ \hm \geq 2$.
We have to prove \equ{2}.
Note first that it is verified for $ m = 2$.
Indeed in this case $ \hm = 2 $ and, hence, $ \ho = \o$.
By definition $ M_j \O_j = a_j \cdot \o $ for $j=1,2.$
Again by the {\sl first} order Melnikov condition
we have, for $j=1,2$, that $M_j\geq 2$ which implies
$M_j\geq 2=\hm$ and \equ{2} holds.
All the remanent cases are covered by the hypotesis $(iii)$.
We can apply Lemma \ref{tau}.

Let $ t_0 \leq \tau \leq t_0 + T_e $ be the time
found in Lemma \ref{tau} 
and consider the interval
$ {\cal J} := 2 \pi \tau + 2 \pi [ - \Theta, \Theta ]
\subset  [t_0 - 2\pi\Theta, t_0 + 2 \pi T_e + 2\pi\Theta] $.
$ \forall \ T = 2 \pi \tau + 2 \pi \theta \in {\cal J} $
(i.e. $ | \theta| \leq \Theta$) formula \equ{OeT} becomes
\beq{aule}
\O_{\eta j} T
=
2\pi \Big(\O_j (\tau+\theta) -(\cQ\cR^{-1})_j
\langle \o(\tau+\theta)\rangle\Big)\qquad \forall\ 1\leq j\leq m.
\eeq
By \equ{disto},  $ \forall\, 1\leq i\leq n $
there exists a $k_i\in\integer$ such that
$\o_i\tau=k_i+\langle \o_i\tau\rangle$
with $|\langle \o_i\tau\rangle|\leq \d$;
moreover $|\o_i \theta|\leq \d/4$
by definition of $ \Theta $.  Hence we have
$ | \o_i ( \tau + \theta ) - k_i | \leq \d + \d/4 < 1 \slash 2 $
being $ \d \leq 1/4 $. So
$ \langle \o_i ( \tau + \theta ) \rangle $ $ = \o_i ( \tau + \theta) - k_i $
and
\beq{nienna}
|\langle \o_i(\tau+\theta)\rangle|\leq \frac{5}{4}\d\ ,\qquad
\forall\, 1\leq i\leq n\ .
\eeq

For $1\leq j\leq m,$
by \equ{distO} we have dist$(\O_j\tau,\integer)\geq d_0,$
hence, by the definition of $\Theta,$
 dist$(\O_j(\tau+\theta),\integer)\geq d_0-d_0/8.$
Collecting the previous inequalities and \equ{nienna},
from \equ{aule} we obtain
$$
{\rm dist} (\O_{\eta j}T,2\pi\integer)
=
2\pi {\rm dist} \Big(
\O_j (\tau+\theta) -(\cQ\cR^{-1})_j
\langle \o(\tau+\theta)\rangle,
\integer
\Big)
\geq
2\pi\left(\!
d_0-\frac{d_0}{8}-\frac{5 d_0}{8}\!\right)\!=\!\frac{\pi d_0}{2}
$$
by definition of $\d$. Finally, we recall
\equ{stimaM} to end the proof.
\hfill$\Box$

\lem{b} Let condition $(b)$ of Theorem \ref{theorem 1} hold,
namely $\a>0$. Let
$$
\theta:=\min_{1\leq i\leq n}\frac{1}{|\o_i|} \qquad {\rm and}
\qquad d_1:=\min \Big\{ \frac{\pi}{8m},\frac{\pi\a\theta}{2nm}
\Big\}.
$$
Then $\forall\, t_0\geq 0$ there exists an open set ${\cal A}
\subset  [t_0, t_0 + 4\pi\theta]$ of measure, at least,
$\pi\theta/n$, such that $\forall\, T\in {\cal A}$,
 $ \cM^{-1} := \cM^{-1}( I_0 (T), T ) $ satisfies
\beq{d1a} | \cM^{-1} | \leq \frac{2}{d_1}\ . \eeq \elem \noindent
{\bf Proof:} The function $ T \to \O_\eta ( I_0 ( T ) ) T $ is a
piecewise smooth function with discontinuities at the points $
T_{i,k} =\frac{2\pi}{\o_i}(k+\frac{1}{2})$, $1\leq i\leq n$,
$k\in\integer.$ Apart these points $ \O_\eta (I_0(T))T$ is
differentiable w.r.t. to $T$ and has constant derivative
$\O-\cQ\cR^{-1}\o=:\xi \in \real^m $.
By the definition of $\theta$, in every interval of the type
$\big( 2\pi\theta(h-1/2), 2\pi\theta(h+1/2) \big)$,
$h\in\integer$, fall at most $n-1$ points of discontinuity and,
hence, there exists an interval ${\cal J}_1\subseteq \big(
2\pi\theta(h-1/2), 2\pi\theta(h+1/2) \big)$ of length, at least,
$\D:=2\pi\theta/n,$ in which $\O_\eta (I_0 ( T ) ) T$ is smooth.
Hence, on ${\cal J}_1$, by \equ{OeT}, $\O_{\eta}T=x+\xi T$ for a
suitable $ x \in \real^m.$

For $ 1\leq j \leq m$, let define the sets
$$ {\cal B}_j:=
\left\{ \matrix{ \{\ T\in{\cal J}_1\ \ {\rm s.t.}\ \ {\rm
dist}(\O_{\eta j}T,2\pi\integer) >
\displaystyle{\frac{\pi}{8m}}\ \}
& \qquad{\rm if}\quad |\xi_j|\D\geq \pi  \cr
& \cr
\{\ T\in{\cal J}_1\ \ {\rm s.t.}\ \ {\rm dist}(\O_{\eta
j}T,2\pi\integer) > \displaystyle{\frac{\D}{4m}}\a\ \} &
\qquad{\rm if}\quad |\xi_j|\D< \pi } \right.\ .
$$
Remembering the definition $\a:=\min_{1\leq j\leq m}|\xi_j|,$
we note that  meas$({\cal B}_j)\geq \D(1-\frac{1}{2m})$,
$\forall 1 \leq j \leq m $.
Let ${\cal A}:=\cap_{j=1}^m {\cal B}_j,$ then ${\cal A}\subseteq {\cal J}_1$
and meas$({\cal J}_1)\geq \D/2.$ By construction $\forall\,
T\in{\cal A}\,$, dist$(\O_{\eta j}T,2\pi\integer)\geq d_1\,$,
$\forall\, 1\leq j\leq m.$ Finally  the Lemma follows from
\equ{stimaM}. \hfill$\Box$

\begin{remark}\label{linint}
{\rm Our non resonance conditions $(a)$--$(b)$ of Theorem
\ref{theorem 1} are sharp:
if $(a)$--$(b)$ are violated,
it is not possible in general to find
a period $ T $ is such a way that the matrix $\cM$
defined in \equ{Mmeno1}--\equ{OeT} is invertible
(clearly it must be $ m \geq 3 $ and $\hm \geq 1 $).
As an example, consider the Hamiltonian $ \cH_* = $
$$
\o \cdot \cI_* +
\O Z_*\cdot \overline Z_*+
\frac{1}{2}\,|\cI_*|^2+
\frac 1 p
\sum_{j=1}^{\hat m}
(a_{j1} {\cI_{*1}}+
a_{j2} {\cI_{*2}})
Z_{*j} \overline Z_{*j}+
\sum_{2|k|+|a + \overline a|\geq 6} R_{k, a, \overline a}^*(\f_*)
\cI_*^k Z_*^a \overline Z_*^{\overline a}\,,
$$
with $ m \geq 4 $ ($ \hm \leq m $) and
$$ a_{j1}:=
\left\{ \matrix{ 1& \qquad{\rm if}\quad 1\leq j\leq \hat m -1  \cr
0& \qquad{\rm if}\quad j= \hat m } \right.\,, \quad\qquad\qquad
a_{j2}:= \left\{ \matrix{ j& \qquad{\rm if}\quad 1\leq j\leq \hat
m -1\cr 1& \qquad{\rm if}\quad j= \hat m } \right.\,.$$
The Hamiltonian $\cH_* $ is in the form \equ{ham (cal H)}--\equ{Mj}. We
choose $ \hm $ so that $ p := \hat m-1$ is a prime integer with
$3 \leq p \leq m-1$
(we require $ p \neq 1,2 $ since $(\o, \O)$ must satisfy
the second order Melnikov non-resonance conditions \equ{CONDIZ NON RISON} and
$(\o, \O)$ are related by \equ{vio} below).
Let also  $\hat \o=(\o,\O_{\hm+1},\dots,
\O_{m})$ be any rationally independent vector in $\real^{\hat n}$
and define \beq{vio} \O_j := \frac{1}{p} (a_{j1}\o_1+a_{j2}\o_2
)\,, \quad \forall 1\leq j\leq \hm\,.\eeq
Note that, by \equ{vio}, the equations
\equ{Mj} are fulfilled with $M_j=p$, $\forall 1 \leq j \leq \hm$.
Hence $M_j = p = \hm -1 < \hm $ and
condition $(a)-(iii)$ is violated. After the
rescaling in \equ{resca}, the new Hamiltonian
is $ \cH = \o \cdot \cI + \O Z \cdot \overline Z + \eta^2
( \frac{1}{2} | \cI |^2 +  \cQ I \cdot Z \bZ ) + O(\eta^4)$ where
$ \cQ \in {\rm Mat}(m \times n, \real) $ is defined by
$$\cQ_{ji}=\left\{
\matrix{
{a_{ji}}/{p}& \qquad{\rm if}\quad 1\leq j\leq \hat m \,,\quad i=1,2\,,\cr
0 & \qquad{\rm elsewhere}\,.
}
\right.$$
Since the Hamiltonian $\cH $ does not contain
terms of the form \equ{pb1}--\equ{pb22}-\equ{pb21}-\equ{pb3},
$\cH $
is yet in the normal form \equ{ham (****)} (i.e. \equ{ham (3)})
with  $\cR={\bf 1}_n$.
Therefore
$ ( \O - \cQ \cR^{-1}\o )_j = ( \O - \cQ \o )_j = 0 $,
$ \forall 1 \leq j\leq \hm $, by \equ{vio},
and also condition $(b)$ is violated. Now,
whatever we choose $ T \in \real $, $I_0 \in \real^n $,
$ n \in \integer^n $ such that
$ \o + \eta^2 \cR I_0 = 2 \pi n / T $, i.e. \equ{otilde} is satisfied,
we get, substituing $ I_0 = \eta^{-2}  \cR^{-1} ((2 \pi n / T) - \o) $
into \equ{omega eta},
$ \O_\eta (I_0) T = (\O - \cQ \o) T + 2 \pi \cQ  n $.
Hence, $\forall 1 \leq j \leq m $
$$ \O_{\eta j}T= \O_{\eta j}(I_0) T = 2 \pi (\cQ_{j1} n_1 + \cQ_{j2}n_2)\,.$$
We claim that, for all $ n_i \in \integer$,
at least one $ \O_{\eta j}T$ is an integer
multiple of $2\pi$ and hence the matrix
$\cM := {\bf 1}_m - e^{\ii \O_\eta  T} $ has a zero eigenvalue.
In fact,
$$\cQ_{j1} n_1 + \cQ_{j2}n_2= \left\{
\matrix{
{(n_1+n_2 j)}/{p}& \qquad{\rm if}\quad 1\leq j\leq p \,,\cr
n_2/p & \qquad{\rm if}\quad j=\hm\,.
}\right.
$$
Thus, if $n_2$ is a multiple of $p$, then
$\O_{\eta \hm}T\in 2\pi \integer$.
Otherwise, let $j^*$ be the unique solution of the linear
congruence $ n_2 j^* \equiv -n_1 $ ${\rm mod} p $
(recall that  $p$ is prime). In this case
$\O_{\eta j^*}T\in 2\pi \integer$. In both cases $\cM$ is not invertible.}
\end{remark}

In the next Lemma we prove the existence of suitable pseudo
$T$-periodic solutions of the
Hamiltonian system \equ{equation of motion} close to the
manifold $ { \cal P } $.
Roughly, by the ``twist condition'' det$\cR \neq 0$ and the nonresonance
property \equ{Mmeno1}, the manifold ${\cal P}$ is ``non-degenerate'',
i.e. the only $T$-periodic solutions of $H_{int} $, close to ${\cal P}$,
are the set ${\cal P}$. This implies, by the Contraction Mapping Theorem,
the existence of a manifold of pseudo $T$-periodic solutions $\z_{\phi_0}$,
close to ${\cal P}$, diffeomorphic to $ \torus^n $.
$\z_{\phi_0}$ are solutions of
\equ{Hequat} for all $ t \in (0, T) $
and satisfy $\phi (T) = \phi(0) = \phi_0 $,
$ z(T) = z(0) $ but it may happen that $ I(T)  \neq I(0)$.

\begin{lemma}  \label{pseudo (1)}
Assume that Condition $(a)$ or $(b)$ of Theorem
\ref{theorem 1} holds. Then,
there exist $ \eta_0 $, $ C_0, C > 0 $, such that: $ \forall
\eta \in (0, \eta_0 ] $, there exist an open set
$\cA_\eta\subseteq [\frac{1}{\eta^2},\frac{1}{\eta^2}+C_0]$ of
measure greater than $1/C_0$ such that $\forall T\in\cA_\eta$ and
$ \forall \phi_0 \in \torus^n $ there exists a unique function
$$
\z_{\phi_0} := \Big( I_{\phi_0}, \phi_{\phi_0}, z_{\phi_0} \Big)
\in  C \Big( [0,T] , \, \real^n \times \torus^n \times \complex^{m} \Big) \cap
C^1 \Big( (0,T), \, \real^n \times \torus^n \times \complex^{m} \Big),
$$
smooth in $ \phi_0 $, such that
\begin{itemize}
\item ($i$)
$ \z_{\phi_0}(t) $ solves \equ{Hequat} for all $t \in (0,T) $;
\item ($ii$)
$\z_{\phi_0}$ satisfies
$ \phi_{\phi_0}(0) = \phi_{\phi_0}(T) = \phi_0 $,
$ z_{\phi_0}(0) = z_{\phi_0}(T) $;
\item ($iii$)
$ \sup_{t\in [0,T] }\Big( |I_{\phi_0}(t) -I_0|+
|\phi_{\phi_0}(t)-\phi_0-\widetilde \o t| +
|z_{\phi_0}(t)| \Big)\leq C \eta^2\, . $
\end{itemize}
where $I_0:=I_0(T)$ is defined in \equ{I0} and
$\widetilde\o:=\widetilde\o(T):=\o+\eta^2 \cR I_0.$
\end{lemma}

\noindent {\bf Proof:}
In order to define the set
$\cA_\eta \subseteq [\frac{1}{\eta^2},\frac{1}{\eta^2}+C_0]$
of ``non-resonant'' periods $ T $,
we use Lemmata \ref{a} or \ref{b}
(according whether condition $(a)$ or $(b)$ of Theorem
\ref{theorem 1} holds) with
$t_0:=\eta^{-2}$. 
$\forall T \in \cA_\eta $
we look for a solution $ \z_{\phi_0} $ of
\equ{equation of motion} (i.e. \equ{Hequat}) of the form
\beq{variation} I_{\phi_0} (t)= I_0 +\eta J(t), \quad
\phi_{\phi_0} (t)  = \phi_0 +\widetilde\o t +\eta \psi(t), \quad
z_{\phi_0} (t) = \eta w(t), \eeq for suitable functions ($J$,
$\psi$, $w$) : $[0,T] \to \real^n \times \torus^n \times
\complex^{m} $  \ satisfying $ \psi (0) = \psi (T) = 0 $ and  $
w(0) = w(T) $. The condition $ \psi (0) = 0 $ is a ``tranversality
condition'': we impose the ``correction'' $(J, \psi, w )$ to
belong to a supplementary linear space to the tangent space of the
unperturbed manifold $ { \cal P} $ (see [BB] for a discussion of
the different supplementary spaces). The functions ($J$, $\psi$,
$w$) must satisfy the system
\begin{eqnarray}
\dot J &=& O(\eta^3)\, , \label{Jdot}
\\ \dot \psi  -\eta^2 \cR J &=& O(\eta^3)\, , \label{psidot} \\
\dot w_j - \ii\O_{\eta j} w_j &=&
\ii\eta^2 \sum_{|a+\ba|=3} R^*_{0,a,\ba}(\phi_0+\widetilde\o t)
\ba_j w^a \bw^{\ba-e_j}+
O(\eta^3), \quad j = 1, \ldots, m  \label{eqofmo}
\end{eqnarray}
where $  \O_\eta = \O_\eta (I_0(T)) \in{\rm Mat}(m\times m, \real)$ is
the real diagonal matrix defined in \equ{omega eta}.

In order to find $(J,\psi,w) $
we use Lemmata \ref{lemma 0} and \ref{lemma 1}.
In connection with the notation of
Lemma \ref{lemma 1} we have here $ M = \eta^2 \cR $ and
$ \cM = {\bf 1}_m - e^{\ii \O_\eta T} $.
Let $ L $ denote the corresponding Green operator in \equ{LGreen}
and $ P := P (J, \psi, w;\phi_0) $ the right hand side of
\equ{Jdot},\equ{psidot},\equ{eqofmo}.
It is sufficient to find a fixed point $(J, \psi, w) \in X $,
space defined in \equ{spazY}, of
\beq{pfisso}
(J,\psi,w)=L \big( P( J,\psi,w;\phi_0 )\big)\ .
\eeq
If $ T $ is obtained in Lemma \ref{a}, resp. Lemma \ref{b}, we have,
by  \equ{d1a}, resp. \equ{d2b}, that
$ | \cM^{-1} | \leq 2 \slash d_1 $,  resp.
$ | \cM^{-1} | \leq 4 \slash \pi d_0 $.
Moreover $ | M | = O(\eta^2 ) $,
$ | e^{\ii \O_\eta} | = |{\rm diag}\{ e^{\ii \O_{\eta j}} \} | \leq m $,
$ T = O ( \eta^{-2} ) $ and, hence, from \equ{eq *}, we get
the estimate on the Green operator
$ | L | \leq C' \eta^{-2} $ for some $ C' > 0 $.
Moreover $ P(0;\phi_0)=O ( \eta^3) $. Then, for a constant
$ C $ large enough we have
$ \d_0 := C \eta \geq  2|L(P(0))|$.  Finally,
since $\sup_{B_{\d_0}}|DP|=O(\eta^3),$
for $ \eta $ small enough, we verify also \equ{DP}. Applying
Lemma \ref{lemma 0} in the ball $B_{\d_0}$ we
prove the existence of a solution  $ \z_{\phi_0} $
satisfying $ || \z_{\phi_0} || \leq C \eta$, and hence, by
\equ{variation}, we get the estimate ($iii$).

Since $ \z_{\phi_0} \in X $ solves the integral
system \equ{pfisso}, $ \z_{\phi_0} (t) $ is actually $ C^1 $
for all $ t \in (0,T) $ and solves \equ{Hequat}.

Finally, since the operator  $\phi_0 \to L P(J,\psi,w;\phi_0 )$
is smooth, by the Implicit Function Theorem, we also deduce that
the function $ \phi_0 \to \z_{\phi_0} $ is $ C^1 $. \hfill$\Box$

\begin{remark}\label{remdt}{\rm
The terms \equ{pb1}, \equ{pb22}, \equ{pb21}, \equ{pb3}, that
have been ``averaged out'' in Proposition 3.1, could not have been
treated perturbatively in the previous Lemma.
Indeed \equ{pb1} (resp. \equ{pb3}) would have given rise
to a term $ O (1)$ (resp. $O(\eta^2 )$) in \equ{eqofmo};
\equ{pb22} to a term $O(\eta )$ in \equ{Jdot};
and \equ{pb21} to a term
$ O(\eta^2 )$ in \equ{eqofmo}.
Roughly speaking, since the period $ T \geq 1/\eta^2 $, only terms
of magnitude $o(\eta^2)$ can be dealt as perturbations.
}\end{remark}

\begin{remark}\label{twist}{\rm
The invertibility of the twist matrix $ \cR $
has been used in two ways:  to modulate the frequency $ \o $
yielding the existence of nearby resonant frequencies,
see \equ{otilde}, and to imply that the manifold
${\cal P}$ is non-degenerate
(together with the non-resonance property \equ{Mmeno1}).
Weaker conditions could in principle be used 
(if the elliptic variables were absent think to [CZ1]). 
However, in this case, when
$ \cR $ is degenerate one is forced to study the higher order terms
in the normal form of Proposition
\ref{averaging} making the perturbative analysis
much more difficult. Moreover in the three body problem
application we will prove that $\cR$ is invertible, see Lemma
\ref{nonde3}. 
}\end{remark}

\begin{remark}\label{intera}
{\rm
When $m=0$, namely in the case of {\sl maximal} tori,
the twist matrix $ \cR $ reduces simply to the
matrix of the second derivatives of $\cH_* $ w.r.t.
the action variables $\cI_* $. This is not the case when
$ m \neq 0 $.
Consider
the Hamiltonian $ \cH_* = $ $ \o \cI_* + $ $ \O Z_* \bZ_* + $ $
\cI_* ( Z_* + \bZ_* ) + $ $  \cI_*^2 / \O $ where $(\cI_*, \f_*)
\in \real \times \torus $, $ (Z_*, \bZ_*)\in \complex^2 $ and $\O>0$.
The twist matrix in \equ{R2} turns out to be $ \cR = 0$ although
the second derivative of $ \cH_*  $ w.r.t. $ \cI_* $ equals to
$ 2 / \O \neq 0 $.
Heuristically, for {\sl lower} dimensional elliptic tori,
the twist matrix $ \cR $ 
takes into account,
with the second addendum in the r.h.s. of \equ{R2},
the interaction between the torus variables $(\cI_*, \f_* )$ and the
elliptic variables $ Z_* $.
Finally observe that the symplectic  map $ \cI_* = I,$
$\f_* = \phi+  (z-\overline z)/(\ii\O),$ $ Z_* = z- I/\O $
(which can be found through the averaging procedure of
section \ref{Sec:NorForm})
transforms $ \cH $ into the
{\it isochronous} Hamiltonian $ H= \o I + \O z \bz $.}
\end{remark}

\subsection{The variational principle}\label{pseudo (2)}

We now define the ``reduced Hamiltonian action functional''
$ \cE : \torus^n \to \real $
as the Hamiltonian action functional ${\cal A}$
evaluated on the pseudo $T$-periodic solutions $ \z_{\phi_0} $
obtained in Lemma \ref{pseudo (1)}, namely\footnote{$\cE (\phi_0 )
\in \real $ since $\int_0^T $ Im
$(\ii z_{\phi_0}\dot {\overline z}_{\phi_0}) \ dt =$
$\int_0^T $ Re $(z_{\phi_0} \dot {\overline z}_{\phi_0}) \ dt $ $= \int_0^T
\frac{1}{2} \frac{d}{dt} (z_{\phi_0} \overline z_{\phi_0}) \ dt = 0$.}
\beq{Reduact}
\cE ( \phi_0 ) := {\cal A} ( I_{\phi_0}, \phi_{\phi_0}, z_{\phi_0},
 {\overline z}_{\phi_0})
= \int_0^T I_{\phi_0}\cdot \dot\phi_{\phi_0}
+ \ii z_{\phi_0} \cdot \dot {\overline z}_{\phi_0} -
H( I_{\phi_0}, \phi_{\phi_0}, z_{\phi_0},
 {\overline z}_{\phi_0} )\,dt\,.
\eeq
By Lemma \ref{pseudo (1)} the reduced action functional
$\cE$ is smooth in $ \phi_0 $.

Critical points of $ \cE : \torus^n \to \real $ give
rise to $T$-periodic solutions
of the Hamiltonian system \equ{Hequat}, according
to the following Lemma:

\lem{LEMMA M}
\noindent
$ \partial_{\phi_0} \cE ( \phi_0 ) = I_{\phi_0}(T) -
I_{\phi_0} ( 0 ) $. Hence, if $\phi^\star \in \torus^n $
is a critical point of $\cE$,
then $ \z_{\phi^\star}(t) $ can be extended
to a $ T $-periodic solution of the Hamiltonian system \equ{Hequat}.
\elem

\noindent {\bf Proof:}
Differentiating w.r.t. $ \phi_0 $ in \equ{Reduact} we get
\begin{eqnarray*}
\partial_{\phi_0} {\cal E}( \phi_0 ) & = &
\int_0^T \partial_{\phi_0} I_{\phi_0} \cdot \dot \phi_{\phi_0} +
I_{\phi_0} \cdot \partial_{\phi_0} {\dot \phi}_{\phi_0}
+ \ii \partial_{\phi_0} z_{\phi_0} {\dot {\overline z}}_{\phi_0}
+ \ii z_{\phi_0} \partial_{\phi_0} {\dot {\overline z}}_{\phi_0} \\
& - & \partial_I H \partial_{\phi_0} I_{\phi_0} -
\partial_\phi H \partial_{\phi_0} \phi_{\phi_0} -
\partial_z H \partial_{\phi_0} z_{\phi_0} -
\partial_\bz H \partial_{\phi_0} \bz_{\phi_0} \ dt
\\
& = & I_{\phi_0} (T) \partial_{\phi_0} \phi_{\phi_0} (T) -
 I_{\phi_0} (0) \partial_{\phi_0} \phi_{\phi_0} (0) +
\ii \partial_{\phi_0} {\overline z}_{\phi_0} (T) z_{\phi_0}(T) -
\ii \partial_{\phi_0} {\overline z}_{\phi_0} (0) z_{\phi_0}(0),
\end{eqnarray*}
by an integration by parts, and since
$\z_{\phi_0}(t) $ satisfies the Hamilton's equations \equ{Hequat}
in $(0,T)$. Moreover, since $ \forall \phi_0 \in \torus^n $,
$ \phi_{\phi_0} (T) = \phi_{\phi_0} (0) = \phi_0 $
and $ {\overline z}_{\phi_0}(T) = {\overline z}_{\phi_0}(0) $,
deriving w.r.t. $ \phi_0 $, we get
$ \partial_{\phi_0} \phi_{\phi_0} (T) =
\partial_{\phi_0} \phi_{\phi_0}(0) = {\bf 1}_n $ and
$ \partial_{\phi_0} {\overline z}_{\phi_0}(T) =
\partial_{\phi_0} {\overline z}_{\phi_0}(0) $.
This implies $ \partial_{\phi_0} \cE ( \phi_0 ) = I_{\phi_0}(T) -
I_{\phi_0} ( 0 ) $. Hence, if $ \phi^\star \in \torus^n $ is a
critical point of $ \cE $, $  \z_{\phi^\star}(T) =
\z_{\phi^\star}(0)$ and we deduce that $ \z_{\phi^\star}(t) $ can be
extended to  a $ T $-periodic solution of the Hamiltonian system
\equ{Hequat}. \hfill$\Box$\giu

The following Lemma, which is a consequence of the
autonomy of the Hamiltonian $ H $, holds.

\lem{invar}
$\forall \phi_0 \in \torus^n $ there exists
$ V ( \phi_0, \eta ) \in \real^n $
with $ V( \phi_0, \eta ) = \widetilde{\o } + O(\eta^3 ) $ such that
\beq{Finva}
V( \phi_0, \eta ) \cdot  \partial_{\phi_0} \cE ( \phi_0 ) = 0,
\qquad \forall \phi_0 \in \torus^n.
\eeq
\elem

\noindent {\bf Proof:}
Since $ \z_{\phi_0} = ( I_{\phi_0}, \phi_{\phi_0}, z_{\phi_0} ) $
satisfies the Hamiltonian system \equ{Hequat} in $ ( 0, T) $
and $ \phi_{\phi_0} (T) =  \phi_{\phi_0} (0) = \phi_0 $,
$ z_{\phi_0}(T) = z_{\phi_0}(0) $,
then
\beq{Consen}
H (I_{\phi_0}(T), \phi_0 , z_{\phi_0}(0), {\overline z}_{\phi_0}(0) ) =
H (I_{\phi_0}(0), \phi_0 , z_{\phi_0}(0), {\overline z}_{\phi_0}(0)).
\eeq
By the mean value Theorem there exists
$ \xi_{\phi_0} $ in the segment between $[I_{\phi_0}(T), I_{\phi_0}(0)]$
such that
\beq{MeanV}
\partial_I H ( \xi_{\phi_0}, \phi_0 , z_{\phi_0}(0),
{\overline z}_{\phi_0}(0)) \cdot [I_{\phi_0}(T) -  I_{\phi_0}(0) ] = 0.
\eeq
Define $ V(\phi_0, \eta) := \partial_I H ( \xi_{\phi_0}, \phi_0,
z_{\phi_0}(0), {\overline z}_{\phi_0}(0) ) $.
By Lemma \ref{LEMMA M} and \equ{MeanV}, formula \equ{Finva} follows.
Finally, since $ I_{\phi_0}(T) - I_0 =
I_{\phi_0}(0) - I_0 = O(\eta^2) $,
$ z_{\phi_0} ( 0 ) = O ( \eta^2 ) $ and
$ \partial_I H =$ $\o +$ $\eta^2 \cR I +$ $\eta^2 {\cal Q}^T z {\overline z} +
$ $O(\eta^3 ) $,
we deduce the estimate
$ V(\phi_0, \eta) = \widetilde{\o} + O(\eta^3) $.
\hfill$\Box$\giu

\giu

\noindent
{\bf Proof of Theorem \ref{theorem 1}:} \label{proof theorem 1}
By Lemma \ref{LEMMA M} the absolute minimum (and maximum)
$ \phi^\star \in \torus ^n $ gives rise to a $T$-periodic solution
$\z_{\phi^\star} $ of \equ{Hequat}.
However one expects the existence of at least
$n$ geometrically distinct $T$-periodic solutions of
\equ{Hequat}, i.e. solutions
not obtained one from each other simply by time-translations.

In order find multiple geometrically distinct periodic solutions
of \equ{Hequat}
we restrict the reduced action function ${\cal E}$
to the plane $ E := [\widetilde \o ]^{\bot} $ orthogonal
to the periodic flow
$ \widetilde \o = (1/T) 2 \pi k $ with $ k \in \integer^n $.
The set $ \integer^n \cap E $ is a lattice of $ E $,
(see for example Lemma 8.2 of [BBB]) and
hence $ \cE $ can be defined on the quotient space
$\Gamma := E \slash (\integer^n \cap E) \sim \torus^{n-1}$.

A critical point $ \phi^\star $ of $ \cE : \Gamma \to \real $
is a critical point of $ \cE : \torus^n \to \real $. Indeed,
since the tangent space
$ T_{\phi_0} \Gamma = [\widetilde \o ]^{\bot} $, then
\beq{critG}
\partial_{\phi_0} {\cal E}( \phi^\star )
= \lambda (\phi^\star ) \widetilde{\o}
\eeq
for some Lagrange multiplier $ \lambda ( \phi^\star ) \in \real $.
\equ{critG} and Lemma \ref{invar} imply that
$ \lambda  ( \phi^\star ) (| \widetilde{\o} |^2 + O(\eta^3 ))
= 0 $
and so, for $ \eta $ small, $ \lambda ( \phi^\star ) = 0 $.

By the Lusternik-Schnirelman category theory, see for example
[Am], since cat$\Gamma$ = cat$\torus^{n-1}$ = $n$, we can define
the $n$ min-max critical values $c_1 \leq c_2 \leq \ldots \leq c_n
$ for the reduced action functional $ \cE_{| \Gamma } $. Let $
\phi_1^\star, \ldots, \phi_n^\star \in \Gamma $ be $n$
corresponding  critical points. If the critical levels $ c_i $ are
all distinct, the $T$-periodic solutions $ \z_{\phi_i^\star} (t)
=$ $( I_0, \phi_i^\star + \widetilde \o t , 0) + O ( \eta^2 )$ of
\equ{Hequat} are geometrically distinct, since their actions
${\cal A}( \z_{\phi_i^\star} ) = {\cal E} ( \phi_i^\star ) = c_i $
are different. On the other hand, if some min-max critical level $
c_i $ coincide, then $ \cE_{| \Gamma } $ possesses infinitely many
critical points. However {\it not} all the corresponding
$T$-periodic solutions
of \equ{Hequat} are necessarily 
geometrically distinct, since two different critical points
could belong to the same orbit. In any case, since one
periodic solution can cross $ \Gamma $ at most
a finite number of times, the existence of infinitely many
geometrically distinct orbits of \equ{Hequat} follows\footnote{
Non-degenerate critical points of the
Poincar\'e-Melnikov primitive
$ \phi_0 \to \int_0^T (H - H_{int}) (I_0, \phi_0 + \widetilde \o t, 0)\ dt $
could be continued to solutions of \equ{Hequat}. However,
it is very difficult to compute this function,
its critical points and check
whether they are non-degenerate
(if ever true).}.

Finally, under the inverse transformation of $ \Phi $, defined in
Proposition \ref{averaging}, $\eta$-close to the identity
and the inverse transformation of \equ{resca},
we find the $T$-periodic solutions $\z_\eta $
of Theorem \ref{theorem 1} satisfying the estimates ($i$)-($ii$).

The statement on the minimal period descends from the following
Proposition. \hfill$\Box$

\begin{proposition}\label{minimo}
Let $ \z(t)=(I(t),\ph(t),z(t))$ be a $T$-periodic solution of
\equ{equation of motion} with
$C \eta^{-2}$ $ \geq $ $ T \geq $ $T_0 $
for some $T_0 $ large, independent of $\eta $. Then
the minimal period $T_{min}$ of $\z $
satisfies $T_{min}\geq {\rm const}$  $
T^{1/(\tau+1)}.$
\end{proposition}
{\bf Proof:} Since $ \z $ is $ T $-periodic,
$ \ph(T) - \ph(0) = 2\pi k$ for some $ k \in \integer^n.$
We claim that $ T_{min} \geq T / g $
where $ g:=\gcd (k_1,\dots k_n)$. Indeed
$ T_{min} \geq T_{min}^\ph $ where $T_{min}^\ph$
is the minimal period of 
$\ph (t)$. Moreover
$ T_{min}^\ph  = T / n $ for some integer $ n \geq 1 $ and
$ \ph ( T_{min} ) - \ph(0) = 2 \pi \widetilde{k} $
for some $ \widetilde{k} \in \integer^n$.
It follows that $ \widetilde{k} n = k $ and so $ n \leq g $,
proving the claim.

By the second equation in \equ{equation of motion}
we also have $\ph(T)-\ph(0)=\o T+O(\eta^2 T)$ and then,
letting $ \hat k:=k/g\in\integer^n$,
$  2 \pi g \hat k = \o T + O( \eta^2 T ) $. We deduce that
\beq{Tmin}
T_{min} \geq \frac{T}{g} \geq
\frac{2\pi}{|\o |} |\hat k| - O \Big( \frac{\eta^2 T }{|\o |g }\Big)
\geq \frac{2\pi}{|\o |} |\hat k| - O \Big( \frac{C}{|\o |}\Big),
\eeq
since $ \eta^2 T \leq C $.
Let now choose
$ h\in\integer^n\setminus \{ 0\}$ with $|h|\leq |\hat k |$ such
that $\hat k\cdot h=0$.
Multiplying by $h $ we get $ 0 = 2 \pi g \hat k \cdot h
= \o \cdot h T + O( \eta^2 T |h| ) $
ans so 
$\o \cdot h  =O(\eta^2 |h|).$ Using the diophantinity of
$\o$ assumed in \equ{CONDIZ NON RISON}, namely $|\o\cdot h|\geq \g
(1+|h|)^{-\tau}$, $\forall h \in \integer^n$,
we obtain $|h|^{\tau+1}\geq {\rm const}$
$\eta^{-2}\geq$ $  {\rm const} \ T$. Hence
$|\hat k|^{\tau+1}\geq | h |^{\tau+1} \geq $
${\rm const} \ T$ and $|\hat k|\geq {\rm const} \ T^{1/(\tau+1)}.$
The proposition follows from \equ{Tmin} for $T \geq T_0 $ large enough.
\hfill$\Box$

\begin{remark}\label{perio}
{\rm
The periodic orbits $ \z_\eta  $ found in
Theorem \ref{theorem 1} are $ O(\eta^2)$ close to the torus
${\cal T} $ since $ I_0 $, defined in \equ{I0}, satisfies $ I_0 = O(1) $.
We could also try to find periodic solutions when
e.g. $ I_0 = O(1 \slash \eta) $
and so $ T= O ( 1 \slash \eta )$. However the terms
$ O(I^k) = O( \eta^{-k} ) $ would be more difficult to control.
}\end{remark}


\section{The planetary spatial three-body problem}\label{physics}

In this section we will prove the existence of periodic orbits
of the planetary non planar three body problem,
with  ``small eccentricities'' and ``small mutual inclinations'',
accumulating onto two-dimensional elliptic invariant tori. We
first discuss the classical Hamiltonian formulation of this problem
which dates back
to Delaunay and Poincar\'e. For a detailed treatment see the Appendix of [BCV].
\giu

The three massive points (``bodies") $ P_0, \,P_1\,,P_2 $, with
masses $ m_0,\,m_1\,,m_2 $, interact one each other through
gravity (with constant of gravitation 1). Assume that the masses
of the bodies satisfy, for some $0<\overline\k\leq 1$, \beq{masse}
\overline\k {\e}\leq\frac{m_1}{m_0}\,,\, \frac{m_2}{m_0} \leq \e
\le 1 \, . \eeq The number $\e>0$ is regarded as a small
parameter: the point $P_0$ represents ``the star'' and the points
$P_1$ and $P_2$ ``the planets''.

We now recall the classical definition of the
``osculating ellipses'' (at time $t_0$)
of the two-body problems associated to the planets $P_i$ ($i=1,2$)
and the star $P_0$.
Let $u^{(0)}$ and $u^{(i)}$ denote the coordinates of the points
$P_0$ and $P_i$ at time $t_0$ and let $\dot u^{(0)}$
and $\dot u^{(i)}$ denote the
respective velocities.  The ``osculating plane" is defined
as the plane spanned by $(u^{(i)}-u^{(0)})$ and
$(\dot u^{(i)}-\dot u^{(0)})$; the ``osculating ellipse'' is defined as
the Keplerian ellipse (lying on the osculating
plane) defined by the Kepler solution, with initial data
$(u^{(0)},u^{(i)})$ and $(\dot u^{(i)}-\dot u^{(0)})$,
of the two-body  problem $(P_0,P_i)$ obtained disregarding
(for $t \ge t_0$) the third body $P_j$ ($j \neq i$). 

We assume that the eccentricities of such ellipses are small
and {that the intersection angle} between the two planes
containing the two osculating ellipses (usually referred as
``mutual inclination") is also small. It is customary in celestial
mechanics to denote the {\it major semi-axes}
of such ellipses by $ a_i $ and their eccentricities by $e_i$.
Let 
$$ \Lambda^*_i := \k_i^*
\sqrt{a_i}\ ,\qquad \k^*_i :=
\frac{m_i}{\e}\,\frac{1}{\sqrt{m_0 (m_0+m_i)}}\ ,$$
($\k^*_i$ is a dimensionless  constant satisfying $\frac{\overline \k}{\sqrt{2}}<\k^*_i<1$).
Since we are interested in small eccentricities,
collisions are avoided by requiring that
the  major semi-axes $ a_i = a_i(\Lambda^*):=(\Lambda_i^*/\k_i^*)^2$,
$i=1, 2 $, are different, and different from zero.
We, therefore, fix, once and for all,
\beq{aminmax}
0<a_{\rm min}<a_{\rm max}
\eeq
and, from now on,
we shall consider (attaching the index 1 to the
``inner planet") values of $\Lambda^*$ in the  set
\beq{compactdomain}
{\cal L} :=
\{\Lambda^*\in\real^2: a_{\rm min}\le a_1<a_2\le a_{\rm max}\}\ .
\eeq
\giu

The  following classical result follows from the
Delaunay-Poincar\'e theory, see [BCV].
\\[1mm]
\noindent
{\bf{Theorem (Delaunay-Poincar\'e)}} 
{\it Fix $\Lambda_0^* \in {\cal L}$ . There exists a symplectic
set of varia\-bles\footnote{With symplectic form $d I_0 \wedge d
\varphi_0 + dp_0 \wedge dq_0 $.}
 $$(I_0,\f_0,p_0,q_0)\in \real^2\times \torus^2
\times \real^2\times\real^2$$
where $ I_0 = \Lambda^* \in  \cI$ and
$ \cI  \subset {\cal L} \subset \real^2 $ is a suitable
two-dimensional open cube centered at $\Lambda_0^*$,
such that the Hamiltonian of the spatial planetary three-body problem
takes the form
\begin{equation}\label{hamform}
\cH_0(I_0,\f_0,p_0,q_0)= h_0(I_0)+f_0(I_0,\f_0,p_0,q_0)\,,
\end{equation}
with
\beqa{specham}
&& h_0:=-\frac{1}{2}\sum_{i=1}^2 \frac{\kappa_i}{I_{0i}^2}\,,\quad
\kappa_i := \Big( \frac{m_i}{\e} \Big)^3 \frac{1}{m_0^2(m_0 + m_i)}, \quad
\Big( \frac{{\bar \kappa}^3}{2} < \kappa_i < 1 \Big),
\nonumber\\
&& f_0:=\e f_1 (I_0,p_0,q_0)+\e f_2 (I_0,\f_0,p_0,q_0)\,,\nonumber\\
&& f_1 := f_{1,0}(I_0) + \sum_{j=1}^2 \overline\O_j (I_0)
(p_{0j}^2+q_{0j}^2) +\widetilde{f}_1 (I_0,p_0,q_0)\ ,
\nonumber\\
&&
\int_{\torus^2} f_2 \,d\f_0=0\ ,\qquad
\sup_{ \cI_{\s_0}^2 }  |
\widetilde{f}_1|\leq {\rm const} |(p_0,q_0)|^4\ .
\eeqa
$f_1$ and $f_2$
are real-analytic and uniformly bounded on\footnote{
$ D_{\r}^d \subset \complex^d $ denotes the open complex
ball of radius $ \r $ centered at the origin.}
\beq{ssr}
{\cal D}_0:=  {\cal I}_{\s_0}^2
\times \torus^2_{s_0}\times  D_{\r_0}^4\subset
\complex^8\ ,
\eeq
where
$ \s_0, s_0, \r_0 $ are suitable positive numbers.
Moreover $ \widetilde f_1 $ is even in $(p_0,q_0)$ and
\beqa{eigenst}
&&\inf_{I_0 \in \cI} \overline\O_j (I_0 ) >
\inf_{I_0 \in \cI_{\s_0}^2 } |\overline\O_j (I_0 )|>\,{\rm const}\,>0\
,\nonumber\\
&& \inf_{I_0 \in \cI} \Big(\overline\O_2(I_0 )-\overline\O_1(I_0 ) \Big)
> \inf_{I_0 \in \cI_{\s_0}^2}
|\overline\O_2 (I_0 ) -\overline\O_1 (I_0) |>\,{\rm const}\,>0
\,.
\eeqa}

We underline that the actions $ I_0 $ are simply $ I_0 =
\Lambda^* $ and we refer to [BCV]-Theorem 1.1  
for the complete expressions
of the other Delaunay-Poincar\'e variables $(\f_0, p_0, q_0)$.

We remark that the eccentricities $ e_i $ are estimated 
as\footnote{See (1.5),(C.10),(C.15),(C.17) of [BCV].} 
\beq{ecce}
c^-_{{\cal I}} |(p_0,q_0)| \leq e_i \leq c^+_{{\cal I}}
|(p_0,q_0)| \ , \qquad i=1,2
\eeq 
for two suitable constants 
$ 0 < c^-_{{\cal I}} < c^+_{{\cal I}}$ 
(depending only on ${\cal I}$). 

\giu

The Hamiltonian ${\cal H}_0$ in \equ{hamform} describes a
nearly-integrable, {\it properly degenerate}, system: the integrable
Hamiltonian $ h_0 $ depends only on the two action variables
$(I_{01}, I_{02})$. 
Such degeneracies are a typical feature of
problems arising in celestial mechanics
and the source of the main difficulties
(the application of standard KAM theory for finding
maximal tori require the Hamiltonian to be non-degenerate).
The frequency vector $ \nabla h_0 ( I_0 ) $
is independent of $ \e $, and then the conjugated angles
$\varphi_0 $ may be regarded as ``fast angles'' and, in
``first approximation'' the ${\cal H}_0 $-motions are governed by the
averaged Hamiltonian  $ h_0 + \e f_1 $.
By \equ{specham}, for any fixed $I_0$,
$\{ \f_0 \in \torus^2 \} \times \{
p_0 = q_0 = 0 \} $ is an elliptic invariant torus for
the averaged Hamiltonian $ h_0 + \e f_1 $ run by the linear flow
$ \f_0 \to \f_0 + \nabla (h_0 + \e f_{1,0})(I_0) \ t $.
These are the quasi-periodic
motions that will persist, for $ \e $ small, and $ I_0 $ in a nearly
full (two-dimensional) measure set, as proved in [BCV], see also
Theorem \ref{quasip}.

\giu

In the next Proposition, through an appropiate averaging procedure,
The Hamiltonian \equ{hamform} is casted into a suitable normal
form.

\begin{proposition}\label{Hstorta1}
Fix $ N \in \natural^+ $.
There exists a $ O ( \sqrt{\e} )$-close to the identity,
real analytic,
symplectic change of variables
$ ( J, \psi, z , \bz ) \in
U_{\widetilde{r}}^2 \times \torus^2_{\widetilde{s}}
\times D^4_{\widetilde{\r}} \to $
$ (I_0, \f_0, p_0, q_0 ) \in {\cal D}_0 \subset \complex^8 $
transforming the Hamiltonian \equ{hamform} into the
real analytic Hamiltonian
\begin{equation}\label{Hstorta}
\cH ( J, \psi, z, \overline z) =
h_\e(J) + \e \widetilde\O(J)z \overline z
+\e g(J,z,\overline z; \e)+\e^N f(J,\psi,z,\overline z; \e )
\end{equation}
where $ h_\e:= h_0 + O(\e) $,
\begin{equation}\label{intornini}
\widetilde{r}={\rm const}\,\sqrt\e, \qquad
\widetilde{s}={\rm const}, \qquad
\widetilde{\r} = {\rm const}\, 
\ ,
\end{equation}
\begin{equation}\label{U}
U:= \left\{ J\in \cI\ \ :\ \ |h'_0(J)\cdot\ell|\geq\a_0\ ,
\forall\ \ell\in\integer^2\ , 0<|\ell |\leq K \right\} \subset
{\cal I},
\end{equation}
is a closed set and 
\beq{gamma bar} 
K := \frac{6 (N-1)}{s_0} 
\,\log\frac1\e\ ,\quad \a_0 :=2\sup_{ {\cal I}_{\s_0}} |h_0''|
\widetilde r  K={\rm const}\sqrt{\e} \log \frac{1}{\e}\ ,
\eeq
It results that
\beq{measIU}
{\rm meas}(\cI\setminus U)\leq {\rm const}\ \a_0 =
O \Big( \sqrt{\e} \log \frac{1}{\e} \Big).
\eeq
Moreover
\begin{equation}\label{stima cubica}
\sup_{J\in U^2_{\widetilde{r}}} |g(J,z , {\overline z})|\leq
\,{\rm const}\,\Big( |z|+|\bz|\Big)^3,
\end{equation}
\beqa{eigenst2}
&&\inf_{J\in U} \widetilde\O_j ( J ) \geq
\inf_{J\in U^2_{\widetilde{r}}} |\widetilde\O_j (J)|
\geq \chi_0 > 0, \nonumber\\
&& \inf_{J\in U}
\Big(\widetilde\O_2 (J) -\widetilde\O_1 (J) \Big) \geq
\inf_{J\in U^2_{\widetilde{r}}}
|\widetilde\O_2 (J) - \widetilde\O_1 (J)| \geq \chi_0 > 0,
\eeqa
for some positive constant $\chi_0 $.
\end{proposition}


{\bf Proof:} The Hamiltonian \equ{Hstorta} has been deduced
in [BCV] in the case $ N = 3 $ and on the smaller domain
$ U' := \{ J \in \cI \ | \
| h'_0(J) \cdot \ell | \geq {\rm const}
(\sqrt{\e} \log^{\tau +1}( 1 / \e )) / (1 + |l|^\tau )$,
$\forall \ell \in \integer^2 \setminus  \{ 0 \} \} \subset U $,
see formula (2.22) of [BCV] and
introduce the complex coordinates \equ{coco}.
Note also that the analyticity constant
$ \widetilde{\rho} = {\rm const } $ is bigger than
 $ \widetilde{\rho} = \log^{-1} (1 / \e )$ given in [BCV] .
For the proof of the Proposition 
see the Appendix.
\hfill$\Box$

\giu

The manifold $ \{ z = 0 \} $ is invariant under the Hamiltonian
system generated by the integrable Hamiltonian
$ h_\e(J) + \e \widetilde \O (J) z \overline z $,
\beq{Hint3}
\dot{J} = 0, \qquad
\dot{\psi} = \partial_J h_\e (J) +
\e \partial_J \widetilde \O(J) z \overline z, \qquad
\dot{z} = \ii \e \widetilde \O(J) z
\eeq
and it is completely filled up by the invariant elliptic tori
$ {\cal T}(J_0) :=$ $\{ J = J_0,$ $ \psi \in \torus^2, $ $ z = 0 \} $,
supporting the linear flow
$ t \to \{ J_0 , \ \psi_0 + \o_\e (J_0) t, 0 \} $
with torus frequency $ \o_\e (J_0) :=$ $ \partial_J h_\e (J_0) $.
This family of $2$-dimensional tori will not persist in its entirety
for the complete Hamiltonian system generated by $ \cH $
\begin{equation}\label{Heq3}
\dot{J} = - \partial_{\psi } \cH,\quad
\dot{\psi} = \partial_{ J } \cH, \quad
\dot{z} = \ii \partial_{\overline z } \cH,\quad
\dot{\overline z} = - \ii \partial_{z} \cH, \
\end{equation}
due to resonances among the oscillations.
However, the persistence of a set of positive ($2$-dimensional)
measure of perturbed elliptic invariant tori -those with diophantine
frequency $ \o_\e (J_0) $ -
has been proved in [BCV], using the KAM Theorem of [P\"o2].

In Theorem \ref{theorem 1 - 3 body} we will prove
the existence of an abundance of periodic solutions
with larger and larger period accumulating
on each perturbed elliptic torus,
applying Theorem \ref{theorem 1}.
For this aim, we will first reprove (Theorem \ref{quasip})
the existence of the [BCV]-elliptic tori,
using the KAM Theorem of [P\"o1] and furnishing also the
KAM-normal form describing the dynamics in its neighborhood
(it will be of the form considered in \equ{ham (cal H)}).
Moreover, for proving that the hypothesys
of Theorem \ref{theorem 1} are satisfied,
we also need precise informations on the form of the
KAM-transformation bringing
into the KAM normal form, see
Theorem \ref{KAMPOESCHL}-$(ii)$.

For the existence result of the elliptic invariant tori (Theorem \ref{quasip})
and their surrounding periodic orbits (Theorem \ref{theorem 1 - 3 body})
it is sufficient, and we will assume, $ N =  3 $.


\giu

First of all we rewrite the Hamiltonian ${\cal H} $ in
(\ref{Hstorta}) in a form suitable to apply the KAM Theorem
of [P\"o1].
Introducing the coordinate $ y \in \real^n $ around each
torus $ { \cal T}(J_0) $ in the usual way
$ J = J_0 + y $, the Hamiltonian $ \cH $ can be developed as
\begin{equation}\label{Hdritta0}
H ( y, \psi ,z, \overline z; \o ):= \cH ( J_0 + y, \psi , z, \overline z)= {\cal N} + P
\end{equation}
where $ {\cal N} := $ $ h_\e ( J_0 ) +$ $ h'_\e( J_0 ) \cdot y +$ $ \e
\widetilde \O ( J_0 ) z \overline z$
and $ P=P(y,\psi , z, \bz; J_0 )= P_1+P_2+P_3+P_4 $ with
\beqa{P1234}
P_1 & := & h_\e(J_0+y)-h_\e(J_0)- h'_\e( J_0 ) \cdot y=O(|y|^2),  \nonumber \\
P_2 & := & \e \Big( \widetilde \O (J_0+y)-
\widetilde\O (J_0)\Big)z\overline z= O\big(\e |y||z||\bz|  \big), \nonumber\\
P_3 &:=& \e g(J_0+y,z,\overline z)=O\Big(\e \big(|z|+|\bz|\big)^3\Big), \nonumber \\
P_4 & := & \e^3 f( J_0 + y, x, z, \overline z )=O(\e^3)\ . \eeqa

By the non-isocronicity property det $ h_\e'' (J_0 ) \neq 0 $
there is a one-to-one correspondence between the actions
$ J_0 \in U\subset\real^2 $ and the torus frequency
$ \o := \o (J_0) = h_\e'(J_0) \in {\cal O} \subset\real^2 $
where
\beq{wO}
{\cal O} :=\left\{\,
\o=h'_\e (J_0)\ : \  J_0\in U\,
\right\}.
\eeq
By \equ{measIU} we have
\beq{measIwO}
{\rm meas}\left(h'_\e(\cI)\setminus{\cal O}
\right)\leq {\rm const}\ \a_0 = O \Big(  \sqrt{\e} \log
\Big( \frac{1}{\e} \Big) \Big).
\eeq
The frequencies $\o:=\o(J_0)$ may be introduced as parameters:
denoting
$ J_0 = J_0 ( \o ) := (h_\e')^{-1} ( \o ) $ its inverse function, the
Hamiltonian \equ{Hdritta0} can be finally written in the form
\begin{equation}\label{Hdritta}
H( y, \psi ,z, \overline z; \o ) :=
{\cal N} (y, z, \bz; \o ) +  P (y,\psi , z, \bz; \o)
\end{equation}
where
$ {\cal N}(y,z,\overline z;\o) :=$
$e(\o ) + \o \cdot y + \O(\o) z \overline z $,
$ e ( \o ) :=  h_\e ( J_0 (\o ) )  $,
$ \O ( \o ) := \e {\widetilde{\O}} ( J_0 ( \o ) ) $
and the perturbation
$ P (y,\psi , z, \bz; \o) $ is obtained by
the one in \equ{Hdritta0} just replacing $ J_0 $ with
$ J_0 ( \o ) := ( h_\e')^{-1} ( \o )$.
Recalling \equ{intornini}
the Hamiltonian $ H $ in \equ{Hdritta} is real analytic in
\begin{equation}\label{dominio1}
(y,\psi ,z,\bz; \o ) \in
 U^2_{\widetilde r/2} \times \torus^2_{\widetilde s}
\times D^4_{\widetilde\r} \times  {\cal O}^2_{\widetilde \d }   ,
\qquad {\rm with}\qquad  \widetilde\d:= {\rm const}\, \sqrt{\e}
\end{equation}
(since $\o \to J_0 ( \o ) := (h_\e')^{-1} $ is analytic).

$H$ is in a suitable form to apply the KAM Theorem of [P\"o1]
that we rewrite in the next subsection.

\subsection{A KAM Theorem for elliptic tori}

The KAM Theorem of [P\"o1] applies to
Hamiltonians like
$$
H := H(y,\psi,z,\bz;\o) := {\cal N} + P = e(\o) +
\o\cdot y + \O(\o) z\bz + P(y,\psi,z,\bz; \o)
$$
where
$ (y,\psi,z,\bz) \in \real^n \times $ $ \torus^n \times $ $ \complex^{m}
\times $  $ \complex^{m} $ and $ \o \in \real^n $ is
regarded as a parameter varying over a compact
subset
$ \cO \subseteq \real^n $.

The functions $ P = P(y,\psi,z,\overline z;\o) $ and  $ \O (\o ) $ are
real analytic on the complex domain
$ \cD_{\br,\bs}$ $ \times \cO_\d $ where
$ \cD_{\br,\bs} := D^n_{\br^2} \times \torus^n_{\bs} \times D^{2m}_{\br} $
and $ \br, \bs, \d $ are suitable positive constants.

The size of the perturbation $ P $ is measured by the
following norm. Taking the Fourier-Taylor expansions
$$
P(y,\psi,z,\overline z; \o ) = \sum_{{\ell\in
\integer^n}\atop{a,\overline a\in \natural^m}}
P_{a,\overline a ,\ell}(y; \o )\,z^a \overline z^{\overline a}\, e^{\ii\ell
\cdot \psi} = \sum_{{\ell \in \integer^n}}
P_{\ell}(y,z,\overline z; \o)\, e^{\ii\ell \cdot \psi},
$$
let define\footnote{
The relation with the usual Fourier norm (used in the Appendix)
$ \| P\|_{\br,\bs,\d} := $ $\sum_{\ell\in\integer^n}|P_\ell|_{\br,\d}
e^{|\ell|\bs}$ is
$ \|| P\||_{\br/2,\bs,\d}\leq $ $2^{2m}\| P\|_{\br,\bs,\d} $ and
$ \|| P\||_{\br,\bs,\d}\geq $ $2^{2m}\| P\|_{\br,\bs,\d}\ .$}
 $$
\|| P \||_{\br,\bs,\d} := \sum_{\ell\in\integer^n} |{\bf M}P_\ell|_{\br,\d}
e^{|\ell|\bs}$$
where
$ {\bf M} P_\ell:=$ $ \sum_{a,\ba\geq 0} |P_{a,\ba,\ell}(y)|z^a \bz^{\ba}$
and
$ | \cdot|_{\br,\d} $ denotes the sup-norm over the
$ (y, z, \bz ) \in $ $ D^n_{\br^2} \times D^{2m}_{\br} $ and $ \o \in $
$ \cO_\d $.
The following Theorem follows from [P\"o1].

\begin{theorem}\label{KAMPOESCHL} {\rm {\bf (P\"oschel [P\"o1])}}
Fix $ \tau > n - 1 $. Suppose \beq{eMMe} \sup_{\o \in \cO_\d } |
\partial_\o \O ( \o ) | \leq M \eeq for some $ 0 <  M < + \infty
$, and that the non-resonance condition \beq{nonresKAM} | \O(\o)
\cdot  h | \geq {\a} \ ,\qquad \forall\ 1\leq |h|\leq 2, \ h \in
\integer^n, \qquad \forall \o \in \cO, \eeq is satisfied for some
$ \a > 0 $. Then, there exists a positive constant $ \kappa :=
\kappa ( {\overline s} )$ such that, if $ P $ is sufficiently
small, \beq{KAMcond} \cP:=\|| P\||_{\br,\bs,\d}
\leq\frac{\kappa}{M+1}\a \br^2\leq\frac{\d \br^2}{16}\ , \eeq
then:
\begin{itemize}
\item[(i)] there exist
a normal form ${\cal N}_*:=e_*+\o\cdot y_*+\O_*(\o) z_*\bz_*$,
a Cantor set $\cO(\a)\subset\cO$ on which
\beq{nonresKAM2}
| \o \cdot \ell + \O_*(\o) \cdot h | \geq
\frac{\a}{2(1+|\ell|^\t)}\ ,\quad \forall\ell\in\integer^n\ ,\
|h|\leq 2\ ,\   |\ell|+|h|\neq 0\ ,\  \o\in\cO(\a)\ ,
\eeq
and a transformation
$ \cF: \cD_{\br/2,\bs/2}\times\cO(\a)\longrightarrow \cD_{\br,\bs}\times\cO_\d$
real analytic and symplectic
for each $ \o $ and Whitney smooth in $ \o $, such that
\beq{eqKAM}
H_*:=H\circ\cF={\cal N}_*+R_* \quad with \quad R_*:=\sum_{2|k|+|a+\ba|\geq 3}
R^*_{ka\ba}(\psi_*)y_*^k z_*^a \bz_*^\ba;
\eeq
\item[(ii)] $\cF$ has the following
form (omitting the dependence on $\o$)
\beqa{tracanKAM}
y &=& y_*+ Y(y_*,\psi_*,z_*,\bz_*)\nonumber\\
\psi &=& \psi_*+ X(\psi_*)\nonumber\\
z &=& z_*+ Z(\psi_*,z_*,\bz_*)\nonumber\\
\bz &=& \bz_*+ \overline Z(\psi_*,z_*,\bz_*)\label{trasfKAM}
\eeqa
where
\beq{Ytras}
Y:= \sum_{2|k|+|a+\ba|\leq 2} Y_{ka\ba}(\psi_*)y_*^k z_*^a \bz_*^\ba,
\qquad
Z:= \sum_{|a+\ba|\leq 1} Z_{a\ba}(\psi_*) z_*^a \bz_*^\ba,
\eeq
and, denoting by $\|\cdot\|_*:=\sup_{\cD_{\br/2,\bs/2}}|\cdot |$,
\beq{perideKAM}
\| Y\|_*\ ,\ \ \frac{\br^2}{\bs}\| X\|_*\ ,\ \ \br \| Z\|_*\ ,\ \ \br \|\overline Z\|_*
\leq {\rm const} \frac{\cP}{\a};
\eeq
\item[(iii)] if $M$ is sufficiently small, i.e., if
\beq{MKAM}
M < 1/4\ ,
\eeq
then
\beq{measKAM}
{\rm meas} \Big( \cO\setminus \cO(\a) \Big)=O(\a d^{n-1})\
\eeq
where $d $ is the exterior diameter of $ \cO $.
\end{itemize}
\end{theorem}

By Theorem \ref{KAMPOESCHL} the Hamiltonian system generated
by $ H $ possesses a family of $n$-dimensional elliptic invariant tori
$ {\cal T} := $ $ \{ (y, \psi, z) =$ $ ( Y (0,\psi_*, 0,0), \  \psi_*, \
z ( \psi_*, 0,0 ))\} $
traveled with frequencies $ \o$,
for each frequency vector in the Cantor set $ \cO_\a $.
The dynamics near each torus ${\cal T} $
is described by the normal form \equ{eqKAM}, which in turn coincide with
\equ{ham (cal H)}.

\giu

Let us make some comments on Theorem \ref{KAMPOESCHL}.
Point $(i)$ is Theorem A of [P\"o1] (we have fixed the constant
$ \r $ in the statement of  Theorem A as
$ \r := \overline s / 4 $).

Concerning point $(ii)$, formula \equ{tracanKAM} follows from Section 4 of
[P\"o1] (see, in particular, page 574 of [P\"o1]).
Formula \equ{perideKAM} follows from Section 7 of [P\"o1]
(see in particular the last estimate on page 592).

Point $(iii)$ follows from Theorem B of [P\"o1].
Indeed, from \equ{eMMe}, \equ{nonresKAM} and \equ{MKAM}
$ \cO $ is {\sl essentially non-resonant} according to
the definition at page 565 of [P\"o1].
For finite dimensional systems this reduces to verify only
that $ \cO $ is {\sl non-resonant} (see again page 565 of [P\"o1]),
namely
\beq{nonre}
\min_{\o\in\cO}|\ell\cdot\o+\O(\o) \cdot h | \geq \a
\eeq
for all $ 0 < |h| \leq 2 $ and $ \ell \in {\cal K}_h,$ the closed
convex hull of the gradient set
$ {\cal G}_h:=\{ \partial_\o (\O(\o) \cdot h) \ :  \ \o \in \cO \}.$
Since, by \equ{eMMe} and $ | h | \leq 2 $,
\footnote{$B^n_R \subset \real^n$ is the closed ball of radius $R$
and center $0$ using  the $| \cdot |_2 $ norm.}
$ { \cal K }_h \subseteq B^n_{2 M},$
taking $ M < 1 / 4 $ as in (\ref{MKAM}),
the unique integer vector $l \in  { \cal K }_h $ is $ l = 0 $.
Hence condition \equ{nonre} must be verified for $\ell=0$ only,
and this is condition (\ref{nonresKAM}).

\subsection{Abundance of periodic solutions in the three-body problem}

In this section we prove the existence of periodic orbits
accumulating on elliptic $2$-dimensional tori of the
three body-problem.

\giu
First of all we show that the KAM-Theorem \ref{KAMPOESCHL} applies
to the Hamiltonian $ H $ in (\ref{Hdritta})
reproving the existence
of elliptic invariant tori in the spatial planetary three body problem as in
[BCV]. We assume the frequency parameter $\o \in \real^2 $ to
vary over the compact subset $ \cO $.

By \equ{dominio1} and \equ{intornini}, the Hamiltonian $ H $ in
(\ref{Hdritta}) is  real analytic on $ \cD_{\br,\bs} \times \cO_\d
:= D^2_{\br^2} \times \torus^2_{\bs} \times D^{4}_{\br} \times
\cO_\d $ for \beq{intornini2} \br:=\sqrt{c_0\e}\ ,\qquad \bs:={\rm
const}\ ,\qquad \d:={\rm const \sqrt\e}\ , \eeq where $ c_0 $ is a
small constant which will be determined later on (we restrict the
domain of $ H $ w.r.t. to the larger one defined in \equ{dominio1}
because with the choice \equ{intornini2} the smallness KAM
condition \equ{KAMcond} is satisfied, see \equ{cP}, \equ{sdfa}).

Applying the KAM-Theorem \ref{KAMPOESCHL}, we get
\begin{theorem}\label{quasip}
Fix $\tau>n-1.$ For $ c_0 $ and $\e $ small enough, there exists a
Cantor set $ \cO ( \e ) \subset \cO \subset h'_\e(\cI)$, with
\beq{misura 2} {\rm meas}\big( h'_\e(\cI)\setminus\cO(\e)\big)
\leq \,{\rm const}\, \sqrt{\e} \log \,\frac{1}{\e} \,,\eeq such
that, for any $ \o = \o(J_0) \in \cO ( \e ),$ there exists a
symplectic trasformation
$$
\Phi : \cD_{\br/2,\bs/2}\longrightarrow \cD_{\br,\bs}\ ,
$$
of the form (\ref{tracanKAM}), transforming the
three-body Hamiltonian $H$ in \equ{Hdritta}
into the normal form $ H_* := H \circ \Phi $
as in (\ref{eqKAM}).
Moreover
\beq{gandalf}
| \o \cdot \ell + \O_*(\o) \cdot h | \geq
\frac{{\rm const} \,\e}{1+|\ell|^\t}\ ,\qquad \forall\ell\in\integer^n\ ,\
|h|\leq 2\ ,\   |\ell|+|h|\neq 0\ ,\  \o\in\cO(\e)\ .
\eeq
In particular, the three-body problem Hamiltonian $ H $
possesses a family of $2$-dimensional elliptic invariant tori
traveled with frequencies $ \o$.
\end{theorem}

\noindent {\bf Proof:}
Recalling that $ \O (\o ) := \e \widetilde \O ( J_0 (\o ))$,
it is easy to see that
$ \sup_{\o \in \cO_\d } | \partial_\o \O ( \o ) | \leq \e  C_0 := M $
for some positive constant $ C_0 $ and
condition \equ{eMMe} holds true. Moreover, for $ 0 <  \e < 1 \slash (4 C_0)$,
$ M < 1 \slash 4 $ and also condition \equ{MKAM} is satisfied.

We claim that condition (\ref{nonresKAM}) holds, with
\beq{alfa}
\a:= \chi_0 \e,
\eeq
namely,
\beq{nonre1}
| \Omega (\o ) \cdot h | \geq \alpha := \chi_0 \e, \qquad
\forall 1 \leq |h| \leq 2, \ h \in \integer^2, \quad \forall
\o \in \cO.
\eeq
Indeed 
$ \forall \o \in \cO $,
$ \O (\o ) = \e \widetilde \O ( J_0 ) $ for some $J_0 \in U $.
Moreover, by \equ{eigenst2}, it is easy to see that
$ | \widetilde \O ( J_0 ) \cdot h | \geq \chi_0 $, $\forall 1 \leq
|h | \leq 2 $, $ h \in \integer^2 $, $ J_0 \in U $ and the
claim follows. 


It remains to check the smallness condition (\ref{KAMcond}), i.e.
$ \cP \leq \e \overline{r}^2 $, since $ \kappa = O(1) $ and $
\alpha := \chi_0 \e $. Note that, from  \equ{P1234} and
\equ{intornini2}, \beq{cP} \cP \leq {\rm const}(\br^4+\e\br^4+\e
\br^3+\e^3)\leq {\rm const}(\br^4+\e^3). \eeq
Hence, by \equ{cP},  
in order to check condition (\ref{KAMcond}),
it is sufficient that
\beq{sdfa}
\br^4+\e^3\leq {\rm const}\ \e \br^2.
\eeq
Since $ \br := \sqrt{c_0 \e } $,  \equ{sdfa}
holds true for $ c_0 $ and $\e$ small enough.
By Theorem \ref{KAMPOESCHL} the result follows.

From \equ{measKAM}, since $ d = O(1) $ and $ \alpha = O(\e )$,
meas $( \cO \setminus \cO( \e ) )=$ $O( \e )$. Moreover, since
by \equ{measIwO}, meas $(h_\e'( {\cal I}) \setminus \cO)=$
$O( \sqrt{\e} \log ( 1 / \e ) )$,
the measure estimate \equ{misura 2} follows.
\hfill$\Box$
\giu

For proving that each elliptic invariant torus
found in Theorem \ref{quasip}
lies in the closure of periodic orbits
of the three body problem, we will apply
Theorem \ref{theorem 1} to the normal form
Hamiltonian $ H_* := H \circ \Phi $ given in  (\ref{eqKAM})
(which has the form \equ{ham (cal H)})
corresponding to the three body problem.

The crucial hypothesys to verify in Theorem \ref{theorem 1}
is the non-degeneracy of the twist matrix $ \cR $ defined in (\ref{R2}).
We have to evaluate, for $1\leq i,i'\leq n = 2 $ and
$ 1\leq j\leq m =2 $
\beqa{R*R*}
R^*_{e_i+e_{i'},0,0}(\psi_*)&=&\frac{1}{1+\d_{i,i'}}
\frac{\partial^2 H_*}{\partial y_{*i}\partial y_{*i'}}(0,\psi_*,0,0)\ ,
\nonumber\\
R^*_{e_i,e_j,0}(\psi_*) & = &
\frac{\partial^2 H_*}{\partial y_{*i}\partial z_{*j}}(0,\psi_*,0,0)\ .
\eeqa
$ R^*_{e_i,0,e_j}(\psi_*) $ can be computed as the
$ R^*_{e_i,e_j,0} (\psi_* ) $ and their calculation is omitted.
Due to the special form (\ref{tracanKAM})-(\ref{Ytras})
of the canonical transformation $ \Phi $, we have
\beq{y*}
\partial_{y_{*i}} \left( H\circ \Phi \right)
=\left[\left(\partial_{y}H\right)\circ\Phi\right]\partial_{y_{*i}} y
\qquad {\rm with}\qquad
\partial_{y_{*i}}y=e_i+Y_{e_i,0,0}(\psi_*)\ ,
\eeq
from which
\beqa{dyy}
\partial^2_{y_{*i'} y_{*i}} \left( H\!\circ\! \Phi \right)
&=&
\left[\left(\partial^2_{yy}H\right)\!\circ\!\Phi\ \partial_{y_{*i'}}y\right]
\partial_{y_{*i}} y
\\
\partial^2_{z_{*j} y_{*i}} \left( H\!\circ\! \Phi \right)
&=&
\Big[
\left(\partial^2_{yy}H\right)\!\circ\!\Phi\ \partial_{z_{*j}}Y+
\left(\partial^2_{zy}H\right)\!\circ\!\Phi\left(e_j+\partial_{z_{*j}}Z\right)+\nonumber \\ &+&
\left(\partial^2_{\overline{z}y}H\right)\!\circ\!\Phi\ \partial_{z_{*j}}\overline{Z}
\Big]
\partial_{y_{*i}} y \ . \label{dzy}
\eeqa

We are now able to prove that

\begin{lemma}\label{nonde3}
For $c_0 $ and $\e $ small enough,
the twist matrix $ \cR $ of the planetary three body problem
is invertible and $ | \cR^{-1}| = O(1) $.
\end{lemma}

\noindent {\bf Proof:}
We need to evaluate \equ{dyy}-\equ{dzy} for
$( y_*, \psi_*, z_*, \bz_* ) = $
$(0,\psi_*,0,0)=:\star$.
By \equ{tracanKAM}, \equ{perideKAM}, \equ{cP} and since $r_0 :=
\sqrt{c_0 \e }$ and $ \alpha := \chi_0 \ \e $, it results
\beq{stima1}
| y_i ( \star ) |= |Y_i(\star)|\leq {\rm const} \, \frac{\cP }{\a }
\leq {\rm const} \frac{ {\overline r}^4 + \e^3 }{ \e } =
{\rm const} \frac{ c_0^2 \e ^2 + \e^3 }{ \e } \leq
{\rm const} \ \e .
\eeq
Moreover, by \equ{perideKAM} and standard 
``Cauchy estimates''\footnote{Cauchy estimates allow to bound $n$-derivatives
of analytic functions on a set $A$ in terms of their sup-norm 
on larger domains $A \subset A' $ divided by 
${\rm dist} (\partial A, \partial A')^n$.}, we have
\beq{stimad}
| \partial _{y_*}Y(\star)|\ ,\ \
\frac{1}{\br}| \partial _{z_*} Y(\star)|\ ,\ \
 | \partial _{z_*} Z(\star)|\ ,\ \
|\partial _{z_*} \overline Z(\star)|
\leq {\rm const} \frac{\cP}{\a \br^2} = O( c_0 ),
\eeq
(for the second estimate note that $ \partial _{z_*} Y $ is 
independent of $ y_* $ due to (\ref{Ytras})).

From \equ{dyy}-\equ{dzy} we deduce, using
\equ{stima1}, \equ{stimad} and
since $ H = h_0 (J_0 + y ) + O(\e) $, with $ h_0 $
defined in (\ref{specham}),
\beqano
\partial^2_{y_{*i} y_{*i}} H_* (\star) & = &
\frac{-3 \k_i}{ ((J_0)_i + y_i ( \star) )^{4}}
+O_\sharp\left(\frac{\cP}{\a\br^2}\right)+O_\sharp(\e)
=\frac{-3 \k_i}{(J_0)_i^{4}}+O_\sharp(c_0), \\
\partial^2_{y_{*i'} y_{*i}} H_* (\star)
&=& O_\sharp( c_0 ) +   O_\sharp(\e) = O_\sharp( c_0 )
\qquad {\rm if } \quad  i\neq i',\\
\partial^2_{z_{*j} y_{*i}} H_* (\star), 
&=&
O_\sharp\left(\frac{\cP}{\a\br}\right)
+O_\sharp(\e)=O_\sharp(c_0^{3/2}\sqrt\e),
\eeqano
where $ O_\sharp(c) $ denotes a function of $ \psi_* $ with sup-norm
for $ \psi_*\in\torus_{\bs/2}$ smaller than a constant multiplied by $c.$
It follows
$$
R^*_{2e_i,0,0}(\psi_*)=\frac{-3 \k_i}{2 (J_0)_i^4}+O_\sharp(c_0),
\ R^*_{e_i+e_{i'},0,0}(\psi_*)=O_\sharp(c_0) \ {\rm if} \ i' \neq i,
\ R^*_{e_i,e_j,0}(\psi_*)=O_\sharp(c_0^{3/2}\sqrt\e)\ .
$$
Evaluating the Fourier coefficients of the above function we obtain
\beqa{resti}
&& R^*_{2e_i,0,0,0}=\frac{-3\k_i}{2(J_0)_i^4}+O(c_0) \ ,
\qquad
R^*_{e_i+e_{i'},0,0,0}=O( c_0 )\ , \ {\rm for} \  i' \neq i \\
&&|R^*_{e_i,e_j,0,\ell}|\ ,\ \ |R^*_{e_i,0,e_j,\ell}| \leq {\rm
const}\ c_0^{3/2}\sqrt\e \ e^{-|\ell |\bs/4}. \label{resti2} \eeqa
From \equ{gandalf}, $ | \o \cdot \ell + \O_{*j} (\o) | $ $\geq$
const $ \e \slash ( 1 + |\ell|^\t )$, $ \forall \ell \in
\integer^n $ and then, from \equ{resti2} the second addendum in
definition of the twist matrix $ \cR $ introduced in (\ref{R2}) is
\beqa{secadd} \sum_{1 \leq j \leq m \atop{\ell \in \integer^n}}
\frac{1}{\o\cdot\ell+\O_{*j}( \o )} \left(
R_{e_i,e_j,0,\ell}^*R_{e_{i'},0,e_j,-\ell}^* +
R_{e_i,0,e_j,-\ell}^*R_{e_{i'},e_j,0,\ell}^* \right) & = &
\nonumber\\O \Big( \sum_{\ell \in \integer^n} \frac{1 + |\ell|^\t
}{\e} c_0^3 \e \ e^{- \frac{|\ell | \bs }{2}} \Big) & = & O (
c_0^3 ). \eeqa Finally \equ{resti} and \equ{secadd}  imply that
the entries of the twist matrix $ \cR $ corresponding to the
spatial three body-problem are 
\beq{cRii} \cR_{i,i} =
\frac{-3\k_i}{2(J_0)_i^4}+O(c_0)+O(c_0^3)\ ,\qquad \cR_{i,i'} =
O(c_0)+O(c_0^3) \ {\rm for} \  i' \neq i. \eeq 
By \equ{cRii}, for $ c_0 $
small enough, the matrix $ \cR $ is invertible and $|\cR^{-1}| =
O(1)$. \hfill$\Box$

\giu


We can finally prove the abundance of periodic solutions
in the spatial planetary three body problem.

\giu

\noindent {\bf Proof of Theorem \ref{theorem 1 - 3 body}:} Since
the number of elliptic variables $ z $ is $ m = 2 $ condition
($a$)-($i$) of Theorem \ref{theorem 1} holds. Moreover, by Lemma
\ref{nonde3}, the twist matrix $ \cR $ is invertible with $ |
\cR^{-1}| = O(1)$, and then we can apply Theorem \ref{theorem 1},
proving the existence of an abundance of periodic solutions of the
three body problem. Theorem \ref{theorem 1 - 3 body} is finally
proved. \hfill$\Box$

\section{Periodic orbits near 
resonant elliptic tori} \label{section theorem 2}

In this section we study the persistence, for $ \e > 0 $, of the
circular decoupled periodic motions of the planets
around the star, once 
suitable conditions on the period and the masses of the bodies
are satisfied.

Let consider 
$ J_0\in {\cal I} $, $ T > 0 $ and $ k=(k_1,k_2) \in \integer^2 $ with
gcd$(k_1,k_2)=1$, such that
\beq{malachia} \o:= h_0' (J_0)= \frac{2\pi k}{T}\ ,\eeq where
$ h_0 $ is the integrable Hamiltonian defined in  \equ{specham}.
The decoupled three-body problem possesses the
family of (circular) periodic solutions
$$
\check\z_{\psi_0} (t):= ( J_0,\,\psi_0+\o  t,\,0,\,0 ),
$$
with minimal period $ T $ and parametrized by 
$ \psi_0 \in \torus^2 $.

Since $ h_0 $ is properly degenerate the persistence of 
these motions 
for $ \e > 0 $ can not be established without further informations
on the perturbation $ f_0 $. 
For continuing some of these
solutions we will exploit 
the normal form ${\cal H}$ defined in \equ{Hstorta}
(with $ N = 4 $).

First of all, through a Lyapunov-Schmidt reduction similar to the one
of Lemma \ref{pseudo (1)}, we obtain the existence of
suitable pseudo-$T$-periodic solutions.

\begin{lemma}\label{ex u e v}
There exists $ c_1, \e_1 > 0 $,  such that, for all $T>0$,
$k\in\integer^2$ with gcd$(k_1,k_2)=1$,  $0<\e \leq \min\left\{
\e_1 , \frac{1}{c_1 T}\right\}$, $J_0\in U$ satisfying
$ \o := h_0' (J_0)= 2\pi k \slash T $, 
and $ \forall \ \psi_0 \in \torus^2 $,
there exists a unique function
$$ \z_{\psi_0} =(J_{\psi_0},\psi_{\psi_0},
z_{\psi_0})\in C^0 ([0,T], \,\real^2\times\torus^2\times
\complex^2) \cap C^1 ( (0,T), \,\real^2\times\torus^2 \times
\complex^2),$$ smooth in $ \psi_0 $, so that, $ \forall t \in
(0,T)$, $ \z_{\psi_0}(t) $ is a solution of the Hamiltonian system
\equ{Heq3}\footnote{For the 
Hamiltonian ${\cal H}$ defined in \equ{Hstorta}
with $ N = 4 $.} with $ \psi_{\psi_0}(0) = \psi_{\psi_0}(T) =
\psi_0 $, $ z_{\psi_0}(0) = z_{\psi_0}(T) $ and
$$ \sup_{t\in [0,T]} |\z_{\psi_0}(t)-\hat\z_{\psi_0} (t)|\leq
c_1 \e^2\ ,
$$
where $\hat\z_{\psi_0} (t):=$ 
$( J_\e,\,\psi_0+\o t,\,0,\,0 )$ and 
$ J_\e = J_0 + O ( \e ) $ is the unique solution of
$ h_\e'(J_\e)= \o $.
\end{lemma}

\noindent {\bf Proof:} We first note that, since det $h_0''(J)
\neq 0$ (see \equ{hamform}--\equ{specham}), by the 
Implicit Function Theorem, 
$J_\e$ exists, is
unique and is $\e$-close to $J_0\in U$. Moreover, for $\e_1$ small
enough, $J_\e\in \real^2\cap U_{{\rm const}\,\e}$ $\subset
U_{\widetilde r/4} $ (since $\widetilde r=O(\sqrt\e)$).

Define $X$ and $Y$ as in Lemma
\ref{lemma 1}, with $n=m=2$. Let also $\O:=\e \widetilde\O(J_\e)$,
$M:=\partial^2_J h_\e (J_\e)$ and $\cM:= {\bf 1}_2 -e^{\ii \e
\widetilde\O(J_\e) T}$. Notice that, if $0 < \e T \leq 1/ c_1 $ is small enough
(namely $ c_1 $ is large enough), then $ | \cM^{-1} | $ is
bounded by a constant (independent on $\e$). 

We look for
$T$-periodic solutions of the Hamiltonian system \equ{Heq3} of
the form $ \z =$ $ \hat{\z}_{\psi_0}+$ $
(\widetilde J,\widetilde\psi,\widetilde z,
\widetilde{\overline{z}})\,.$ Set
$$ P
\left(
\matrix{
\widetilde J \cr
\widetilde \psi \cr
\widetilde z
}
\right) =
\left(\matrix{
-\e^4 \partial_\psi f(\star)\cr
\cr
-\partial^2_J
h_\e(J_\e)\,\widetilde J +
\partial_J
h_\e (J_\e+\widetilde J)-
\partial_J
h_\e (J_\e)+\e\partial_J\widetilde \O (J_\e+\widetilde J)
\widetilde z\cdot \overline{\widetilde z} +\cr \,\,\,\,\,\,+\Big(
\e \partial_J g +\e^4 \partial_J f\Big) (\star)\cr \cr \ii \Big[
\e \Big( \widetilde\O(J_\e+\widetilde J)- \widetilde\O(J_\e)
\Big)\cdot\widetilde z + \Big( \e \partial_{\overline z} g +\e^4
\partial_{\overline z} f\Big) (\star)\Big] } \right)\,,$$ where
the star above denotes, for short, $ \star := (J_\e+\widetilde
J,\psi_0+\o t+ \widetilde\psi,\widetilde z,
\widetilde{\overline{z}}) $.  
We note that $P(0)=O(\e^4)$ by \equ{stima cubica}. We want
to apply Lemma \ref{lemma 0}. From \equ{eq *} of Lemma \ref{lemma
1} we find $|L|=O(1/(c_1^2\e^2))$, hence $\d_0:=2\,|L(P(0))| =
O(\e^2/c_1^2)$ and $\sup_{x\in B_{\d_0}}|DP (x)|\leq \, {\rm
const}\,(\d_0+\e^4)$. Then, taking $c_1$ big enough, we can apply
 Lemma \ref{lemma 0} proving the existence and uniqueness of a
$\z^\star_{\psi_0} \in B_{\r_0}$ so that
\begin{equation}\label{equation 007}
\z^\star_{\psi_0} = L(P( \z^\star_{\psi_0} )).
\end{equation}
Let $\z_{\psi_0}(t):=\hat \z_{\psi_0} +\z^\star_{\psi_0}$.
By means of \equ{equation 007} and \equ{eq **},
it follows that $\z_{\psi_0}$ satisfies
the Hamilton equations \equ{Heq3} in $ ( 0, T) $
and the boundary conditions
$ \psi_{\psi_0}(0) = \psi_{\psi_0}(T) = \psi_0 $,
$ z_{\psi_0}(0) = z_{\psi_0}(T) $.
As usual, by the Implicit Function Theorem, 
$ \psi_0 \to \z_{\psi_0} $ is smooth.
\hfill$\Box$

\giu

Finally, critical points of the ``reduced action functional''
$\cE: \torus^2 \to \real $ defined by
$$
\cE( \psi_0 ) := \int_0^T J_{\psi_0}\cdot \dot\psi_{\psi_0}
+\ii z_{\psi_0}\cdot \dot\bz_{\psi_0}-\cH(\z_{\psi_0})\,dt\
$$
give rise, arguing as in Lemma \ref{LEMMA M}, to periodic
solutions of \equ{Heq3}.  As in the proof of Theorem \ref{theorem
1} (see page \pageref{proof theorem 1}) we deduce the existence of
at least two geometrically distinct periodic solutions of
\equ{Heq3} (corresponding to the points of maximum  $\psi_0^+$ and
 minimum $\psi_0^-$  of $ \cE $).

We can finally state the following result which, in particular,
will imply Theorems \ref{theorem 2} and \ref{cor6} :

\begin{theorem}\label{Theo}
There exists $ c_1, \e_1 > 0 $,  such that, for all $T>0$,
$k\in\integer^2$ with gcd$(k_1,k_2)=1$, $0<\e \leq \min\left\{
\e_1 , \frac{1}{c_1 T}\right\}$, $J_0\in U$ satisfying
$ \o := h_0' (J_0)= 2\pi k \slash T $, 
there exist at least two geometrically distinct
$T$-periodic solutions $ \z_{\psi_0^\pm} =
(J_{\psi_0^\pm},\psi_{\psi_0^\pm}, z_{\psi_0^\pm}) $ of the
Hamiltonian system \equ{Heq3}\footnote{For the
Hamiltonian ${\cal H}$ defined in \equ{Hstorta} with $ N = 4 $.}
with
$$ 
\sup_{t\in \real} |\z_{\psi_0^\pm}(t)-\hat\z_{\psi_0^\pm} (t)|\leq
c_1 \e^2\ ,
$$
where $ \hat\z_{\psi_0^\pm} (t):= ( J_\e,\,\psi_0^\pm+\o t,\,0,\,0)$ 
and $ J_\e = J_0 + O ( \e ) $ is the unique solution of
$ h_\e'(J_\e)=\o $ $=2\pi k/T$.
Note that $ \sup_{t\in [0,T]} |\z_{\psi_0^\pm }(t)-
( J_0,\, \psi_0^\pm + \o t, \, 0 , \, 0 )| = O( \e )$. 
\end{theorem}

\begin{remark}
{\rm The condition $0< \e T < 1/ c_1 $ (for $c_1 $ large enough)  
required in the Theorem is sharp. Indeed, for
$ \e T = O ( 1 ) $ some further resonance phenomenon
can appear, destroying {\it any} periodic solution:
note that $ e^{\ii\e \widetilde{\Omega}(J_\e) T} - {\bf 1}_2 $
could be degenerate.
It is for reaching this estimate that we take $ N = 4 $
in the Hamiltonian \equ{Hstorta}.}
\end{remark}

We now show that Theorems \ref{theorem 2} and  \ref{cor6} follow
by Theorem \ref{Theo}.
\\[1mm]
Consider $ J_0 \in {\cal I} $, 
$ T > 0$, $ k=(k_1,k_2) \in \integer^2 $ 
with gcd$(k_1,k_2)=1$, such that 
$ \o := h_0' (J_0)= 2\pi k \slash T $, i.e. \equ{malachia} holds. 
Under suitable restrictions on 
$ \e $ and $T$ we shall prove that 
$ J_0 \in U $ where $ U = U ( \e ) \subset {\cal I} $ 
is the domain defined in \equ{U}.

First of all we claim that there exists a constant $c_2>0$ large
enough such that \beq{osea}\e\geq e^{-T / c_2 } \quad\Longrightarrow
\quad \o\cdot\ell\neq 0\ ,\qquad \forall\   0 < | \ell | \leq K\eeq where
$ K = K ( \e ) = {\rm const} \ \log (1\slash  \e )$ is 
defined in \equ{gamma bar}. 
In fact, since gcd$(k_1,k_2)=1$, \beq{lk}
\o \cdot \ell = 0 \quad \Longrightarrow \quad
k \cdot \ell = 0  \quad \Longrightarrow 
\quad \ell \in \integer (-k_2, k_1) \Longrightarrow 
\quad |\ell|\geq |k|.\eeq
Moreover \beq{kcI} |k| = T \frac{|h_0'(J_0) |}{2 \pi}  
\geq c_{{\cal I}} T \eeq 
where
$ c_{{\cal I}} := (1\slash 2\pi) \min_{J_0 \in  {\cal I}} | h_0'(J_0) | > 0 $ 
(note that 
$ h_0'(J_0)= ( \kappa_1 J_{01}^{-3}, \kappa_2 J_{02}^{-3} ) \neq 0 $
and $ {\cal I} $ is compact).

From \equ{lk}-\equ{kcI} and the 
hypothesis $ \e \geq e^{-T / c_2 } $,  
we obtain that $ k \cdot \ell = 0$,
$ \ell \neq 0 $, imply $ | \ell | \geq c_{{\cal I}} T > 
c_{{\cal I}} c_2 \log ( 1 \slash \e ) $. 
For $ c_2 $ large enough, 
$ c_{{\cal I}} c_2 \log (1 \slash \e ) >$ $K =$ $ {\rm const}\log(1/\e)$,
and we conclude that 
$ k \cdot \ell \neq 0$ for all $ 0< |\ell| \leq K$. The claim is proved.

It follows that $ \forall \ 0 < | \ell | \leq K $, 
$k \cdot \ell \geq 1 $ and 
hence $ |\o \cdot \ell| = $
$ 2 \pi |k\cdot\ell| T^{-1} \geq  2 \pi T^{-1}$.
Moreover, if $ \e \log^2 ( 1 / \e ) \leq  c_3 T^{-2} $, 
then $ 2 \pi T^{-1} \geq 2 \pi \sqrt{\e / c_3} \log (1 \slash \e )
\geq \a = $ $ {\rm const} \sqrt\e \log(1 / \e ) $ 
(the constant $ \a $ is defined in  \equ{gamma bar}) 
for a suitable $ c_3 > 0 $ large enough.

We have proved that conditions $ \e \geq e^{-T/ c_2 } $ and
$ \e \log^2 (1/\e) \leq c_3 T^{-2}$ imply $ J_0 \in U$.
In conclusion, defining the functions
\beq{ee} 
\underline \e(T) := \min\Big\{ e^{-T/ c_2 },\ (c_1
T)^{-1},\ \e_1 \Big\}\quad {\rm and}\quad\overline\e(T)
:=\min\Big\{ F(c_3 T^{-2}),\ (c_1 T)^{-1},\ \e_1 \Big\}\ , 
\eeq
where $F$ is the inverse 
of $ G(\e):=\e\log^2(1/\e)$ and
$c_1,\e_1$ are defined in Lemma \ref{ex u e v}, we have proved:
\begin{lemma}\label{lemmaee}
Let $ T > 0 $, $ k=(k_1,k_2) \in \integer^2 $, with
gcd$(k_1,k_2)=1$,  and  $ J_0\in {\cal I} $  such that
$  \o:= h_0' (J_0)= \frac{2\pi k}{T} $ holds. 
Then 
\beq{eee} \underline\e(T)\leq
\e\leq\overline\e(T)\qquad\Longrightarrow\qquad J_0\,\in\, U\ .
\eeq
\end{lemma}

We can finally deduce Theorems \ref{theorem 2} and \ref{cor6}.
\giu

\noindent {\bf Proof of Theorem \ref{theorem 2}:} It is a direct
consequence of Lemma \ref{lemmaee} and of Theorem
\ref{Theo} taking $\underline\e(T)$, $\overline\e(T)$ as in
\equ{ee} and $T_0>0$ as the last instant for which
$\underline\e(T_0)=\overline\e(T_0)$ . \hfill$\Box$

\giu 

\noindent {\bf Proof of Theorem \ref{cor6}:} Define
$$
\underline{T} (\e):=\min \left\{ c_2 \log\frac{1}{\e}\ ,\frac{1}{c_1
\e_1} \right\}\qquad {\rm and}\qquad \overline{T} (\e):=\min
\left\{ \frac{\sqrt{c_3}}{\sqrt\e \log(1/\e)}\ ,\frac{1}{c_1 \e_1}
\right\}
$$
with $\e_1,c_2,c_3$ defined above. Theorem \ref{cor6} is a
direct consequence of Lemma \ref{lemmaee} and of Theorem
\ref{Theo}. \hfill$\Box$

\Giu

\section{Appendix: Proof of Proposition \ref{Hstorta1}}

We first recall the following Averaging Theorem, the proof of which can
be found in [BCV].
We first introduce some notations. For $r,\r,s>0$
and $U\subset \real^2$
we denote the complex set
$W_{r,\r,s}:=U_r\times\torus^2\times D^4_\r.$
For a function
$f:=f(I,\f,p,q)$
real analytic for  $(I,\f,p,q)\in W_{r,\r,s}$
we denote by $\| f\|_{r,\r,s}$ its ``sup-Fourier'' norm given by
$$
\| f\|_{r,\r,s}:=
\sum_{k\in\integer^2}\left(
\sup_{(I,p,q)\in U_r\times D^4_\r}|f_k(I,p,q)|
\right)e^{|k|s}\ .
$$

\protwo{pro:averaging}{Averaging Theorem}
Let ${\cal H}_0:=h_0(I_0)+f_0(I_0,\f_0,p_0,q_0)$
be a real-ana\-ly\-tic Hamiltonian
on $W_{r,\r, s}$. Assume
that there exist $\a_0,K>0$, satisfying $Ks\ge 6$, such that
\beq{picden}
|h'_0 (I_0)\cdot k|\geq \a_0/2\ ,\qquad
\forall\, k\in\{k\in\integer^2:\ 0<|k|\leq K\}\ , \ \
\forall\, I_0\in U_r \ .
\eeq
Assume, also, that, if $d:=\min\{ rs,\r^2\}$, then
\beq{f}
\|f_0\|_{r,\r, s} < \frac{\a_0 d}{c\  Ks}
\ ,
\eeq
where $c>1$ is a suitable (universal) constant.
Then, there exists a real-analytic symplectic transformation
\beq{Psi}\Psi:(I',\f',p',q')\in W_{r/2,\r/2,s/6}\mapsto
(I_0,\f_0,p_0,q_0)=\Psi(I',\f',p',q')\in W_{r,\r, s}
\eeq
and a real-analytic function $g_0=g_0(I',p',q')$
such that
\beq{genri}
{\cal H}':= {\cal H}_0 \circ \Psi =h_0+g_0+f_*\ ,
\eeq
and the following bounds hold\footnote{$f_{0,0}$
denotes the 0-Fourier coefficient
of $f_0$, i.e., its $\f$-average.}:
\beq{g}
\sup_{(I',p',q')\in U_{r/2}\times V_{\r/2}} \
| g_0(I',p',q')-f_{0,0}(I',p',q')|\leq
\frac{c}{\a_0 d}\ \|f\|_{r,\r, s}^2
\ ,
\eeq
\beq{f*}
\| f_*\|_{{r/2,\r/2,s/6}}
\leq e^{-Ks/6}\ \|f\|_{r,\r, s}\ .
\eeq
Furthermore, for each $(I',\f',p',q')\in W_{r/2,\r/2,s/6}$,
$(I_0,\f_0,p_0,q_0)=\Psi(I',\f',p',q')$
satisfies
\beq{spost}
s\ | I_0-I'  |,\ \
r\ | \f_0-\f'  |,\ \
\r \ | p_0-p'  |,\ \
\r \ | q_0-q'  |
\leq \frac{c \|f\|_{r,\r, s}}{\a_0}\ .
\eeq
\epro

\noi
Let, now, ${\cal H}_0=h_0+f_0$ be as in  \equ{hamform}.
and $U$ as in \equ{U}.
The estimate \equ{measIU} directly follows by the definition of
$U.$
Next, let us  choose the sets and the parameters involved in
Proposition~\ref{pro:averaging} as follows:
\beq{list}
K := \frac{6 (N-1)}{s_0}
\,\log\frac1\e\ ,\quad
\a_0 :=2\sup_{ {\cal I}_{\s_0}} |h_0''| r K\ ,
\quad r=c_*\sqrt{\e}\leq \s_0\ ,\quad s:=s_0\ ,
\quad \r := \r_0\ ,
\eeq
where $c_*$ is a suitable large constant to be
fixed later and
$\s_0,s_0,\r_0$ were defined in
\equ{ssr}.
Moreover we better specify the definition of
 $U$  in \equ{U}:
$$
U:= \left\{ J\in \cI\ \ :\ \
|h'_0(J)\cdot\ell|\geq\a_0\ ,
\forall\ \ell\in\integer^2\ , 0<|\ell |\leq K  \right\} \subset \real^2\ .
$$
The estimate \equ{measIU} directly follows by the
previous definition.
Notice that, from these definitions,
there follows (for $\e$ small enough) that
\beq{orders}
\a_0=\,{\rm const}\,  \sqrt{\e} \log\frac{1}{\e}\ ,\qquad
d=\,{\rm const}\,  r\ ,\qquad
\frac{\a_0 d}{c K s}=\,{\rm const}\,  \e c_*^2\ ,
\eeq
(clearly, in the last evaluation, ``$\,{\rm const}\, $" does not involve $c_*$).

\nl
Now, it is not difficult to check that, choosing $c_*$
big enough and letting $\e$ be small enough, assumptions
\equ{picden} and \equ{f} are met. In fact, observing that  $f_0$ in
\equ{hamform}-\equ{specham} is such that
$$
\|f_0\|_{r,\r, s}\le \,{\rm const}\,  \e\ ,
$$
\equ{f} follows from last equality in \equ{orders}, by choosing $c_*$ large enough.
As for \equ{picden},
observe that for any point in $I_0\in U_r$ there is a point $I_0^*\in
U$ at distance less than $r$ from it.
Hence by the
definitions of
$\a_0$ and $r$ and by
Cauchy estimates\footnote{``Cauchy estimates" allow to
bound
$n$-derivatives of analytic
functions on a set $A$ in terms of their sup-norm on larger
domains $A'\supset A$ divided by
$\,{\rm dist}\,(\partial A,\partial A')^n$.},
for any $I_0\in U_r$ and any
$0<|k|\le K$,
\beqano
|h_0'(I_0)\cdot k|&\ge& |h'_0(I_0^*)\cdot k|- |h'_0(I_0^*)-h'_0(I_0)|\ |k|\\
&\ge & \a_0 - \sup_{\cI_{\s_0}} |h''_0| \ r\ K\\
&= & \a_0/2\ ,
\eeqano
which proves also \equ{picden}.
Thus, by Proposition~\ref{pro:averaging}, there exists
a real analytic symplectic transformation
$$\Psi: (I',\f',p',q')\in
\cD_1:=  U_{\frac{r}2}\times \torus^2_{\frac{s}6}\times
D_{\frac{\r}2}^4\to (I_0,\f_0,p_0,q_0)\in
U_{r}\times \torus^2_{s_0}\times D_{\r_0}^4\subset \cD_0  ,
$$
such that
\beqa{eaclose}
&& |I'-I_0|\le \,{\rm const}\,  \frac\e\a_0=\,{\rm const}\,
\frac{\sqrt{\e}}{\log\frac{1}{\e}}\ \nonumber\\
&&|p'-p_0|\ ,\ |q'-q_0|  \le \,{\rm const}\,  \frac\e{\a_0 \r}=\,{\rm const}\,
\frac{\sqrt{\e}}{\log\frac{1}{\e}}\ ,
\eeqa
and which casts the Hamiltonian $\cH_0$ into $\cH':=\cH_0\circ\Psi$ with
\beq{44bis-}
\cH'(I',\f',p',q'):= h(I')+ g(I',p',q')
+f_*(I',\f',p',q')\,,
\eeq
where (since, by \equ{specham},  $f_{0,0}$ coincides with $\e f_1(I,p,q)$)
\beqa{cons1}
\sup_{(I',p',q')\in U_{r/2}\times V_{\r/2}}
|g_0-\e f_1|&\le& \,{\rm const}\,  \frac{\e^2}{\a_0
r}=\,{\rm const}\,
\frac{\e}{\log(1/\e)}
\ ,\nonumber\\
\| f_*\|_{{r/2,\r/2,s/6}}&\le & \,{\rm const}\,  \e e^{-K s/6}\leq
\,{\rm const}\,\e^N\ . \eeqa Thus, setting $g_0=:\e\overline g$,
$f_*=:\e^N \overline f$,
 we see that $\cH'$ can be rewritten as
\beqa{44bis} \cH'&:=& h(I') + \e \overline g(I',p',q')+\e^N
\overline f
(I',\f',p',q')\,,\nonumber\\
\overline g &=& f_1(I',p',q')+\frac{1}{\log(1/\e)}\overline f_1
(I',p',q')\eeqa
with $\overline f$ and $\overline f_1$  real-analytic on $\cD_1$
(compare \equ{list}).

\giu

We now look for {\sl elliptic equilibria of the Hamiltonian
$\overline g$} in \equ{44bis}. Set
$$ G(I', p',q'):= \Big(
\partial_{p'} \overline{g}(I', p',q'),\,
\partial_{q'} \overline{g}(I', p',q')
\Big)\,.$$
Recalling \equ{44bis} and the definition of $f_1$ in \equ{specham},
we see that, by \equ{eigenst},
$$ G(I',0,0)\Big|_{\e=0}=0\quad{\rm and}\quad \det \partial_{(p',q')}
G(I',0,0)\Big|_{\e=0}= 16 (
\overline\O_1 \overline\O_2)^2>0\, ,\ \forall\ I' \in U_{{r}/2}.$$
Therefore, by the Implicit Function Theorem, we infer that,
for any $I'
\in U_{{r}/2}$ and for $\e$ small enough,
there exist  real-analytic functions, so that
$$
I'\in U_{{r}/2} \to
\Big( p'(I',\e),\, q'(I',\e) \Big)\in B_{{\rm const} / \log\frac{1}{\e}}
\subset B_{\r/2}
\,,$$
and
\begin{equation}\label{bjsdjksdkj1shh}
\partial_{p'} \overline{g}\Big( I', p'(I',\e) ,q'(I',\e)\Big)=0=
\partial_{q'} \overline{g}\Big( I', p'(I',\e) ,q'(I',\e)\Big)\,.
\end{equation}
For $\e$ small enough, we can consider the
following analytic symplectic transformation, which  leaves fixed the $I'$-variable and
is $O(\frac{1}{\log\frac{1}{\e}})-$close to the
identity\footnote{$\Phi'$ has generating function
$J'\cdot\f'+
\Big(v'+p'(J',\e) \Big)\cdot\Big( q'-q'(J',\e)\Big)$.},
$$ \Phi':(J',\psi',v',u')\in
U_{ r/2}\times
\torus^2_{s/7}
\times
D_{\r/3}\,\mapsto\,
(I',\f',p',q')\in U_{ r/2}\times\torus^2_{s/6}
\times
D_{\r/2}^4\,,$$
given by
\beqano
I'&=& J'\,,\\
\f' &=& \psi'+p'(J',\e)\,
\partial_{I'}q'(J',\e)+
\partial_{I'}q'(J',\e)\,v'-
\partial_{I'}p'(J',\e)\,u'\,,\\
p' &=& v'+p'(J',\e)\,,\\
q'&=& u'+q'(J',\e)\,.
\eeqano
In view of \equ{bjsdjksdkj1shh}, the new Hamiltonian
$\hat\cH:=\cH'\circ \Phi'$ has the form
$$ \hat\cH(J',\psi',v',u')=h(J')+\e \hat{g}(J',v',u')+
\e^N\hat{f}(J',\psi',v',u')\,,$$
with $\hat f$ and $\hat g$ analytic in
$U_{ r/2}\times\torus^2_{s/7}
\times
D_{\r/3}$ and
$$ \partial_{v',u'}\hat{g}(J',0,0)=
\partial_{p',q'}\overline{g}(I',p'(I',\e),q'(I',\e))=
0\,,\qquad
\forall\ I'\in U_{r/2}\ .
$$
Also, the eigenvalues of the symplectic quadratic part of
$\hat{g}$ are given by $\pm \ii
\hat{\O}_j (J')$, for $j=1,2$,
where
\beq{eigv*}
\hat{\O}_j\in \real
\qquad{\rm and\qquad}|\hat{\O}_j-\overline{\O}_j|\leq
{\rm const}\frac{1}{\log\frac{1}{\e}}\,.\eeq
Thus,
by a well known result by
Weierstrass on the symplectic diagonalization of qua\-dra\-tic
Hamiltonians, we can find an analytic transformation
$O(\frac{1}{\log\frac{1}{\e}})-$close to the identity
$$ \check{\Phi}:( {J}, {\psi}, {v}, {u})\in
U_{ r/2}\times\torus^2_{s/8}
\times
D_{\r/4}^4\,\mapsto\,
(J',\psi',v',u')\in U_{ r/2}\times\torus^2_{s/7}
\times
D_{\r/3}^4\,,$$
so that $J'= {J}$ and
the quadratic part of $\hat g$ becomes, simply,
$\sum_{j=1}^2 \hat{\O}_j ( {J})\,( {u}_j^2+ {v}_j^2)$.
Whence, the Hamiltonian $\hat\cH$ takes the form
$\check\cH:=\hat\cH\circ\check{\Phi}$,
with
\beq{hamKAMc}
\check\cH( {J}, {\psi}, {v}, {u})= h_\e( {J})+
\e \sum_{j=1}^2\hat{\O}_j
( {J})\,( {u}_i^2+ {v}_i^2)+
\e\widetilde{g} ( {J}, {v}, {u})+
\e^N \widetilde f( {J}, {\psi}, {v}, {u})\,,
\eeq
where
\beq{h0***}
h_\e( {J}):= h_0( {J})+ \e\hat g( {J},0,0)\ ,
\eeq
$\widetilde g$, $\widetilde f$, $\hat{\O}_j$ are real-analytic
for
$( {J}, {\psi}, {v}, {u})$ in
\beq{cD2}
\cD_2:= U_{ r/2}\times\torus^2_{s/8}
\times
D_{\r/4}^4
\eeq
and
\begin{equation}\label{ordine3}
\sup_{ {J}\in U_{ r/2}}|\widetilde g
( {J}, {v}, {u})|\leq \,{\rm const}\,  |( {v}, {u})|^3
\,.\end{equation}
Because of \equ{eigv*},
the non-degeneracy condition \equ{eigenst} implies (for
$\e$ small enough)
\beqa{eigenst**}
&&\inf_{ {J}\in U} \hat{\O}_i \geq
 \inf_{ {J}\in U_{r/2}} |\hat{\O}_i|\geq \,{\rm const}\, >0\ ,\nonumber\\
&& \inf_{ {J}\in U} \Big(\hat{\O}_2-\hat{\O}_1\Big) \geq
\inf_{ {J}\in U_{r/2}}
|\hat{\O}_2-\hat{\O}_1|\geq\,{\rm const}\, >0
\,.
\eeqa

Setting $\widetilde{\O}_j:=2\hat{\O}_j$ for $j=1,2,$
introducing complex coordinates
\beq{coco}
z = \frac{ v + \ii u}{  \sqrt 2}, \qquad
\overline z = \frac{ v - \ii u}{\sqrt 2 }
\eeq
and defining
$$
g( J, z, \overline z) :=
\widetilde g \left(J, \frac{z + \overline z }{ \sqrt 2},
\frac{ z - \overline z }{ \ii \sqrt 2}\right) \quad
{\rm and}  \quad
f(J,\psi,z,\overline z):= \widetilde
f\left(J,\psi, \frac{z+\overline z}{ \sqrt 2},
\frac{ z-\overline z}{ \ii \sqrt 2}\right)\ ,
$$
we obtain \equ{Hstorta}.
Finally \equ{intornini},\equ{stima cubica},\equ{eigenst2},
follows from
\equ{list},\equ{ordine3},\equ{eigenst**}, respectively.

\section*{References}

\footnotesize

\noindent[Am] Ambrosetti, A.: {\sl Critical points and nonlinear
variational problems} Mém. Soc. Math. Fr., Nouv. Sér. 49, 1992.
\smallskip

\noindent[ACE]
Ambrosetti, A.; Coti-Zelati, V., Ekeland, I.:
{\sl Symmetry breaking in Hamiltonian systems},
Journal Diff. Equat. 67, 1987, p. 165-184.
\smallskip

\noindent[AB] Ambrosetti, A.; Badiale, M.:
{\sl Homoclinics: Poincar\'e-Melnikov type results
via a variational approac},
Annales I. H. P. - Analyse nonlin., vol. 15, n.2, 1998, p. 233-252.
\smallskip

\noindent[A]
Arnold, V. I.: {\sl
Small denominators and problems of stability of motion in classical and
celestial mechanics},
Uspehi Mat. Nauk {\bf 18} (1963), no. 6 (114), 91-192.
\smallskip

\noindent[BK]
Bernstein D., Katok, A.: { \sl
Birkhoff periodic orbits for small perturbations
of completely integrable Hamiltonian systems with convex Hamiltonians},
Invent. Math. {\bf 88} (1987), no. 2, 225--241.
\smallskip

\noindent[BBB]
Berti, M., Biasco, L., Bolle, P.: {\sl Drift in phase space,
a new variational mechanism with optimal diffusion time},
to appear on Journal des Mathematiques Pures et Appliques.
\smallskip

\noindent[BB]
Berti, M., Bolle, P. : {\sl A functional analysis approach to
Arnold Diffusion}, Annales de l'I.H.Poincar\'e,
analyse non lineaire, 19, 4, 2002, 395-450.
\smallskip

\noindent[BCV] Biasco, L.; Chierchia, L.; Valdinoci, E.: {\sl
Elliptic two-dimensional invariant tori for the planetary
three-body problem}, to appear in Arch. Ration. Mech.
Anal.\footnote{A preliminary version is
available on-line at {{\tt
http://www.math.utexas.edu/mp\_arc }}}.
\smallskip

\noindent[B]
Birkhoff, G.D.,
{\sl Une generalization \'a $n$-dimensions du dernier th\'eor\`eme
de g\'eometri\'e
de Poincar\'e}, Compt. Rend. Acad. Sci., 192, 1931, 196-198.
\smallskip

\noindent[BL]
Birkhoff, G.D., Lewis, D.C.,
{\sl On the periodic motions near a given periodic
motion of a dynamical system},
Ann. Mat. Pura Appl., IV. Ser. 12, 117-133, 1933.
\smallskip

\noindent[Bo]
Bourgain J.:
{\sl On Melnikov's persistency problem}, Math. Res. Lett. 4, 1997,
445-458.
\smallskip

\noindent[BHS]
Broer H.W., Huitema G.B., Sevriuk M.B., {\sl Quasi periodic
motions in families of dynamical systems}, Lecture
Notes in Math. 1645, Springer, 1996.
\smallskip

\noindent[CZ]
Conley, C.; Zehnder, E.:
{\sl
An index theory for periodic solutions of a Hamiltonian system},
Lecture Notes in Mathematics 1007, Springer, 1983, 132-145.
\smallskip

\noindent[CZ1]
Conley, C.; Zehnder, E.:
{\sl The Birkhoff-Lewis fixed point theorem
and a conjecture of V. I. Arnold},
Invent. Math. 73 (1983), no. 1, 33--49.
\smallskip


\noi[E]
Eliasson, L.:
{\sl Perturbations of stable invariant tori for Hamiltonian systems},
Ann. Scuola Norm. Sup. Pisa, Cl. Sci., {\bf 15} (1988), 115-147.
\smallskip

\noindent[J]
Jefferys, W. H.: {\sl Periodic orbits in the
three-dimensional three-body problem},
Astronom. J. {\bf 71} (1966), no. 7, 566-567.
\smallskip

\noindent[JM] Jefferys, W. H.; Moser, J. : {\sl Quasi-periodic
solutions for the three-body problem}, Astronom. J. {\bf 71}
(1966),  568-578.
\smallskip

\noindent[JV]
Jorba, \'A.; Villanueva J.: {\sl
On the Normal Behaviour of Partially Elliptic Lower Dimensional Tori of Hamiltonian Systems},
Nonlinearity {\bf 10}
(1997), 783-822.
\smallskip

\noi
[K] Kuksin, S. B.:
{\sl Perturbation theory of conditionally periodic solutions of
infinite-dimensional Hamiltonian systems and its applications to the Korteweg-de
Vries equation},
Mat. Sb. (N.S.) {\bf 136 (178)} (1988), no. 3, 396-412,
431; translation in Math. USSR-Sb. {\bf 64} (1989), no. 2, 397-413.
\smallskip

\noindent[LR]
Laskar, J.; Robutel, P.: {\sl Stability of the planetary three-body
problem. I. Expansion of the planetary Hamiltonian}, Celestial Mech.
Dynam. Astronom. {\bf 62} (1995),  no. 3, 193-217.
\smallskip

\noi [L] Lewis, D.C. {\sl Sulle oscillazioni periodiche di un
sistema dinamico}, Atti Acc. Naz. Lincei, Rend. Cl. Sci. Fis. Mat.
Nat., 19, 1934, pp.234-237.
\smallskip

\noi [M] Melnikov, V. K.: {\sl On certain cases of conservation of almost
periodic motions with a small change of the Hamiltonian function},
Dokl. Akad. Nauk SSSR {\bf 165} (1965), 1245-1248.
\smallskip

\noindent[Mo] Moser, J.:{ \sl Proof
of a generalized form of a fixed point Theorem due to G. D. Birkhoff},
Geometry and topology (Proc.
III Latin Amer. School of Math.,
Inst. Mat. Pura Aplicada CNP, Rio de Janeiro, 1976),
pp. 464--494. Lecture Notes in Math.,
Vol. 597, Springer, Berlin, 1977.
\smallskip


\noindent[P] Poincar\'e, H.:
{\sl Les M\'ethodes nouvelles de
la M\'ecanique C\'eleste}, Gauthier 
Villars, Paris,  1892.
\smallskip

\noindent[P\"o] P\"oschel, J.:
{\sl Integrability of Hamiltonian systems on Cantor sets},
Comm. Pure Appl. Math. 35, 1982, 653-696.
\smallskip

\noindent[P\"o1] P\"oschel, J.:
{\sl
On elliptic lower dimensional tori in Hamiltonian system},
Math. Z. {\bf 202} (1989),
559-608.
\smallskip

\noindent[P\"o2] P\"oschel, J.:
{\sl A KAM-Theorem for some nonlinear PDEs},  Ann. Scuola Norm. Sup.
Pisa, Cl. Sci., {\bf 23} (1996), 119-148.
\smallskip

\noindent[R] Robutel, P.:
{\sl Stability of the planetary three-body problem.
II. KAM theory and existence of quasiperiodic motions},
Celestial Mech. Dynam. Astronom. {\bf 62} (1995),  no. 3, 219-261.
\smallskip


\noindent[W]
Wayne, C. E.: {\sl
Periodic and quasi-periodic solutions of nonlinear wave equations
via KAM theory}, Commun. Math. Phys.
{\bf 127} (1990), 479-528.
\smallskip

\noindent[XY]
Xu, J., You, J.: {\sl
Persistence of lower-dimensional tori
under the first Melnikov's non-resonance condition},
J. Math. Pures Appl.  80, 10, 2001, 1045-1067.
\smallskip

\end{document}